\documentclass[11pt,oneside,a4paper]{amsart}

\usepackage{mathptmx}
\usepackage[arrow,curve,matrix,tips,2cell]{xy}

\usepackage[hmargin=2cm,vmargin=1.75cm]{geometry}

\usepackage{fancyhdr}

\usepackage{amsmath}
\usepackage{amscd}
\usepackage{amssymb}
\usepackage{latexsym}
\usepackage{url}

\usepackage{graphicx}

\newtheorem{theorem}{Theorem}[section]
\newtheorem*{theorem*}{Theorem}
\newtheorem{lemma}[theorem]{Lemma}
\newtheorem*{lemma*}{Lemma}
\newtheorem{corollary}[theorem]{Corollary}
\newtheorem{proposition}[theorem]{Proposition}

\newtheorem{remark}[theorem]{Remark}
\newtheorem{definition}[theorem]{Definition}
\newtheorem*{definition*}{Definition}

\newtheorem{question}[theorem]{Question}
\newtheorem*{question*}{Question}

\newtheorem{example}[theorem]{Example}
\newtheorem{examples}[theorem]{Examples}

\makeatletter
\def\revddots{\mathinner{\mkern1mu\raise\p@
\vbox{\kern7\p@\hbox{.}}\mkern2mu
\raise4\p@\hbox{.}\mkern2mu\raise7\p@\hbox{.}\mkern1mu}}
\makeatother 
\newcommand{\bgl}{\begin{equation}} 
\newcommand{\egl}{\end{equation}}
\newcommand{\bgloz}{\begin{equation*}} 
\newcommand{\egloz}{\end{equation*}}
\newcommand{\bgln}{\begin{eqnarray}} 
\newcommand{\egln}{\end{eqnarray}}
\newcommand{\bglnoz}{\begin{eqnarray*}} 
\newcommand{\eglnoz}{\end{eqnarray*}}
\newcommand{\btheo}{\begin{theorem}}
\newcommand{\etheo}{\end{theorem}}
\newcommand{\btheooz}{\begin{theorem*}}
\newcommand{\etheooz}{\end{theorem*}}

\newcommand{\blemma}{\begin{lemma}}
\newcommand{\elemma}{\end{lemma}}
\newcommand{\blemmaoz}{\begin{lemma*}}
\newcommand{\elemmaoz}{\end{lemma*}}
\newcommand{\bproof}{\begin{proof}}
\newcommand{\eproof}{\end{proof}}
\newcommand{\bbew}{\begin{beweis}}
\newcommand{\ebew}{\end{beweis}}
\newcommand{\bremark}{\begin{remark}\em}
\newcommand{\eremark}{\end{remark}}
\newcommand{\bdefin}{\begin{definition}}
\newcommand{\edefin}{\end{definition}}
\newcommand{\bdefinoz}{\begin{definition*}}
\newcommand{\edefinoz}{\end{definition*}}
\newcommand{\bex}{\begin{example}}
\newcommand{\eex}{\end{example}}
\newcommand{\bexs}{\begin{examples}}
\newcommand{\eexs}{\end{examples}}

\newcommand{\bprop}{\begin{proposition}}
\newcommand{\eprop}{\end{proposition}}
\newcommand{\bcor}{\begin{corollary}}
\newcommand{\ecor}{\end{corollary}}
\newcommand{\bfa}{\begin{cases}} 
\newcommand{\efa}{\end{cases}}

\newcommand{\bquestion}{\begin{question}}
\newcommand{\equestion}{\end{question}}
\newcommand{\bquestionoz}{\begin{question*}}
\newcommand{\equestionoz}{\end{question*}}
\newcommand{\cA}{\mathcal A}

\newcommand{\cF}{\mathcal F}
\newcommand{\cG}{\mathcal G}
\newcommand{\cH}{\mathcal H}

\newcommand{\cL}{\mathcal L}

\newcommand{\cS}{\mathcal S}

\def\Cz{\mathbb{C}}

\def\Qz{\mathbb{Q}}
\def\Rz{\mathbb{R}}

\def\Zz{\mathbb{Z}}

\def\1z{\mathbb{1}}
\newcommand{\fL}{\mathfrak L}
\newcommand{\fM}{\mathfrak M}

\newcommand{\fX}{\mathfrak X}
\newcommand{\fY}{\mathfrak Y}
\newcommand{\an}[1]{``#1''} 
\newcommand{\ti}{\tilde}

\newcommand{\lori}{\longrightarrow}

\newcommand{\ma}{\mapsto} 
\newcommand{\onto}{\twoheadrightarrow} 
\newcommand{\into}{\hookrightarrow} 
\newcommand{\Rarr}{\Rightarrow} 
\newcommand{\Larr}{\Leftarrow} 
\newcommand{\LRarr}{\Leftrightarrow} 

\def\SEMI{\mbox{$\times\kern-2pt\vrule height5pt width.6pt \kern3pt $}}

\newcommand{\Hom}{{\rm Hom}\,}

\newcommand{\Spec}{{\rm Spec\,}} 

\newcommand{\id}{{\rm id}}

\newcommand{\Ind}{\mathrm{ Ind}\,}

\newcommand{\img}{{\rm im\,}}
\renewcommand{\ker}{{\rm ker}\,}

\newcommand{\abs}[1]{\lvert#1\rvert} 
\newcommand{\norm}[1]{\left\|#1\right\|} 
\newcommand{\defeq}{\mathrel{:=}} 

\newcommand{\dop}{\text{: }} 
\newcommand{\supp}{{\rm supp}} 
\newcommand{\res}{{\rm res}}
\newcommand{\lge}{\left\{} 
\newcommand{\rge}{\right\}} 
\newcommand{\lru}{\left(} 
\newcommand{\rru}{\right)} 
\newcommand{\lsp}{\left\langle} 
\newcommand{\rsp}{\right\rangle} 
\newcommand{\rukl}[1]{\lru #1 \rru} 
\newcommand{\gekl}[1]{\lge #1 \rge} 
\newcommand{\spkl}[1]{\lsp #1 \rsp} 
\newcommand{\menge}[2]{\gekl{ #1 \dop #2 }} 
\newcommand{\vecg}{\vec{g}}
\newcommand{\vech}{\vec{h}}

\newcommand{\vecgamma}{\vec{\gamma}}

\newcommand{\veceta}{\vec{\eta}}

\newcommand{\climg}{\overline{\rm im}\,}
\begin{document}

\setlength{\parindent}{0cm} \setlength{\parskip}{0.5cm}

\title{Dynamic characterizations of quasi-isometry, \\ and applications to cohomology}

\thispagestyle{fancy}

\author{Xin Li}

\address{Xin Li, School of Mathematical Sciences, Queen Mary University of London, Mile End Road, London E1 4NS, United Kingdom}
\email{xin.li@qmul.ac.uk}


\begin{abstract}
We build a bridge between geometric group theory and topological dynamical systems by establishing a dictionary between coarse equivalence and continuous orbit equivalence. As an application, we give conceptual explanations for previous results of Shalom and Sauer on coarse invariance of homological and cohomological dimensions and Shalom's property $H_{FD}$. As another application, we show that group homology and cohomology in a class of coefficients, including all induced and co-induced modules, are coarse invariants. We deduce that being of type $FP_n$ (over arbitrary rings) is a coarse invariant, and that being a (Poincar{\'e}) duality group over a ring is a coarse invariant among all groups which have finite cohomological dimension over that ring. Our results also imply that every self coarse embedding of a Poincar{\'e} duality group over an arbitrary ring must be a coarse equivalence.
\end{abstract}

\thanks{Research supported by EPSRC grant EP/M009718/1.}

\maketitle

\section{Introduction}
\label{intro}

The philosophy of geometric group theory is to study groups not merely as algebraic objects but from a geometric point of view. There are two ways of developing a geometric perspective, by viewing groups themselves as geometric objects (for instance with the help of their Cayley graphs, which leads to the notion of quasi-isometry) or by studying groups by means of \an{nice} group actions on spaces which carry some topology or geometry. Once a geometric point of view is taken, an immediate question is: How much of the original algebraic structures is still visible from our new perspective? Or more precisely: Which algebraic invariants of groups are quasi-isometry invariants?

Our goals in this paper are twofold. First, we want to connect the two geometric perspectives mentioned above by giving dynamic characterizations of quasi-isometry, or more generally, coarse equivalence. It turns out that for topological dynamical systems, the concept corresponding to coarse equivalence is given by (modified versions of) continuous orbit equivalence, as introduced in \cite{Li1,Li2}. The latter means that we can identify the orbit structure of our dynamical systems in a continuous way. The idea of developing dynamic characterizations of coarse equivalence goes back to Gromov's notion of topological couplings and has been developed further in \cite{Sha,Sau}. Recently, independently from the author, a dynamic characterization of bilipschitz equivalence for finitely generated groups was obtained in \cite{MST}, which is a special case of our result.
\setlength{\parindent}{0.5cm} \setlength{\parskip}{0cm}

Secondly, we want to study the behaviour of algebraic invariants of groups under coarse equivalence. More precisely, we consider invariants of (co)homological nature. Using our dynamic characterizations of coarse equivalence, we give conceptual explanations of the results in \cite{Sha,Sau} on coarse invariance of homological and cohomological dimensions and Shalom's property $H_{FD}$. Moreover, using a refined, more concrete version of our dynamic characterizations, we produce many new coarse invariants of (co)homological nature. We generalize the result in \cite{Ger} that among groups $G$ satisfying the finiteness condition $F_n$ (i.e., there exist models for Eilenberg-MacLane spaces with finite $n$-skeleton), the cohomology groups $H^n(G,RG)$ are coarse invariants for all commutative rings $R$ with unit. We show that for a class of coefficients (called $\res$-invariant modules), including all induced and co-induced modules, group homology and cohomology are coarse invariants. In particular, $H^*(G,RG)$ is always a coarse invariant. This answers a question in \cite{MSW} (see \cite[Questions after Theorem~2.7]{MSW}). Our results imply that being of type $FP_n$ over $R$ (i.e., the trivial $RG$-module $R$ admits a projective resolution which is finitely generated up to level $n$) is a coarse invariant. This is a partial generalization of \cite[Theorem~1.7]{Sha}. A different approach is mentioned in \cite[Theorem~9.61]{DK}, and the case $R = \Zz$ has been treated in \cite{Alo}. As a consequence, we obtain that for an arbitrary commutative ring $R$ with unit, the property of being a duality or Poincar{\'e} duality group over $R$ is a coarse invariant among all groups which have finite cohomological dimension over $R$. A group $G$ is called a duality group over $R$ if there is a right $RG$-module $C$ and an integer $n \geq 0$ with natural isomorphisms $H^k(G,A) \cong H_{n-k}(G,C \otimes_R A)$ for all $k \in \Zz$ and all $RG$-modules $A$ (see \cite[\S~9.2]{Bie}, \cite{Bie75}, and \cite[Chapter~VIII, \S~10]{Bro}). $G$ is called a Poincar{\'e} duality group over $R$ if $C \cong R$ as $R$-modules. $C$ is called the dualizing module; note that we must have $C \cong H^n(R,RG)$ as right $RG$-mdules. Our result generalizes \cite[Corollary~3]{Ger}, as we do not need the finiteness condition $F_{\infty}$ (i.e., $F_n$ for all $n$) and can work over arbitrary rings. Examples of groups which are not duality groups over $\Zz$ but over some other ring can be found in \cite{Dav}, and examples of (Poincar{\'e}) duality groups which are not of type $F_{\infty}$ appear in \cite{Dav,Lea}. Combined with Sauer's result \cite[Theorem~1.2~(ii)]{Sau}, we obtain that among amenable groups, being a (Poincar{\'e}) duality group over a divisible ring is a coarse invariant. This generalizes \cite[Theorem~3.3.2]{Sau}. We also prove a rigidity result for coarse embeddings into Poincar{\'e} duality groups. If a group $G$ with ${\rm hd}_R \, G < \infty$ coarsely embeds into a Poincar{\'e} duality group $H$, then ${\rm hd}_R \, G < {\rm cd}_R \, H$. In particular, self coarse embeddings of Poincar{\'e} duality groups over an arbitrary ring must be coarse equivalences.
\setlength{\parindent}{0cm} \setlength{\parskip}{0.5cm}

Let us now formulate and explain our main results in more detail. At the same time, we fix some notations. Throughout this paper, all our groups are countable and discrete. First, we recall the notion of coarse maps (see \cite[Definition~2.21]{Roe}). Note that coarse embeddings in our sense are called uniform embeddings in \cite{Sha,Sau}.
\setlength{\parindent}{0cm} \setlength{\parskip}{0cm}

\bdefin
\label{Def:CM}
A map $\varphi: \: G \to H$ between two groups $G$ and $H$ is called a coarse map if $\varphi^{-1}(\gekl{y})$ is finite for all $y \in H$, and for every $S \subseteq G \times G$ with $\menge{st^{-1}}{(s,t) \in S}$ finite, $\menge{\varphi(s) \varphi(t)^{-1}}{(s,t) \in S}$ is finite. 
\setlength{\parindent}{0.5cm} \setlength{\parskip}{0cm}

$\varphi: \: G \to H$ is called a coarse embedding if for every subset $S \subseteq G \times G$, $\menge{st^{-1}}{(s,t) \in S}$ is finite if and only if $\menge{\varphi(s) \varphi(t)^{-1}}{(s,t) \in S}$ is finite. 

Two maps $\varphi, \, \phi: \: G \to H$ are called close if $\menge{\varphi(x) \phi(x)^{-1}}{x \in G}$ is finite. We write $\varphi \sim \phi$ in that case.

A coarse map $\varphi: \: G \to H$ is called a coarse equivalence if it is coarsely invertible, i.e., there is a coarse map $\psi: \: H \to G$ such that $\psi \circ \varphi \sim \id_G$ and $\varphi \circ \psi \sim \id_H$.

We say that two groups $G$ and $H$ are coarsely equivalent if there is a coarse equivalence $G \to H$.
\edefin
\setlength{\parindent}{0cm} \setlength{\parskip}{0cm}

Clearly, coarse embeddings are coarse maps. Examples of coarse embeddings are subgroup embeddings and quasi-isometric embeddings. For finitely generated groups, coarse equivalences coincide with quasi-isometries (see \cite{Sha}). Note that unlike in \cite{Sha,Sau}, in our definition, we use $st^{-1}$ and not $s^{-1}t$ (see Remark~\ref{Rem:left--right}).
\setlength{\parindent}{0cm} \setlength{\parskip}{0.5cm}

Let us explain our dynamic characterizations of coarse embeddings and equivalences. Let $G \curvearrowright X$ and $H \curvearrowright Y$ be topological dynamical systems, where the groups act by homeomorphisms on locally compact Hausdorff spaces. A continuous orbit couple is a pair of continuous maps $p: \: X \to Y$ and $q: \: Y \to X$ which both preserve orbits in a continuous way, such that $p$ and $q$ are inverses up to orbits (i.e., $q(p(x))$ lies in the same $G$-orbit of $x$ and similarly for $p \circ q$). \an{Preserving orbits in a continuous way} is made precise by continuous maps $a: \: G \times X \to H$ such that $p(g.x) = a(g,x).p(x)$ for all $g \in G$ and $x \in X$. If $p$ and $q$ are actual inverses (i.e., $q \circ p = \id_X$ and $p \circ q = \id_Y$), then our dynamical systems are called continuously orbit equivalent.
\setlength{\parindent}{0.5cm} \setlength{\parskip}{0cm}

Our first main result establishes the following dictionary: The existence of a coarse embedding $G \to H$ corresponds to the existence of a continuous orbit couple for topologically free systems $G \curvearrowright X$ and $H \curvearrowright Y$, where $X$ is compact. The existence of a coarse equivalence $G \to H$ corresponds to the existence of a continuous orbit couple for topologically free systems $G \curvearrowright X$ and $H \curvearrowright Y$, where both $X$ and $Y$ are compact, and we can find a bijective coarse equivalence $G \to H$ if and only if we can find a continuously orbit equivalence for $G \curvearrowright X$ and $H \curvearrowright Y$. We refer to Theorem~\ref{u,q,l--DS} for precise statements.

It turns out that for compact $X$, the existence of a continuous orbit couple for $G \curvearrowright X$ and $H \curvearrowright Y$ is equivalent to saying that $G \curvearrowright X$ and $H \curvearrowright Y$ are Kakutani equivalent, i.e., there are clopen subspaces $A \subseteq X$ and $B \subseteq Y$ which are $G$- and $H$-full such that the partial actions $G \curvearrowright A$ and $H \curvearrowright B$ are continuously orbit equivalent (in the sense of \cite{Li2}). This implies that the transformation groupoids of $G \curvearrowright X$ and $H \curvearrowright Y$ are Morita equivalent. Building on this observation, we show that the results in \cite{Sha,Sau} on coarse invariance of (co)homological dimension and Shalom's property $H_{FD}$ are immediate consequences of Morita invariance of groupoid (co)homology. This gives a conceptual explanation for the results in \cite{Sha,Sau}, and at the same time, our work isolates precise conditions on the dynamical systems which are needed to show coarse invariance.
\setlength{\parindent}{0cm} \setlength{\parskip}{0.5cm}

The dynamic characterizations we described so far are abstract as the dynamical systems are not specified. It is striking that even such abstract characterizations suffice to derive the results in \cite{Sha,Sau}. However, to show coarse invariance of group (co)homology in particular coefficients, we need more concrete versions of our dynamic characterizations. Inspired by \cite{SW}, we first observe that in place of abstract dynamical systems, we may always take the canonical action $G \curvearrowright \beta G$ of groups $G$ on their Stone-{\v C}ech compactifications $\beta G$. The appearance of $G \curvearrowright \beta G$ is not surprising because of its universal property. But now, our crucial observation is that we can go even further and consider the actions $G \curvearrowright G$ of groups acting on themselves by left multiplication. By doing so, it seems that we are losing all the information as any two actions $G \curvearrowright G$ and $H \curvearrowright H$ are continuously orbit equivalent as long as $G$ and $H$ have the same cardinality. The problem is that the spaces on which our groups act are no longer compact. However, we can replace compactness by asking for finiteness conditions on the maps $a$, which -- as in the definition of continuous orbit couples -- make precise that orbits are preserved in a continuous way: We require that for every $g \in G$, the map $a(g,\cdot)$ should have finite image. It is this finiteness condition which singles out \an{controlled} orbit equivalences which behave well in (co)homology. The point is that every coarse equivalence $G \to H$ gives rise to a \an{controlled} orbit equivalence between $G \curvearrowright G$ and $H \curvearrowright H$. This change of perspective, putting the emphasis on this finiteness condition, turns out to be crucial.

These ideas lead to the following results: Let $R$ be a commutative ring with unit and $W$ an $R$-module. The set $C(G,W)$ of functions $G \to W$ carries a natural $RG$-module structure. An $RG$-submodule $L \subseteq C(G,W)$ is called $\res$-invariant if for every $f \in L$ and $A \subseteq G$, the restriction of $f$ to $A$ (viewed as a function on $G$ by extending it by $0$) still lies in $L$. Examples include $C(G,W)$, the submodule $C_f(G,W)$ of $f \in C(G,W)$ taking only finitely many values, $RG \otimes_R W$, and for $W = R = \Rz$ or $\Cz$, $c_0(G,W) = \menge{f: \: G \to W}{\lim_{x \to \infty} \abs{f(x)} = 0}$, $\ell^p(G,W) = \menge{f: \: G \to W}{\sum_{x \in G} \abs{f(x)}^p < \infty}$ ($0 < p \leq \infty$), $H^{s,p}(G,W) = \menge{f: \: G \to W}{f \cdot (1 + \ell)^s \in \ell^p (G,W)}$ ($s \in \Rz \cup \gekl{\infty}$, $1 \leq p < \infty$), where $G$ is finitely generated and $\ell$ is the word length on $G$, and $H^{\infty,p}(G,W) = \bigcap_{s \in \Rz} H^{s,p}(G,W)$. We show that a coarse equivalence $\varphi: \: G \to H$ induces a one-to-one correspondence between $\res$-invariant submodules of $C(G,W)$ and $\res$-invariant submodules of $C(H,W)$, denoted by $L \ma \varphi_* L$, together with isomorphisms $H_*(\varphi): \: H_*(G,L) \cong H_*(H,\varphi_* L)$ for all $L$. Similarly, $\varphi$ induces a one-to-one correspondence between $\res$-invariant submodules of $C(H,W)$ and $\res$-invariant submodules of $C(G,W)$, say $M \ma \varphi^* M$, together with isomorphisms $H^*(\varphi): \: H^*(H,M) \cong H^*(G,\varphi^* M)$ for all $M$. In particular, we obtain
\setlength{\parindent}{0cm} \setlength{\parskip}{0cm}

\btheooz[Corollary~\ref{RGlpc0}]
Among all countable discrete groups $G$, the following (co)homology groups are coarse invariants: $H_*(G,C(G,W))$, $H_*(G,C_f(G,W))$, $H^*(G,C_f(G,W))$, $H^*(G,RG \otimes_R W)$ for every commutative ring $R$ with unit and every $R$-module W;
\setlength{\parindent}{1cm} \setlength{\parskip}{0cm}

$H_*(G,c_0 (G,R))$, $H^*(G,c_0 (G,R))$, ${\bar H}_*(G,c_0 (G,R))$, ${\bar H}^*(G,c_0 (G,R))$; $H_*(G,\ell^p (G,R))$, $H^*(G,\ell^p (G,R))$, ${\bar H}_*(G,\ell^p (G,R))$, ${\bar H}^*(G,\ell^p (G,R))$, for all $0<p \leq \infty$; and for finitely generated groups $G$, $H_*(G,H^{s,p}(G,R))$, $H^*(G,H^{s,p}(G,R))$, ${\bar H}_*(G,H^{s,p}(G,R))$, ${\bar H}^*(G,H^{s,p}(G,R))$, for all $s \in \Rz \cup \gekl{\infty}$, $1 \leq p \leq \infty$, where $R = \Rz$ or $\Cz$.
\etheooz
\setlength{\parindent}{0cm} \setlength{\parskip}{0cm}

Some of these (co)homology groups can be identified with existing (co)homology theories: $H^*(G,RG)$ is coarse cohomology \cite[\S~5.1]{Roe}, $H_*(G,C_f(G,\Zz))$ and $H_*(G,\ell^{\infty}(G,\Rz))$ coincide with uniformly finite homology \cite{BW, BNW, BD}, and for $\ell^p$ coefficients, we obtain $L^p$-cohomology \cite{Pan,Ger,El}. Actually, we show that every coarse map $\varphi: \: G \to H$ induces a map $H_*(\varphi): \: H_*(G,L) \to H_*(H,\varphi_* L)$ such that $H_*(\varphi) = H_*(\phi)$ if $\varphi \sim \phi$ and $H_*(\psi \circ \varphi) = H_*(\psi) \circ H_*(\varphi)$. It is then evident that coarse equivalences induce isomorphisms as they are precisely those coarse maps which are invertible modulo $\sim$. A similar remark applies to cohomology. Thus, not only these (co)homology groups, but, by functoriality, the actions of the groups of coarse equivalences (modulo $\sim$) on these (co)homology groups are coarse invariants as well. We obtain analogous results for coarse embeddings in the topological setting, i.e., for topological $\res$-invariant modules and reduced (co)homology. It turns out that coarse embeddings always induce isomorphisms in (co)homology and reduced (co)homology.
\setlength{\parindent}{0.5cm} \setlength{\parskip}{0cm}

The aforementioned results on coarse invariance of type $FP_n$ and being a (Poincar{\'e}) duality group are immediate consequences, as is our rigidity result for coarse embeddings into Poincar{\'e} duality groups. We also deduce that vanishing of $\ell^2$-Betti numbers is a coarse invariant, as observed in \cite{Pan,Ogu,MOSS}, and generalized by Sauer and Schr{\"o}dl to all unimodular locally compact second countable groups \cite{SauSch}.

This is a good point to formulate an interesting and natural question, which we elaborate on in \S~\ref{ss:Consequences}:
\bquestionoz[Question~\ref{Q:dim}]
Are homological and cohomological dimension over a commutative ring $R$ with unit always coarse invariants among all countable discrete groups with no $R$-torsion?
\equestionoz
\setlength{\parindent}{0cm} \setlength{\parskip}{0cm}

We refer to \S~\ref{sec:HII} for more details. \S~\ref{sec:HI} and \S~\ref{sec:HII} are independent from each other. Thus readers interested in this last set of results on coarse invariance of group (co)homology may go directly from \S~\ref{sec:DynChar} to \S~\ref{sec:HII}.
\setlength{\parindent}{0cm} \setlength{\parskip}{0.5cm}

As far as our methods are concerned, we use groupoid techniques as in \cite{Sha,Sau,Ogu}, but instead of working with abstract dynamical systems, we base our work on concrete dynamic characterizations of coarse equivalence. The difference between our work and \cite{Ger} is that we do not work with descriptions of group (co)homology in terms of Eilenberg-MacLane spaces, as these descriptions require finiteness conditions (like $F_n$ or $F_{\infty}$) on our groups and have to be modified whenever we change coefficients. Instead, since coarse embeddings automatically lead to \an{controlled} orbit equivalences satisfying the finiteness condition mentioned above, we can work directly with complexes coming from bar resolutions.

I thank Andreas Thom for informing me about \cite{MST}, Roman Sauer for explaining certain parts of his work \cite{Sau} to me, and for informing me about \cite{SauSch,MOSS}, and Roman Sauer and Clara L{\"o}h for very helpful comments.

\section{Dynamical characterizations of quasi-isometry}
\label{sec:DynChar}

\subsection{Preliminaries}

The central notions of coarse maps, embeddings and equivalences have been introduced in \S~\ref{intro}. We remark that it is easy to see that a coarse embedding $\varphi: \: G \to H$ is coarsely invertible if and only if $H$ can be covered by finitely many translates of $\varphi(G)$, i.e., there is a finite set $F \subseteq H$ such that $H = \bigcup_{h \in F} h \varphi(G)$.
\setlength{\parindent}{0.5cm} \setlength{\parskip}{0cm}

\bremark
\label{Rem:left--right}
Note that unlike in \cite{Sha}, our definition of coarse maps is right-invariant, not left-invariant (i.e., we use $st^{-1}$ instead of $s^{-1} t$). For finitely generated groups, this amounts to considering right-invariant word lengths and word metrics. We do so because in the following, we will consider left actions of groups, in particular the action of a group by left multiplication on itself. Of course, this is merely a matter of convention.
\eremark
\setlength{\parindent}{0cm} \setlength{\parskip}{0cm}

The following concept, due to Gromov, builds a bridge between geometric group theory and topological dynamical systems.
\bdefin
For two groups $G$ and $H$, a $(G,H)$ topological coupling consists of a locally compact space $\Omega$ with commuting free and proper left $G$- and right $H$-actions which admit clopen $H$- and $G$-fundamental domains $\bar X$ and $\bar Y$. Our $(G,H)$ topological coupling is called $G$-cocompact if $\bar Y$ is compact, $H$-cocompact if $\bar X$ is compact, and cocompact if it is both $G$- and $H$-cocompact. It is called topologically free (or free) if the combined action $G \times H \curvearrowright \Omega$ is topologically free (or free).
\edefin
\setlength{\parindent}{0cm} \setlength{\parskip}{0cm}
All our spaces are Hausdorff. Also, being only concerned with the topological setting, we simply write \an{coupling} (without prefix \an{topological}). We often write $G \curvearrowright {}_{\bar{Y}}\Omega {}_{\bar X} \curvearrowleft H$ to keep track of all the relevant data.
\setlength{\parindent}{0cm} \setlength{\parskip}{0.5cm}

The following result goes back to ideas of Gromov and is proven in \cite{Sha} and \cite{Sau}.
\setlength{\parindent}{0cm} \setlength{\parskip}{0cm}
\btheo
\label{THM:Gromov}
Let $G$ and $H$ be countable discrete groups.
\begin{enumerate}
\item[(i)] There exists a coarse embedding $G \to H$ if and only if there exists a $H$-cocompact $(G,H)$ coupling.
\item[(ii)] There exists a coarse equivalence $G \to H$ if and only if there exists a cocompact $(G,H)$ coupling.
\item[(iii)] There is a bijective coarse equivalence $G \to H$ if and only if there is a cocompact $(G,H)$ coupling $G \curvearrowright {}_{\bar{Y}}\Omega {}_{\bar X} \curvearrowleft H$ with $\bar X = \bar Y$.
\end{enumerate}
\etheo
\bproof
For (i), see \cite[Theorem~2.2, (i) $\LRarr$ (ii)]{Sau}. For (ii), see \cite[Theorem~2.2, (iii) $\LRarr$ (iv)]{Sau}. For (iii), see \cite[Remark after Theorem~2.1.2]{Sha}.
\eproof
\setlength{\parindent}{0cm} \setlength{\parskip}{0.5cm}

\bremark
\label{TOTDISCONN}
The proofs in \cite{Sau,Sha} show that the underlying space $\Omega$ of the $(G,H)$ couplings can be chosen to be second countable and totally disconnected in the above statements.
\eremark
\setlength{\parindent}{0cm} \setlength{\parskip}{0cm}

Let us now isolate an idea from \cite{MST} which will be useful later on.
\blemma
\label{Lem:Coupling-topfree}
If there exists a $(G,H)$ coupling $G \curvearrowright \Omega \curvearrowleft H$, then there exists a topologically free $(G,H)$ coupling $G \curvearrowright \Omega' \curvearrowleft H$. If $G \curvearrowright \Omega \curvearrowleft H$ is $G$-cocompact, $H$-cocompact or cocompact, $G \curvearrowright \Omega' \curvearrowleft H$ may be chosen with the same property. If $\Omega$ is second countable and totally disconnected, we may choose $\Omega'$ with the same property.
\elemma
\bproof
The idea of the proof appears in the proof of \cite[Theorem~3.2]{MST}. Let $G \times H \curvearrowright Z$ be a free action on the Cantor space $Z$. It is easy to see that $\Omega' = \Omega \times Z$ with diagonal $G$- and $H$-actions is a $(G,H)$ coupling which is topologically free (even free). As $Z$ is compact and totally disconnected, our additional claims follow.
\eproof
\setlength{\parindent}{0cm} \setlength{\parskip}{0.5cm}

\subsection{Topological couplings and continuous orbit couples}
\label{sec_TC-COC}

We explain the connection between topological couplings and continuous orbit couples. First of all, a topological dynamical system $G \curvearrowright X$ consists of a group $G$ acting on a locally compact space $X$ via homeomorphisms. We write $g.x$ for the action.
\bdefin
\label{coc}
Let $G \curvearrowright X$ and $H \curvearrowright Y$ be topological dynamical systems.
\setlength{\parindent}{0.5cm} \setlength{\parskip}{0cm}

A continuous map $p: \: X \to Y$ is called a continuous orbit map if there exists a continuous map $a: \: G \times X \to H$ such that $p(g.x) = a(g,x).p(x)$ for all $g \in G$ and $x \in X$.

A continuous orbit couple for $G \curvearrowright X$ and $H \curvearrowright Y$ consists of continuous orbit maps $p: \: X \to Y$ and $q: \: Y \to X$ such that there exist continuous maps $g: \: X \to G$ and $h: \: Y \to H$ such that $q(p(x)) = g(x).x$ and $p(q(y)) = h(y).y$ for all $x \in X$ and $y \in Y$.
\setlength{\parindent}{0cm} \setlength{\parskip}{0.5cm}
\edefin

\bdefin
A $(G,H)$ continuous orbit couple consists of topological dynamical systems $G \curvearrowright X$ and $H \curvearrowright Y$ and a continuous orbit couple for $G \curvearrowright X$ and $H \curvearrowright Y$. If $G \curvearrowright X$ and $H \curvearrowright Y$ are topologically free, then the $(G,H)$ continuous orbit couple is called topologically free. We call $X$ the $G$-space and $Y$ the $H$-space of our $(G,H)$ continuous orbit couple.
\edefin

\bremark
\label{X=Y-->coc=coe}
In this language, a continuous orbit equivalence for $G \curvearrowright X$ and $H \curvearrowright Y$ in the sense of \cite{Li1} is the same as a continuous orbit couple for $G \curvearrowright X$ and $H \curvearrowright Y$ with $g \equiv e$ and $h \equiv e$, i.e., $p = q^{-1}$.
\eremark

\bdefin
A $(G,H)$ continuous orbit equivalence consists of topological dynamical systems $G \curvearrowright X$ and $H \curvearrowright Y$ and a continuous orbit equivalence for $G \curvearrowright X$ and $H \curvearrowright Y$.
\edefin
\setlength{\parindent}{0cm} \setlength{\parskip}{0cm}

\btheo
\label{1-1}
Let $G$ and $H$ be groups. There is a one-to-one correspondence between isomorphism classes of topologically free $(G,H)$ couplings and isomorphism classes of topologically free $(G,H)$ continuous orbit couples, with the following additional properties:
\begin{enumerate}
\item[(i)] A $(G,H)$ coupling $G \curvearrowright {}_{\bar{Y}}\Omega {}_{\bar X} \curvearrowleft H$ corresponds to a $(G,H)$ continuous orbit couple with $G$-space homeomorphic to $\bar X$ and $H$-space homeomorphic to $\bar Y$.
\item[(ii)] A $(G,H)$ coupling $G \curvearrowright {}_{\bar{Y}}\Omega {}_{\bar X} \curvearrowleft H$ with $\bar X = \bar Y$ corresponds to a $(G,H)$ continuous orbit equivalence.
\end{enumerate}
\etheo
Here, the notions of isomorphisms are the obvious ones: Topological couplings $G \curvearrowright {}_{\bar{Y}_1} \Omega_1 {}_{\bar{X}_1} \curvearrowleft H$ and $G \curvearrowright {}_{\bar{Y}_2} \Omega_2 {}_{\bar{X}_2} \curvearrowleft H$ are isomorphic if there exists a $G \times H$-equivariant homeomorphism $\Omega_1 \cong \Omega_2$ sending $\bar{X}_1$ to $\bar{X}_2$ and $\bar{Y}_1$ to $\bar{Y}_2$. Continuous orbit couples $(p_i,q_i)$ for $G \curvearrowright X_i$ and $H \curvearrowright Y_i$, $i=1,2$, are isomorphic if there exist $G$- and $H$-equivariant homeomorphisms $X_1 \cong X_2$ and $Y_1 \cong Y_2$ such that we obtain commutative diagrams
$$
  \xymatrix{
  X_1 \ar[d]_{\cong} \ar[r]^{p_1} & Y_1 \ar[d]^{\cong}
  \\
  X_2 \ar[r]^{p_2} & Y_2
  }
\ \ \ \ \ \ 
  \xymatrix{
  Y_1 \ar[d]_{\cong} \ar[r]^{q_1} & X_1 \ar[d]^{\cong}
  \\
  Y_2 \ar[r]^{q_2} & X_2
  }
$$
For the proof of Theorem~\ref{1-1}, we now present explicit constructions of continuous orbit couples out of topological couplings and vice versa. The constructions are really the topological analogues of those in \cite[\S~3]{Fur} (see also \cite{Sha,Sau}). In the following, we write $gx$ ($g \in G, x \in \Omega$) and $xh$ ($x \in \Omega, h \in H$) for the left $G$- and right $H$-actions in topological couplings, and $g.x$, $h.y$ for the actions $G \curvearrowright X$, $H \curvearrowright Y$ from continuous orbit couples.

\subsubsection{From topological couplings to continuous orbit couples}
\label{tc->coc}

Let $G \curvearrowright {}_{\bar{Y}}\Omega {}_{\bar X} \curvearrowleft H$ be a $(G,H)$ coupling. Set $X \defeq \bar{X}$ and $Y \defeq \bar{Y}$. Define a map $p: \: X \to Y$ by requiring $Gx \cap Y = \gekl{p(x)}$ for all $x \in X$. The intersection $Gx \cap Y$, taken in $\Omega$, consists of exactly one point because $Y$ is a $G$-fundamental domain. By construction, there is a map $\gamma: \: X \to G$ such that $p(x) = \gamma(x) x$. For $g \in G$, $\gamma$ takes the constant value $g$ on $X \cap g^{-1}Y$. As $X \cap g^{-1}Y$ is clopen, because $X$ and $Y$ are, $\gamma$ is continuous. $p$ is continuous as it is so on $X \cap g^{-1}Y$ for all $g \in G$.
\setlength{\parindent}{0cm} \setlength{\parskip}{0.5cm}

We now define a $G$-action, denoted by $G \times X \to X, \, (g,x) \ma g.x$, as follows: For every $g \in G$ and $x \in X$, there exists a unique $\alpha(g,x) \in H$ such that $gx \in X \alpha(g,x)$. For fixed $g \in G$ and $h \in H$, we have $\alpha(g,x) = h$ for all $x \in X \cap g^{-1}Xh$. As $X \cap g^{-1}Xh$ is clopen because $X$ is, $\alpha: \: G \times X \to H$ is continuous. Set $g.x \defeq gx\alpha(g,x)^{-1}$. It is easy to check that $\alpha$ satisfies the cocycle identity $\alpha(g_1 g_2,x) = \alpha(g_1,g_2.x) \alpha(g_2,x)$. Using this, it is easy to see that $G \times X \to X, \, (g,x) \ma g.x$ defines a (left) $G$-action on $X$ by homeomorphisms.

Similarly, we define a continuous map $q: \: Y \to X$ by requiring $X \cap yH = \gekl{q(y)}$ for all $y \in Y$, and let $\eta: \: Y \to H$ be the continuous map satisfying $q(y) = y\eta(y)$. To define an $H$-action on $Y$, let $\beta(y,h) \in G$ be such that $yh \in \beta(y,h) Y$. Again, $\beta: \: Y \times H \to G$ is continuous. Set $h.y \defeq \beta(y,h^{-1})^{-1}yh^{-1}$. It is easy to check that $\beta$ satisfies $\beta(y,h_1 h_2) = \beta(y,h_1) \beta(h_1^{-1}.x h_2)$. Using this, it is again easy to see that $H \times Y \to Y, \, (h,y) \ma h.y$ defines an $H$-action on $Y$ by homeomorphisms.

Let us check that $(p,q)$ is a $(G,H)$ continuous orbit couple. We need to identify $G g x \alpha(g,x)^{-1} \cap Y$ in order to determine $p(g.x) = p(gx\alpha(g,x)^{-1})$. We have
$$G g x \alpha(g,x)^{-1} \ni \beta(\gamma(x)x,\alpha(g,x)^{-1})^{-1} \gamma(x) x \alpha(g,x)^{-1} \in Y,$$
so $p(g.x) = \beta(\gamma(x)x,\alpha(g,x)^{-1})^{-1} \gamma(x) x \alpha(g,x)^{-1} = \alpha(g,x).(\gamma(x)x) = \alpha(g,x).p(x)$. Similarly, in order to identify $q(h.y) = q(\beta(y,h^{-1})^{-1}yh^{-1})$, we need to determine $X \cap \beta(y,h^{-1})^{-1}yh^{-1}H$. As
$$X \ni \beta(y,h^{-1})^{-1} y \eta(y) \alpha(\beta(y,h^{-1})^{-1},y \eta(y))^{-1} \in \beta(y,h^{-1})^{-1}yh^{-1}H,$$
we deduce $q(y.h) = \beta(y,h^{-1})^{-1} y \eta(y) \alpha(\beta(y,h^{-1})^{-1},y \eta(y))^{-1} = \beta(y,h^{-1})^{-1}.(y\eta(y)) = \beta(y,h^{-1})^{-1}.q(y)$. Finally, $qp(x) = q(\gamma(x)x) = \gamma(x)x\alpha(\gamma(x),x)^{-1} = \gamma(x).x$ and $pq(y) = p(y\eta(y)) = \beta(y,\eta(y))^{-1} y \eta(y) = \eta(y)^{-1}.y$. All in all, we see that $p$ and $q$ give rise to a continuous orbit couple for $G \curvearrowright X$ and $H \curvearrowright Y$, with $g(x) = \gamma(x)$ and $h(y) = \eta(y)^{-1}$.

Note that our coupling does not need to be topologically free for this construction. However, it is clear that $G \curvearrowright \Omega \curvearrowleft H$ is topologically free (i.e., $G \times H \curvearrowright \Omega$ is topologically free) if and only if $G \curvearrowright X$ and $H \curvearrowright Y$ are topologically free.

\bremark
Our notation differs slightly from the one in \cite{Sha} and \cite{Sau}. Our $\alpha(g,x)$ is $\alpha(g^{-1},x)^{-1}$ in \cite[\S~2.2, Equation~(3)]{Sha} and \cite[\S~2.2, Equation~(2.2)]{Sau}. This is closely related to Remark~\ref{Rem:left--right}.
\eremark

\bremark
The dynamical system $G \curvearrowright X$ we constructed above can be canonically identified with $G \curvearrowright \Omega / H$. Similarly, our system $H \curvearrowright Y$ can be identified with $G \backslash \Omega \curvearrowleft H$ in a canonical way.
\eremark
\setlength{\parindent}{0cm} \setlength{\parskip}{0cm}

\subsubsection{From continuous orbit couples to topological couplings}
\label{coc->tc}

Let $G \curvearrowright X$ and $H \curvearrowright Y$ be topologically free systems on locally compact spaces $X$ and $Y$. Assume that we are given a continuous orbit couple for $G \curvearrowright X$ and $H \curvearrowright Y$, and let $p$, $q$, $a$, $g$ and $h$ be as in Definition~\ref{coc}, and let $b: \: H \times Y \to G$ be a continuous map with $q(h.y) = b(h,y).q(y)$ for all $h \in H$ and $y \in Y$. Define commuting left $G$- and right $H$-actions on $X \times H$ by $g(x,h) = (g.x,a(g,x)h)$, $(x,h)h' = (x,hh')$. Furthermore, define commuting left $G$- and right $H$-actions on $G \times Y$ by $g'(g,y) = (g'g,y)$ and $(g,y)h = (gb(h^{-1},y)^{-1},h^{-1}.y)$.
\setlength{\parindent}{0cm} \setlength{\parskip}{0.5cm}

A straightforward computation, using the cocycle identities (\cite[Lemma~2.8]{Li1}) for $a$ and $b$, shows that $\Theta: \: X \times H \to G \times Y, \, (x,h) \ma (g(x)^{-1}b(h^{-1},p(x))^{-1},h^{-1}.p(x))$ is a $G$- and $H$-equivariant homeomorphism with inverse $\Theta^{-1}: \: G \times Y \to X \times H, \, (g,y) \ma (g.q(y),a(g,q(y))h(y))$. Thus, the $G \times H$-space $\Omega = X \times H$ and $\bar{X} = X \times \gekl{e}$, $\bar{Y} = \Theta^{-1}(\gekl{e} \times Y)$ yield the desired topologically free $(G,H)$ coupling $G \curvearrowright {}_{\bar{Y}}\Omega {}_{\bar X} \curvearrowleft H$.

Note that topological freeness of $G \curvearrowright X$ and $H \curvearrowright Y$ ensures that $a$ and $b$ satisfy the cocycle identities (as in \cite[Lemma~2.8]{Li1}), which are needed in the preceding computations.
\setlength{\parindent}{0cm} \setlength{\parskip}{0cm}

\subsubsection{One-to-one correspondence}

\bproof[Proof of Theorem~\ref{1-1}]
It is straightforward to check that the constructions described in \S~\ref{tc->coc} and \S~\ref{coc->tc} are inverse to each other up to isomorphism. If we start with a topologically free $(G,H)$ coupling $G \curvearrowright {}_{\bar{Y}}\Omega {}_{\bar X} \curvearrowleft H$, construct a continuous orbit couple and then again a $(G,H)$ coupling, we end up with a $(G,H)$ coupling of the form $G \curvearrowright {}_{\ti{Y}}\ti{\Omega} {}_{\ti{X}} \curvearrowleft H$ where $\ti{\Omega} = \bar{X} \times H \cong G \times \bar{Y}$, $\ti{X} = \bar{X} \times \gekl{e}$ and $\ti{Y} \cong \gekl{e} \times \bar{Y}$. It is then obvious that $\ti{\Omega} = \bar{X} \times H \to \Omega, \, (x,h) \ma xh$ is an isomorphism of the couplings $G \curvearrowright {}_{\ti{Y}} \ti{\Omega} {}_{\ti{X}} \curvearrowleft H$ and $G \curvearrowright {}_{\bar{Y}}\Omega {}_{\bar{X}} \curvearrowleft H$. Conversely, if we start with a continuous orbit couple for topologically free systems $G \curvearrowright X$ and $H \curvearrowright Y$, construct a $(G,H)$ coupling and then again a $(G,H)$ continuous orbit couple, we end up with a continuous orbit couple for $G \curvearrowright \ti{X}$ and $H \curvearrowright \ti{Y}$ where $\ti{X} = X \times \gekl{e}$ and $\ti{Y} \cong \gekl{e} \times Y$. The canonical isomorphisms $X \cong X \times \gekl{e}$ and $Y \cong \gekl{e} \times Y$ yield an isomorphism between the original $(G,H)$ continuous orbit couple and the one we obtained at the end.
\setlength{\parindent}{0.5cm} \setlength{\parskip}{0cm}

Additional property (i) is clear from our constructions. For (ii), take $\bar X = \bar Y$ in the construction of \S~\ref{coc->tc}. Then it is clear that our maps $p$ and $q$ become the identity map on $\bar X = \bar Y$, that $\gamma$ becomes the constant function with value $e \in G$ and $\eta$ the constant function with value $e \in H$. Hence it is obvious that our construction yields a $(G,H)$ continuous orbit equivalence (see also Remark~\ref{X=Y-->coc=coe}).
\eproof
\setlength{\parindent}{0cm} \setlength{\parskip}{0.5cm}

\bremark
The maps $p$, $q$ constructed in \S~\ref{tc->coc} are open. Thus the maps $p$, $q$ appearing in a continuous orbit couple (Definition~\ref{coc}) are automatically open. This is also easy to see directly from the definition.
\eremark
\setlength{\parindent}{0cm} \setlength{\parskip}{0cm}

\subsection{Continuous orbit couples and Kakutani equivalence}
\label{coc-K}

\bdefin{\rm (Compare also \cite[Definition~4.1]{Mat}.)}
Topological dynamical systems $G \curvearrowright X$ and $H \curvearrowright Y$ are Kakutani equivalent if there exist clopen subsets $A \subseteq X$ and $B \subseteq Y$ such that $G.A = X$, $H.B = Y$ and $(X \rtimes G) \vert A \cong (Y \rtimes H) \vert B$ as topological groupoids. Here $(X \rtimes G) \vert A = s^{-1}(A) \cap r^{-1}(A)$ and $(Y \rtimes H) \vert B = s^{-1}(B) \cap r^{-1}(B)$. 
\edefin
\bremark
\label{Kaku=coe}
$(X \rtimes G) \vert A$ is (isomorphic to) the transformation groupoid attached to the partial action $G \curvearrowright A$ which is obtained by restricting $G \curvearrowright X$ to $A$. Similarly, $(Y \rtimes H) \vert B$ is (isomorphic to) the transformation groupoid attached to the partial action $H \curvearrowright B$ which is obtained by restricting $H \curvearrowright Y$ to $B$. In view of this, two topologically free systems $G \curvearrowright X$ and $H \curvearrowright Y$ are Kakutani equivalent if and only if there exist clopen subsets $A \subseteq X$ and $B \subseteq Y$ with $G.A = X$, $H.B = Y$ such that the partial actions $G \curvearrowright A$ and $H \curvearrowright B$ are continuously orbit equivalent in the sense of \cite{Li2}. This follows from \cite[Theorem~2.7]{Li2}.
\eremark
\setlength{\parindent}{0cm} \setlength{\parskip}{0cm}
The reader may find more about partial actions in \cite[\S~2]{Li2} and the relevant references in \cite{Li2}.
\setlength{\parindent}{0cm} \setlength{\parskip}{0cm}

\btheo
\label{COC--Kaku}
Let $G \curvearrowright X$ and $H \curvearrowright Y$ be topologically free systems. There exists a continuous orbit couple for $G \curvearrowright X$ and $H \curvearrowright Y$ with $p(X)$ closed if and only if $G \curvearrowright X$ and $H \curvearrowright Y$ are Kakutani equivalent.
\etheo
\setlength{\parindent}{0cm} \setlength{\parskip}{0cm}

Here $p: \: X \to Y$ is as in Definition~\ref{coc}. The assumption that $p(X)$ is closed always holds if $X$ is compact. This will be the case of interest later on.
\bproof
By Remark~\ref{Kaku=coe}, we have to show that there exists a continuous orbit couple for $G \curvearrowright X$ and $H \curvearrowright Y$ if and only if there exist clopen subspaces $A \subseteq X$ and $B \subseteq Y$ with $X = G.A$ and $Y = H.B$ such that the partial actions $G \curvearrowright A$ and $H \curvearrowright B$ are continuously orbit equivalent.
\setlength{\parindent}{0cm} \setlength{\parskip}{0.5cm}

For \an{$\Rarr$}, suppose we are given a continuous orbit couple for  $G \curvearrowright X$ and $H \curvearrowright Y$, and let $p$, $q$, $a$, $b$, $g$ and $h$ be as in Definition~\ref{coc} and \S~\ref{coc->tc}. For $g \in G$, let $U_g = \menge{x \in X}{g(x) = g}$. Then $U_g$ is clopen, and $X = \bigsqcup_{g \in G} U_g$. For every $g \in G$, $V_g \defeq p(U_g)$ is clopen, and $p: \: U_g \to V_g$ is a homeomorphism, whose inverse is given by $V_g \to U_g, \, y \ma g^{-1}.q(y)$. Set $B \defeq p(X)$. By assumption, $B$ is closed, hence clopen. We have $B = \bigcup_{g \in G} V_g$. As $G$ is countable, this is a countable union. Hence by inductively choosing compact open subspaces $B_g$ of $V_g$, we can arrange that $B$ is the disjoint union $B = \bigsqcup_{g \in G} B_g$. Let $A_g \defeq U_g \cap p^{-1}(B_g)$ and  $A \defeq \bigsqcup_{g \in G} A_g$. As every $A_g$ is clopen, $A = \bigsqcup_{g \in G} A_g$ is clopen in $X = \bigsqcup_{g \in G} U_g$. Set $\varphi \defeq p \vert_A = \bigsqcup_{g \in G} p \vert_{A_g}$. By construction, $\varphi$ is a homeomorphism with inverse $\varphi^{-1} = \bigsqcup_{g \in G} (p \vert_{A_g})^{-1} = \bigsqcup_{g \in G} (g^{-1}.q) \vert_{B_g}$.
\setlength{\parindent}{0.5cm} \setlength{\parskip}{0cm}

We have $\varphi(g.x) = p(g.x) = a(g,x).p(x)$ for all $x \in A$, $g \in G$ with $g.x \in A$. Moreover, take $y \in B_{g_1}$ and $h \in H$ with $h.y \in B_{g_2}$. Then $\varphi^{-1}(h.y) = g_2^{-1}.q(h.y) = g_2^{-1} b(h,y).q(y) = g_2^{-1} b(h,y) g_1. \varphi^{-1}(y)$. Define a map $b'$ by setting $b'(h,y) = g_2^{-1} b(h,y) g_1$ if $y \in B_{g_1} \cap h^{-1}.B_{g_2}$. Then $b'$ is continuous, and we have $\varphi^{-1}(h.y) = b'(h,y).\varphi^{-1}(y)$ for all $y \in B$, $h \in H$ with $h.y \in B$. This shows that $\varphi$ gives rise to a continuous orbit equivalence for $G \curvearrowright A$ and $H \curvearrowright B$. To see that $G.A = X$, take for $x' \in X$ an $x \in A$ such that $p(x) = p(x')$. Then $g(x).x = q(p(x)) = q(p(x')) = g(x').x'$, and therefore $x' \in G.x$. To see $H.B = Y$, take $y \in Y$ arbitrary. Then $p(q(y)) = h(y).y$ shows that $y = h(y)^{-1}.p(q(y)) \in H.B$. This shows \an{$\Rarr$}.
\setlength{\parindent}{0cm} \setlength{\parskip}{0.5cm}

For \an{$\Larr$}, suppose that $G \curvearrowright X$ and $H \curvearrowright Y$ are Kakutani equivalent, i.e., there are clopen subsets $A \subseteq X$ and $B \subseteq Y$ with $X = G.A$, $Y = H.B$ and the partial actions $G \curvearrowright A$ and $H \curvearrowright B$ are continuously orbit equivalent via a homeomorphism $\varphi: \: A \cong B$. By definition of continuous orbit equivalence (see \cite{Li2}), there exist continuous maps $a'$ and $b'$ satisfying $\varphi(g.x) = a'(g,x).\varphi(x)$ and $\varphi^{-1}(h.y) = b'(h,y).\varphi^{-1}(y)$ whenever this makes sense. As $X = G.A$, we can find clopen subsets $X_{\gamma} \subseteq \gamma.A$, $\gamma \in G$, such that $X = \bigsqcup_{\gamma \in G} X_{\gamma}$ and $X_e = A$. Define $p: \: X \to Y$ by setting $p(x) \defeq \varphi(\gamma^{-1}.x)$ for $x \in X_{\gamma}$. $p$ is continuous, and $p(X) = B$ is clopen. Similarly, there are clopen subsets $Y_{\eta} \subseteq \eta.B$ such that $Y = \bigsqcup_{\eta \in H} Y_{\eta}$ and $Y_e = B$. We define $q: \: Y \to X$ by setting $q(y) = \varphi^{-1}(\eta^{-1}.y)$ if $y \in Y_{\eta}$. By construction, $q$ is continuous.
\setlength{\parindent}{0.5cm} \setlength{\parskip}{0cm}

We have $p(g.x) = \varphi(\gamma_2^{-1}g.x) = \varphi(\gamma_2^{-1} g \gamma_1.(\gamma_1^{-1}.x)) = a'(\gamma_2^{-1} g \gamma_1, \gamma_1^{-1}.x).\varphi(\gamma_1^{-1}.x) 
  = a'(\gamma_2^{-1} g \gamma_1, \gamma_1^{-1}.x).p(x)$ for $x \in X_{\gamma_1}$ and $g \in G$ with $g.x \in X_{\gamma_2}$. Set $a: \: G \times X \to H$, $a(g,x) = a'(\gamma_2^{-1} g \gamma_1, \gamma_1^{-1}.x)$ for $x \in X_{\gamma_1} \cap g^{-1}.X_{\gamma_2}$. Then $a$ is continuous and $\varphi(g.x) = a(g,x).\varphi(x)$ for all $g \in G$ and $x \in X$.

For $y \in Y_{\eta_1}$ and $h \in H$ such that $h.y \in Y_{\eta_2}$, we have
$$
  q(h.y) = \varphi^{-1}(\eta_2^{-1} h.y) = \varphi^{-1}(\eta_2^{-1} h \eta_1 (\eta_1^{-1}.y)) = b'(\eta_2^{-1} h \eta_1, \eta_1^{-1}.y).\varphi^{-1}(\eta_1^{-1}.y)
  = b'(\eta_2^{-1} h \eta_1, \eta_1^{-1}.y).q(y).
$$
Set $b: \: H \times Y \to G$, $b(h,y) = b'(\eta_2^{-1} h \eta_1, \eta_1^{-1}.y)$ for $y \in Y_{\eta_1} \cap h^{-1}.Y_{\eta_2}$. Then $b$ is continuous and $\varphi^{-1}(h.y) = b(h,y).\varphi^{-1}(y)$ for all $h \in H$ and $y \in Y$.

Moreover, for $x \in X_{\gamma}$, $q(p(x)) = q(\varphi(\gamma^{-1}.x)) = \varphi^{-1}(\varphi(\gamma^{-1}.x)) = \gamma^{-1}.x$. Set $g: \: X \to G$, $g(x) = \gamma^{-1}$ if $x \in X_{\gamma}$. Then $g$ is continuous and $q(p(x)) = g(x).x$ for all $x \in X$. For $y \in Y_{\eta} \cap \eta.\varphi(X_{\gamma})$, we have $p(q(y)) = p(\varphi^{-1}(\eta^{-1}.y)) = \varphi(\gamma^{-1}.\varphi^{-1}(\eta^{-1}.y)) = \varphi(\gamma^{-1} b'(\eta^{-1},y).\varphi^{-1}(y) = a'(\gamma^{-1} b'(\eta^{-1},y), \varphi^{-1}(y)).y$. Set $h: \: Y \to H$, $h(y) \defeq a'(\gamma^{-1} b'(\eta^{-1},y), \varphi^{-1}(y))$ if $y \in Y_{\eta} \cap \eta.\varphi(X_{\gamma})$. Then $h$ is continuous and $p(q(y)) = h(y).y$ for all $y \in Y$.

So $p$ and $q$ give a continuous orbit couple for $G \curvearrowright X$ and $H \curvearrowright Y$. This shows \an{$\Larr$}.
\setlength{\parindent}{0cm} \setlength{\parskip}{0.5cm}
\eproof

\subsection{Dynamic characterizations of coarse embeddings, equivalences and bijections}

Putting together Theorem~\ref{THM:Gromov}, Lemma~\ref{Lem:Coupling-topfree}, Theorems~\ref{1-1} and \ref{COC--Kaku}, we obtain the following
\btheo
\label{u,q,l--DS}
\setlength{\parindent}{0cm} \setlength{\parskip}{0cm}

Let $G$ and $H$ be countable discrete groups. 
\begin{itemize}
\item The following are equivalent: 
\subitem- There exists a coarse embedding $G \to H$.
\subitem- There exist Kakutani equivalent topologically free $G \curvearrowright X$ and $H \curvearrowright Y$, with $X$ compact.
\subitem- There is a continuous orbit couple for topologically free $G \curvearrowright X$ and $H \curvearrowright Y$, with $X$ compact.
\item The following are equivalent:
\subitem- There is a coarse equivalence $G \to H$.
\subitem- There are Kakutani equivalent topologically free $G \curvearrowright X$ and $H \curvearrowright Y$ on compact spaces $X$, $Y$.
\subitem- There is a continuous orbit couple for topologically free $G \curvearrowright X$ and $H \curvearrowright Y$, with $X$, $Y$ compact.
\item There is a bijective coarse equivalence $G \to H$ if and only if there exist continuously orbit equivalent topologically free systems $G \curvearrowright X$ and $H \curvearrowright Y$ on compact spaces $X$ and $Y$.
\end{itemize}
In all these statements, the spaces $X$ and $Y$ can be chosen to be totally disconnected and second countable.
\etheo
This is a generalization of \cite[Theorem~3.2]{MST}, where the authors independently prove the last item of our theorem in the special case of finitely generated groups.
\setlength{\parindent}{0cm} \setlength{\parskip}{0.5cm}

\bremark
\label{Kaku->sCOE}
The last observation in Theorem~\ref{u,q,l--DS} says that we can always choose our spaces $X$, $Y$ to be totally disconnected. In that case, \cite[Theorem~3.2]{CRS} tells us that we can replace Kakutani equivalence in the theorem above by stable continuous orbit equivalence. Two topological dynamical systems $G \curvearrowright X$ and $H \curvearrowright Y$ are called stably continuously orbit equivalent if $\Zz \times G \curvearrowright \Zz \times X$ and $\Zz \times H \curvearrowright \Zz \times Y$ are continuously orbit equivalent. Here the integers $\Zz$ act on themselves by translation. 
\eremark

\subsection{Dynamic characterizations of coarse embeddings, equivalences and bijections in terms of actions on Stone-{\v C}ech compactifications}

Inspired by \cite{SW}, we characterize coarse embeddings, equivalences and bijections in terms of Kakutani equivalence (or stable continuous orbit equivalence) and continuous orbit equivalence of actions on Stone-{\v C}ech compactifications.

Let $G$, $H$ be two countable discrete groups. Let $\varphi: \: G \to H$ be a coarse embedding. Consider the Stone-{\v C}ech compactification $\beta G$ of $G$. It is homeomorphic to the spectrum $\Spec(\ell^{\infty}(G))$, and can be identified with the space of all ultrafilters on $G$. We will think of elements in $\beta G$ as ultrafilters on $G$. Given any subset $X \subseteq G$, we obviously have the identification $\menge{\cF \in \beta G}{X \in \cF} \cong \beta X, \, \cF \ma \cF \cap X \defeq \menge{F \cap X}{F \in \cF}$.

Now suppose that $X \subseteq G$ is a subset such that $\varphi \vert_X$ is injective. Setting $Y \defeq \varphi(X) \subseteq H$, we obtain a bijection $X \cong Y, \, x \ma \varphi(x)$, which we again denote by $\varphi$. Let us consider the topological dynamical systems $G \curvearrowright \beta G$ and $H \curvearrowright \beta H$. We identify $\beta X$ and $\beta Y$ with clopen subsets of $\beta G$ and $\beta H$, respectively, in the way explained above. $\varphi$ induces a homeomorphism $\beta \varphi: \: \beta X \cong \beta Y, \, \cF \ma \varphi(\cF)$. The dynamical systems $G \curvearrowright \beta G$ and $H \curvearrowright \beta H$ restrict to partial dynamical systems $G \curvearrowright \beta X$ and $H \curvearrowright \beta Y$.
\bprop
\label{phiX}
$\beta \varphi$ induces a continuous orbit equivalence between $G \curvearrowright \beta X$ and $H \curvearrowright \beta Y$, in the sense of \cite[Definition~2.6]{Li2}.
\eprop
\setlength{\parindent}{0cm} \setlength{\parskip}{0cm}

\bproof
For all $g \in G$, we need to find a continuous map $a: \: \gekl{g} \times U_{g^{-1}} \to H$ with $\beta \varphi(g.\cF) = a(g,\cF).\beta \varphi(\cF)$. Here $U_{g^{-1}} = \beta X \cap g^{-1}. \beta X = \menge{\cF \in \beta X}{g.\cF \in \beta X} = \menge{\cF \in \beta G}{X \in \cF, \, g^{-1}X \in \cF} \cong \beta(X \cap g^{-1}X)$. For $x \in X \cap g^{-1}X$, define the ultrafilter $\cF_x$ by saying that $Z \in \cF_x$ if and only if $x \in Z$. Define a map $\ti{a}: \: \gekl{g} \times \menge{\cF_x}{x \in X \cap g^{-1}X} \to H$ by setting $\ti{a}(g,\cF_x) \defeq \varphi(gx) \varphi(x)^{-1}$. Then
\begin{align}
\label{agFx}
  \ti{a}(g,\cF_x).\beta \varphi(\cF_x) = \varphi(gx) \varphi(x)^{-1}. \beta \varphi(\cF_x) = \varphi(gx) \varphi(x)^{-1}. \cF_{\varphi(x)} = \cF_{\varphi(gx)} = \beta \varphi(\cF_{gx}) = \beta \varphi(g.\cF_x)
\end{align}
for all $g \in G$, $x \in X \cap g^{-1}X$. Let us fix $g \in G$. Set $S = \menge{(gx,x)}{x \in G}$. As $\varphi$ is a coarse embedding and $\menge{st^{-1}}{(s,t) \in S} = \gekl{g}$ is finite, $\menge{\varphi(s) \varphi(t)^{-1}}{(s,t) \in S} = \menge{\varphi(g,x) \varphi(x)^{-1}}{x \in G}$ is finite. Hence
$\img(\ti{a}) \subseteq \menge{\varphi(gx) \varphi(x)^{-1}}{x \in G}$ is finite, hence a compact subset of $H$. By universal property of $\beta(X \cap g^{-1}X)$, there exists a continuous extension of $\ti{a}$ to $\gekl{g} \times U_{g^{-1}}$ which we denote by $a$. We claim that $\beta \varphi(g.\cF) = a(g,\cF).\beta \varphi(\cF)$ for all $\cF \in U_{g^{-1}}$. Let $x_i \in X \cap g^{-1}X$ be a net such that $\lim_i \cF_{x_i} = \cF$. Then $a(g,\cF_{x_i}) = \varphi(gx_i) \varphi(x_i)^{-1}$ converges to $a(g,\cF)$ by construction. Hence
\begin{equation*}
  a(g,\cF).\beta \varphi(\cF) = \lim_i a(g,\cF_{x_i}).\beta \varphi(\cF_{x_i})
  \overset{\eqref{agFx}}{=} \lim_i \beta \varphi(g.\cF_{x_i}) = \beta \varphi(\lim_i g.\cF_{x_i}) = \beta \varphi(g.\cF). 
\qedhere
\end{equation*}
\eproof
\setlength{\parindent}{0cm} \setlength{\parskip}{0cm}

The following observation will be used several times. 
\blemma
\label{Lem:UE:X->Y}
Let $\varphi: \: G \to H$ be a coarse embedding. Set $Y \defeq \varphi(G)$. For every $y \in Y$, choose $x_y \in G$ with $\varphi(x_y) = y$. Set $X \defeq \menge{x_y}{y \in Y}$.
\setlength{\parindent}{0.5cm} \setlength{\parskip}{0cm}

Then $\varphi$ restricts to a bijection $X \cong Y$, and there is finite subset $F \subseteq G$ with $G = \bigcup_{g \in F} gX$. 
\elemma
\bproof
Clearly, the restriction of $\varphi$ to $X$ is a bijection onto $Y$. To prove that $G$ can be covered by finitely many translates of $X$, set $S \defeq \menge{(g,x_{\varphi(g)})}{g \in G}$.
Then $\menge{\varphi(s) \varphi(t)^{-1}}{(s,t) \in S} = \gekl{e}$, where $e$ is the identity in $H$. Since $\varphi$ is a coarse embedding, $\menge{g x_{\varphi(g)}^{-1}}{g \in G} = \menge{st^{-1}}{(s,t) \in S}$ must be finite. Hence there is finite subset $F \subseteq G$ with $G = \bigcup_{g \in F} gX$.
\eproof

We now obtain the following characterizations of coarse embeddings, equivalences and bijections.
\bcor
\label{qibilip_StoneCech}
\setlength{\parindent}{0cm} \setlength{\parskip}{0cm}

Let $G$ and $H$ be countable discrete groups. 
\begin{itemize}
\item[(i)]
The following are equivalent:
\subitem- There is a coarse embedding $G \to H$.
\subitem- There is an open, dense, $H$-invariant subspace $\ti{Y} \subseteq \beta H$ such that $G \curvearrowright \beta G$ and $H \curvearrowright \ti{Y}$ are Kakutani equivalent.
\subitem- There is an open, dense, $H$-invariant subspace $\ti{Y}\subseteq \beta H$ such that there is a continuous orbit couple for $G \curvearrowright \beta G$ and $H \curvearrowright \ti{Y}$.
\subitem- There is an open, dense, $H$-invariant subspace $\ti{Y} \subseteq \beta H$ such that $G \curvearrowright \beta G$ and $H \curvearrowright \ti{Y}$ are stably continuously orbit equivalent.
\item[(ii)]
There is a coarse equivalence $G \to H$ if and only if $\ti{Y} = \beta H$ works in the statements in (i).
\item[(iii)] There is a bijective coarse equivalence $G \to H$ if and only if $G \curvearrowright \beta G$ and $H \curvearrowright \beta H$ are continuously orbit equivalent.
\end{itemize}
\ecor
\bproof
(i): Let $\varphi: \: G \to H$ be a coarse embedding. Let $Y$ and $X$ be as in Lemma~\ref{Lem:UE:X->Y}. As the restriction of $\varphi$ to $X$ is a bijection onto $Y$, Proposition~\ref{phiX} yields that $G \curvearrowright \beta X$ and $H \curvearrowright \beta Y$ are continuously orbit equivalent. As there is finite subset $F \subseteq G$ with $G = \bigcup_{g \in F} gX$, we have $\beta G = G. \beta X$. Let $\ti{Y} \defeq H. \beta Y$. Then $G \curvearrowright \beta G$ and $H \curvearrowright \ti{Y}$ are Kakutani equivalent. $\ti{Y}$ is $H$-invariant by construction, and it is easy to see that $\ti{Y}$ is open and dense. Now (i) follows from Theorem~\ref{u,q,l--DS}, Theorem~\ref{COC--Kaku} and Remark~\ref{Kaku->sCOE}.
\setlength{\parindent}{0.5cm} \setlength{\parskip}{0cm}

(ii): A coarse embedding $\varphi: \: G \to H$ is coarsely invertible if and only if there is a finite subset $F \subseteq H$ such that $H = \bigcup_{h \in F} h \varphi(G)$. This happens if and only if in the proof of (i), we get $\ti{Y} = \beta H$.

(iii): If $\varphi: \: G \to H$ is a bijective coarse equivalence, then we can take $X = G$, $Y = H$ in the above proof of (i) and obtain that $G \curvearrowright \beta G$ and $H \curvearrowright \beta H$ are continuously orbit equivalent. The reverse implication \an{$\Larr$} in (ii) is proven in Theorem~\ref{u,q,l--DS}.
\eproof
\setlength{\parindent}{0cm} \setlength{\parskip}{0.5cm}

\bremark
In combination with \cite{SW}, Corollary~\ref{qibilip_StoneCech} implies that nuclear Roe algebras have distinguished Cartan subalgebras, as explained in \cite{LR}.
\eremark

\bremark
Corollary~\ref{qibilip_StoneCech} shows that quasi-isometry rigidity can be interpreted as a special case of continuous orbit equivalence rigidity (in the sense of \cite{Li1}), applied to actions on Stone-{\v C}ech compactifications. This points towards an interesting connection between these two types of rigidity phenomena and would be worth exploring further.
\eremark

\section{Applications to (co)homology I}
\label{sec:HI}

We now show how the results in \cite{Sha,Sau} on coarse invariance of (co)homological dimensions and property $H_{FD}$ follow from Morita invariance of groupoid (co)homology. Let us first define groupoid (co)homology. We do this in a concrete and elementary way which is good enough for our purposes. We refer to \cite{CM} for a more general and more conceptual approach, and for more information about groupoids. Let $\cG$ be an {\'e}tale locally compact groupoid with unit space $X = \cG^{(0)}$, and $R$ a commutative ring with unit. A $\cG$-sheaf of $R$-modules is a sheaf $\cA$ of $R$-modules over $X$, i.e., we have a locally compact space $\cA$ with an {\'e}tale continuous surjection $\pi: \: \cA \onto X$ whose fibres are $R$-modules, together with the structure of a right $\cG$-space on $\cA$. In particular, ever $\gamma \in \cG$ induces an isomorphism of $R$-modules $\cA_{r(\gamma)} \to \cA_{s(\gamma)}, \, a \ma a * \gamma$. To pass from right to left actions, we write $\gamma.a \defeq a * \gamma^{-1}$ if $\pi(a) = s(\gamma)$.

Let $\cG^{(n)} = \menge{(\gamma_1, \dotsc, \gamma_n) \in \cG^n}{s(\gamma_i) = r(\gamma_{i+1}) \ {\rm for} \ {\rm all} \ 1 \leq i \leq n-1}$, and set $r(\gamma_1, \dotsc, \gamma_n) = r(\gamma_1)$. We write $\vecgamma$ for elements in $\cG^{(n)}$. Given a $\cG$-sheaf of $R$-modules $\cA$ with projection $\pi: \: \cA \onto X$, let $\Gamma_c(\cG^{(n)},\cA)$ be the $R$-module of continuous functions $f: \: \cG^{(n)} \to \cA$ with compact support such that $\pi(f(\vecgamma)) = r(\vecgamma)$. Now we define a chain complex
$\dotso \overset{d_{n+1}}{\lori} \Gamma_c(\cG^{(n)},\cA) \overset{d_n}{\lori} \Gamma_c(\cG^{(n-1)},\cA) \overset{d_{n-1}}{\lori} \dotso \overset{d_2}{\lori} \Gamma_c(\cG,\cA) \overset{d_1}{\lori} \Gamma_c(X,\cA) \to 0$, with $d_1(f)(x) = \sum_{\substack{\gamma \in \cG \\ s(\gamma) = x}} \gamma^{-1}.f(\gamma) - \sum_{\substack{\gamma \in \cG \\ r(\gamma) = x}} f(\gamma)$ for $f \in \Gamma_c(\cG,\cA)$, and for $n \geq 1$: $d_n(f) = \sum_{i=0}^n (-1)^i d_n^{(i)}(f)$ for $f \in \Gamma_c(\cG^{(n)},\cA)$, where
\begin{align*}
  d_n^{(0)}(f)(\gamma_1, \dotsc, \gamma_{n-1}) &= \sum_{\substack{\gamma_0 \in \cG \\ s(\gamma_0) = r(\gamma_1)}} \gamma_0^{-1}.f(\gamma_0, \gamma_1, \dotsc, \gamma_{n-1}), \\
  d_n^{(i)}(f)(\gamma_1, \dotsc, \gamma_{n-1}) &= \sum_{\substack{\eta, \xi \in \cG \\ \eta \xi = \gamma_i}} f(\dotsc, \gamma_{i-1}, \eta, \xi, \gamma_{i+1}, \dotsc) \ \ \ \ \ \ {\rm for} \ 1 \leq i \leq n-1,\\
  d_n^{(n)}(f)(\gamma_1, \dotsc, \gamma_{n-1}) &= \sum_{\substack{\gamma_n \in \cG \\ r(\gamma_n) = s(\gamma_{n-1})}} f(\gamma_1, \dotsc, \gamma_{n-1}, \gamma_n).  
\end{align*}
We then define the $n$-th homology group $H_n(\cG,\cA) \defeq \ker(d_n) / \img(d_{n+1})$. In the case $R = \Zz$ and where $\cA$ is a constant sheaf with trivial $\cG$-action, we recover \cite[Definition~3.1]{Mat}.

Let us also introduce cohomology. Let $\cG$, $R$ and $\cA$ be as above, and let $\Gamma(\cG^{(n)},\cA)$ be the $R$-module of continuous functions $f: \: \cG^{(n)} \to \cA$ with $\pi(f(\vecgamma)) = r(\vecgamma)$. 
We define a cochain complex $0 \to \Gamma(X,\cA) \overset{d^0}{\lori} \Gamma(\cG,\cA) \overset{d^1}{\lori} \dotso \overset{d^{n-1}}{\lori} \Gamma(\cG^{(n)},\cA) \overset{d^n}{\lori} \Gamma(\cG^{(n+1)},\cA) \overset{d^{n+1}}{\lori} \dotso $ with $d^0(f)(\gamma) = \gamma.f(s(\gamma)) - f(r(\gamma))$, and for $n \geq 1$: $d^n(f) = \sum_{i=0}^{n+1} (-1)^i d_{(i)}^n(f)$, where
\begin{align*}
  d_{(0)}^n(f)(\gamma_0, \dotsc, \gamma_n) &= \gamma_0.f(\gamma_1, \dotsc, \gamma_n);\\
  d_{(i)}^n(f)(\gamma_0, \dotsc, \gamma_n) &= f(\gamma_0, \dotsc, \gamma_{i-1} \gamma_i, \dotsc, \gamma_n) \ \ \ \ \ \ {\rm for} \ 1 \leq i \leq n;\\
  d_{(n+1)}^n(f)(\gamma_0, \dotsc, \gamma_n) &= f(\gamma_0, \dotsc, \gamma_{n-1}).
\end{align*}
We set $H^n(\cG,\cA) \defeq \ker(d^n) / \img(d^{n-1})$.

Now let $G \curvearrowright X$ be a topological dynamical system. For notational purposes, and to keep the conventions in the literature, let us pass to the right action $X \curvearrowleft G$, $x.g = g^{-1}.x$, and consider the corresponding transformation groupoid $X \rtimes G$ with source and range maps given by $s(x,g) = x.g$, $r(x,g) = x$. We note that the transformation groupoid $G \ltimes X$ attached to the original action, as in \cite{Li1,Li2}, is isomorphic to $X \rtimes G$ via $G \ltimes X \to X \rtimes G, \, (g,x) \ma (g.x,g)$. It is easy to see that a ($X \rtimes G$)-sheaf of $R$-modules is nothing else but a sheaf $\cA$ of $R$-modules over $X$, $\pi: \: \cA \onto X$, together with a left $G$-action on $\cA$ via homeomorphisms (denoted by $G \times \cA \to \cA, \, (g,a) \ma g.a$) such that $\pi$ becomes $G$-equivariant, and $\cA_x \to \cA_{g.x}, \, a \ma g.a$ is an isomorphism of $R$-modules. We call these $G$-sheaves of $R$-modules over $X$.

\subsection{Isomorphisms in homology and cohomology}

First of all, let us prove
\setlength{\parindent}{0cm} \setlength{\parskip}{0cm}
\btheo
\label{hom,cohom}
Let $G \curvearrowright X$ and $H \curvearrowright Y$ be topologically free systems, where $G$ and $H$ are countable discrete groups. Suppose that $G \curvearrowright X$ and $H \curvearrowright Y$ are Kakutani equivalent. Then there is an equivalence of categories between $G$-sheaves of $R$-modules over $X$ and $H$-sheaves of $R$-modules over $Y$, denoted by $\cS_X \ma \cS_Y$ on the level of objects, such that $H_*(G,\Gamma_c(X,\cS_X)) \cong H_*(H,\Gamma_c(Y,\cS_Y))$ and $H^*(G,\Gamma(X,\cS_X)) \cong H^*(H,\Gamma(Y,\cS_Y))$.
\etheo
Here $\Gamma$ stands for continuous sections and $\Gamma_c$ for those with compact support.

\bproof
It is easy to see that $H_*(G,\Gamma_c(X,\cA)) \cong H_*(X \rtimes G,\cA)$ and $H^*(G,\Gamma(X,\cA)) \cong H^*(X \rtimes G,\cA)$ for topological dynamical systems $G \curvearrowright X$ and $G$-sheaves $\cA$ of $R$-modules over $X$.
\setlength{\parindent}{0.5cm} \setlength{\parskip}{0cm}

Now, by assumption, there are clopen subspaces $A \subseteq X$ and $B \subseteq Y$ with $X = G.A$, $Y = H.B$ and an isomorphism of topological groupoids $\chi: \: (X \rtimes G) \vert A \cong (Y \rtimes H) \vert B$. Let $\iota_A: \: (X \rtimes G) \vert A \into X \rtimes G$ and $\iota_B: \: (Y \rtimes H) \vert B \into Y \rtimes H$ be the canonical inclusions. As $A$ is $G$-full and $B$ is $H$-full, $\iota_A$ and $\iota_B$ induce equivalences of categories of sheaves. So we obtain an equivalence of categories between $G$-sheaves of $R$-modules over $X$ and $H$-sheaves of $R$-modules over $Y$, denoted by $\cS_X \ma \cS_Y$ on the level of objects, such that $\cS_Y$ is uniquely determined by $\chi^*(\cS_Y \vert B) = \cS_X \vert A$. Our theorem now follows from Morita invariance of groupoid (co)homology.
\eproof

For every topological dynamical system $G \curvearrowright X$, we have $\sup \menge{n}{H_n(G, \Gamma_c(X,\cA)) \ncong \gekl{0}} \leq {\rm hd}_R(G)$ and $\sup \menge{n}{H^n(G,\Gamma(X,\cA)) \ncong \gekl{0}} \leq {\rm cd}_R(G)$ by the definitions of homological and cohomological dimensions. Here the suprema are taken over all $G$-sheaves $\cA$ of $R$-modules over $X$.
\bdefin
A $(G,H)$ continuous orbit couple is called $H_{*,R}G$-full if $\sup \menge{n}{H_n(G, \Gamma_c(X,\cA)) \ncong \gekl{0}} = {\rm hd}_R(G)$ holds for its topological dynamical system $G \curvearrowright X$. It is called $H^{*,R}G$-full if its topological dynamical system $G \curvearrowright X$ satisfies $\sup \menge{n}{H^n(G,\Gamma(X,\cA)) \ncong \gekl{0}} = {\rm cd}_R(G)$.
\edefin
The following is an immediate consequence of Theorem~\ref{hom,cohom}.
\bcor
\label{Hfull}
If there exists an $H_{*,R}G$-full topologically free $(G,H)$ continuous orbit couple, then ${\rm hd}_R(G) \leq {\rm hd}_R(H)$. If there exists an $H^{*,R}G$-full topologically free $(G,H)$ continuous orbit couple, then ${\rm cd}_R(G) \leq {\rm cd}_R(H)$.
\ecor
\bremark
Together with Theorem~\ref{u,q,l--DS}, Corollary~\ref{Hfull} can be viewed as an explanation and generalization of the results in \cite{Sha,Sau} concerning coarse invariance of (co)homological dimension. In our terminology, the conditions from \cite{Sha,Sau} that the topological dynamical system $G \curvearrowright X$ of a $(G,H)$ continuous orbit couple admits a $G$-invariant probability measure and $\Qz \subseteq R$ ensure that the $(G,H)$ continuous orbit couple is $H_{*,R}G$-full and $H^{*,R}G$-full (see \cite[\S~3.3]{Sha} and \cite[\S~4]{Sau}). Existence of a $G$-invariant probability measure is guaranteed if $G$ is amenable and the $G$-space of our continuous orbit couple is compact. Moreover, again in our terminology, it is shown in \cite[\S~4]{Sau} that a $(G,H)$ continuous orbit couple with compact $G$-space is $H_{*,R}G$-full if ${\rm hd}_R(G) < \infty$ and $H^{*,R}G$-full if ${\rm cd}_R(G) < \infty$. Once we know this, \cite[Theorem~1.5]{Sha} and \cite[Theorem~1.2]{Sau} are immediate consequences of Theorem~\ref{u,q,l--DS} and Corollary~\ref{Hfull}. In \S~\ref{ss:Consequences}, we present an alternative approach to these results.
\eremark
\setlength{\parindent}{0cm} \setlength{\parskip}{0.5cm}

\subsection{Isomorphisms in reduced cohomology}
\label{IsomRedCohom}

Let $\cG$ be an \'{e}tale locally compact groupoid and $\fL = (\mu,\cH,L)$ a (unitary) representation of $\cG$ as in \cite[Chapter~II, Definition~1.6]{Ren}. Here $\mu$ is a quasi-invariant measure on $\cG^{(0)}$, $\cH$ a Hilbert bundle over $(\cG^{(0)},\mu)$, and $L$ a representation of $\cG$, i.e., for each $\gamma \in \cG$, $L(\gamma)$ is a unitary $\cH_{s(\gamma)} \cong \cH_{r(\gamma)}$, and the conditions in \cite[Chapter~II, Definition~1.6]{Ren} are satisfied ($\sigma$ in \cite[Chapter~II, Definition~1.6]{Ren} is the trivial cocycle in our case). Let $D$ be the modular function attached to $\mu$, as in \cite[Chapter~I, Definition~3.4]{Ren}. In particular, we are interested in the case $\cG = X \rtimes G$ of a transformation groupoid attached to a topological dynamical system $G \curvearrowright X$ on a compact space $X$. A representation $\fL$ of $X \rtimes G$ gives rise -- through its integrated form -- to a *-representation of $C(X) \rtimes G$, which in turn corresponds in a one-to-one way to a covariant representation $(\pi_{\fL},\sigma_{\fL})$ of $G \curvearrowright X$ (or rather of $(C(X),G)$). 

Now let $\cG = X \rtimes G$ be as above and $\fL$ a representation of $\cG$. We define cohomology groups $H^n(\cG,\fL)$ and reduced cohomology groups ${\bar H}^n(\cG,\fL)$ with coefficients in $\fL$. Let us write $\fL = (\mu,\cH,L)$. Let $\cG^{(n)} = \menge{(\gamma_1, \dotsc, \gamma_n) \in \cG^n}{s(\gamma_i) = r(\gamma_{i+1}) \ {\rm for} \ {\rm all} \ 1 \leq i \leq n-1}$, and set $r(\gamma_1, \dotsc, \gamma_n) = r(\gamma_1)$. We will write $\vecgamma$ for elements in $\cG^{(n)}$. Let $\Gamma(\cG^{(n)},\cH)$ be the set of all Borel functions $f: \: \cG^{(n)} \to \cH$ with $f(\vecgamma) \in \cH_{r(\vecgamma)}$ such that for every compact subset $K \subseteq \cG^{(n)}$, $\int_{\cG^{(0)}} \sum_{\substack{\vecgamma \in K \\ r(\vecgamma) = x}} \norm{f(\vecgamma)}^2 d\mu(x) < \infty$, divided by the equivalence relation saying that $f_1 \sim f_2$ if for every compact subset $K \subseteq \cG^{(n)}$, $\int_{\cG^{(0)}} \sum_{\substack{\vecgamma \in K \\ r(\vecgamma) = x}} \norm{f_1(\vecgamma) - f_2(\vecgamma)}^2 d\mu(x) = 0$. The topology on $\Gamma(\cG^{(n)},\cH)$ is given by the following notion of convergence: A net $(f_i)_i$ converges to an element $f$ in $\Gamma(\cG^{(n)},\cH)$ if for every compact subset $K \subseteq \cG^{(0)}$, $\lim_{i \to \infty} \int_{\cG^{(0)}} \sum_{\substack{\vecgamma \in K \\ r(\vecgamma) = x}} \norm{f(\vecgamma) - f_i(\vecgamma)}^2 d\mu(x) = 0$. We define a cochain complex $0 \to \Gamma(\cG^{(0)},\cH) \overset{d^0}{\lori} \Gamma(\cG^{(1)},\cH) \overset{d^1}{\lori} \dotso$ with $d^0(f)(\gamma) = D^{-\frac{1}{2}}(\gamma) L(\gamma) f(s(\gamma)) - f(r(\gamma))$, and for $n \geq 1$: $d^n = \sum_{i=0}^{n+1} (-1)^i d_{(i)}^n$, where
\begin{align*}
  d_{(0)}^n(f)(\gamma_0, \dotsc, \gamma_n) &= D^{-\frac{1}{2}}(\gamma_0) L(\gamma_0) f(\gamma_1, \dotsc, \gamma_n);\\
  d_{(i)}^n(f)(\gamma_0, \dotsc, \gamma_n) &= f(\gamma_0, \dotsc, \gamma_{i-1} \gamma_i, \dotsc, \gamma_n) \ \ \ \ \ \ {\rm for} \ 1 \leq i \leq n;\\
  \delta_{(n+1)}^n(f)(\gamma_0, \dotsc, \gamma_n) &= f(\gamma_0, \dotsc, \gamma_{n-1}).
\end{align*}
It is easy to check that $d^n \circ d^{n-1} = 0$ for all $n \geq 1$. Thus $\img(d^{n-1}) \subseteq \ker(d^n)$. Since all the $d^n$ are continuous, we also have $\climg(d^{n-1}) \subseteq \ker(d^n)$. We set $H^n(\cG,\fL) \defeq \ker(d^n) / \img(d^{n-1})$ and ${\bar H}^n(\cG,\fL) \defeq \ker(d^n) / \climg(d^{n-1})$.

Our goal is to prove the following
\btheo
\label{redcohom}
Suppose there is a continuous orbit couple for topological dynamical systems $G \curvearrowright X$ and $H \curvearrowright Y$ on compact spaces $X$ and $Y$. Then there is a one-to-one correspondence between representations of $X \rtimes G$ and $Y \rtimes H$, denoted by $\fL \leftrightarrow \fM$, with $H^*(G,\sigma_{\fL}) \cong H^*(H,\sigma_{\fM})$ and ${\bar H}^*(G,\sigma_{\fL}) \cong {\bar H}^*(H,\sigma_{\fM})$.
\etheo
\setlength{\parindent}{0cm} \setlength{\parskip}{0cm}

For the definition of reduced cohomology ${\bar H}^*$, we refer to \cite[Chapitre~III]{Gui}.
\bproof
Clearly, $H^*(X \rtimes G,\fL) \cong H^*(G,\sigma_{\fL})$ and $\bar{H}^*(X \rtimes G,\fL) \cong \bar{H}^*(G,\sigma_{\fL})$.
\setlength{\parindent}{0.5cm} \setlength{\parskip}{0cm}

Now, if there is a continuous orbit couple for topological dynamical systems $G \curvearrowright X$ and $H \curvearrowright Y$ on compact spaces $X$ and $Y$, then by Theorem~\ref{COC--Kaku}, $G \curvearrowright X$ and $H \curvearrowright Y$ are Kakutani equivalent. So there exist clopen subspaces $A \subseteq X$ and $B \subseteq Y$ with $G.A = X$, $H.B = Y$, together with an isomorphism of topological groupoids $\chi: \: (X \rtimes G) \vert A \cong (Y \rtimes H) \vert B$. As $A$ is $G$-full and $B$ is $H$-full, we get one-to-one correspondences $\fL \leftrightarrow \fL \vert A$ and $\fM \leftrightarrow \fM \vert B$ between representations of $X \rtimes G$ and $(X \rtimes G) \vert A$, and between representations of $Y \rtimes H$ and $(Y \rtimes H) \vert B$, respectively. Thus we obtain a one-to-one correspondence between representations of $X \rtimes G$ and $Y \rtimes H$, denoted by $\fL \leftrightarrow \fM$, where $\fM$ is uniquely determined by $\chi^*(\fM \vert B) = \fL \vert A$. The theorem now follows from Morita invariance of groupoid (co)homology.
\eproof

\bremark
\label{Rem:Ind-Phi}
If the topological dynamical system $G \curvearrowright X$ is on a second countable space $X$, then every *-representation of $C_c(X \rtimes G)$ on a Hilbert space is the integrated form of a representation of $X \rtimes G$. Actually, *-representations of $C_c(X \rtimes G)$ and representations of $X \rtimes G$ are in one-to-one correspondence (see \cite[Chapter~II, Theorem~1.21 and Corollary~1.23]{Ren}). Thus we obtain a reformulation of Theorem~\ref{redcohom}: Suppose there is a continuous orbit couple for topological dynamical systems $G \curvearrowright X$ and $H \curvearrowright Y$ on second countable compact spaces $X$ and $Y$. 
By Theorem~\ref{COC--Kaku}, $G \curvearrowright X$ and $H \curvearrowright Y$ are Kakutani equivalent, so there exist clopen subspaces $A \subseteq X$ and $B \subseteq Y$ with $G.A = X$, $H.B = Y$, together with an isomorphism of topological groupoids $\chi: \: (X \rtimes G) \vert A \cong (Y \rtimes H) \vert B$. Let $\Phi: \: C^*((X \rtimes G) \vert A) \cong C^*((Y \rtimes H) \vert B)$ be the corresponding isomorphism of groupoid C*-algebras. Then the one-to-one correspondence $\fL \leftrightarrow \fM$ from Theorem~\ref{redcohom} translates to a one-to-one correspondence $(\pi,\sigma) \leftrightarrow (\rho,\tau)$ between covariant representations of $G \curvearrowright X$ and $H \curvearrowright Y$, where $(\rho,\tau)$ is uniquely determined (up to unitary equivalence) by the requirement that
$\rukl{\rho \rtimes \tau \vert_{C^*((Y \rtimes H) \vert B)}} \circ \Phi
  = \pi \rtimes \sigma \vert_{C^*((X \rtimes G) \vert A)}$. Here we view $C^*((Y \rtimes H) \vert B)$ and $C^*((X \rtimes G) \vert A)$ as full corners in $C(Y) \rtimes H$ and $C(X) \rtimes G$. We write $(\rho,\tau) = \Ind_{\Phi^{-1}}(\pi,\sigma)$ and $(\pi,\sigma) = \Ind_{\Phi}(\rho,\tau)$.
\eremark

\bcor
\label{2nd,redcohom}
Let $G \curvearrowright X$ and $H \curvearrowright Y$ be topological dynamical systems on second countable compact spaces $X$ and $Y$, and assume that there is a continuous orbit couple for $G \curvearrowright X$ and $H \curvearrowright Y$. Let $(\pi,\sigma) \leftrightarrow (\rho,\tau)$ be as in Remark~\ref{Rem:Ind-Phi}. Then we have $H^*(G,\sigma) \cong H^*(H,\tau)$ and ${\bar H}^*(G,\sigma) \cong {\bar H}^*(H,\tau)$.
\ecor

\bremark
Theorem~\ref{redcohom} and Corollary~\ref{2nd,redcohom} have natural analogues in homology, i.e., for $H_*$ and ${\bar H}_*$.
\eremark

\subsection{Coarse invariance of property $H_{FD}$}
\label{ss:HFD}

As a consequence of Theorem~\ref{redcohom}, we discuss coarse invariance of Shalom's property $H_{FD}$ from \cite{Sha}. In this section (\S~\ref{ss:HFD}), we assume that our spaces are second countable. Let us start with the following
\blemma
\label{XG=FullCorner}
Let $G \curvearrowright {}_{\bar{Y}}\Omega {}_{\bar{X}} \curvearrowleft H$ be a topological coupling, let $\alpha$ and $\beta$ be as in \S~\ref{tc->coc}, let $G \curvearrowright \bar{Y}$, $H \curvearrowright \bar{X}$ be the actions given by $g.x = g x \alpha(g,x)^{-1}$, $h.y = \beta(y,h^{-1})^{-1} x h^{-1}$, and let $\bar{X} \rtimes G$, $\bar{Y} \rtimes H$ be the corresponding transformation groupoids. Then
\begin{align*}
  & \bar{X} \rtimes G \to (\Omega \rtimes (G \times H)) \vert \bar{X}, \, (x,g) \ma (x,g,\alpha(g^{-1},x)^{-1}) \\
  & \bar{Y} \rtimes H \to (\Omega \rtimes (G \times H)) \vert \bar{Y}, \, (y,h) \ma (y,\beta(y,h),h)
\end{align*}
are isomorphisms of topological groupoids.
\elemma
\bproof
As $r(x,g) = x = r(x,g,\alpha(g^{-1},x)^{-1})$, $s(x,g) = x.g = g^{-1}.x = g^{-1} x \alpha(g^{-1},x)^{-1} = s(x,g,\alpha(g^{-1},x)^{-1})$ and $(x,g,\alpha(g^{-1},x)^{-1}) (g^{-1} x \alpha(g^{-1},x)^{-1}, \bar{g}, \alpha(\bar{g
}^{-1}, g^{-1} x \alpha(g^{-1},x)^{-1})^{-1}) = (x, g \bar{g}, \alpha((g \bar{g})^{-1}, x)^{-1})$, $\bar{X} \rtimes G \to (\Omega \rtimes (G \times H)) \vert \bar{X}, \, (x,g) \ma (x,g,\alpha(g^{-1},x)^{-1})$ is a groupoid homomorphism. It is clearly continuous, and $(\Omega \rtimes (G \times H)) \vert \bar{X} \to \bar{X} \rtimes G, \, (x,g,h) \ma (x,g)$ is its continuous inverse. The proof of the second claim is analogous.
\eproof

Given a topologically free $(G,H)$ continuous orbit couple which corresponds to the $(G,H)$ coupling $G \curvearrowright {}_{Y}\Omega {}_{X} \curvearrowleft H$ with compact $X$ and $Y$, the proof of Theorem~\ref{u,q,l--DS} provides a concrete way to construct Kakutani equivalent dynamical systems $G \curvearrowright X$ and $H \curvearrowright Y$ together with clopen subspaces $A \subseteq X$ and $B \subseteq Y$ such that $(X \rtimes G) \vert A \cong (Y \rtimes H) \vert B$. We need the following
\blemma
\label{A=B}
We can modify our $(G,H)$ continuous orbit couple above, without changing its topological dynamical system $G \curvearrowright X$, so that the described process yields a topological coupling and subspaces $A$, $B$ with $A = B$ as subspaces of $\Omega$.
\elemma
\bproof
In the proof of Theorem~\ref{COC--Kaku}, we had constructed $A$ and $B$ as disjoint unions $A = \bigsqcup_g A_g$ and $B = \bigsqcup_g B_g$. Following the construction of the continuous orbit couple out of our topological coupling in \S~\ref{tc->coc}, we see that these subspaces $A_g$ and $B_g$ were related by $g A_g = B_g$ in $\Omega$. Set $Y' \defeq (Y \setminus B) \sqcup A$. Then $X$ and $Y'$ are still fundamental domains for the $H$- and $G$-actions on $\Omega$. So we obtain a new topologically free $(G,H)$ coupling $G \curvearrowright {}_{Y'}\Omega {}_{X} \curvearrowleft H$. The construction in \S~\ref{tc->coc} yields a continuous orbit couple with new continuous orbit map $p': \: X \to Y'$ satisfying $p'(X) = A$. Hence our construction in the proof of Theorem~\ref{COC--Kaku} gives us the subspaces $A \subseteq X$ and $A \subseteq Y'$ implementing the Kakutani equivalence between $G \curvearrowright X$ and $H \curvearrowright Y'$.
\eproof
Let $G \curvearrowright {}_Y \Omega {}_X \curvearrowleft H$ and $G \curvearrowright X$, $H \curvearrowright Y$ be as above, with a clopen subspace $A \subseteq X \cap Y$ such that $G.A = X$, $H.A = Y$ and $(X \rtimes G) \vert A \cong (Y \rtimes H) \vert A$. Let $\Phi: \: C^*((X \rtimes G) \vert A) \cong C^*((Y \rtimes H) \vert A)$ be the induced C*-isomorphism. Lemma~\ref{XG=FullCorner} yields an isomorphism of C*-algebras $C(X) \rtimes G \cong 1_X (C_0(\Omega) \rtimes (G \times H)) 1_X$, and $1_X$ is a full projection. Therefore, $C(X) \rtimes G$ is Morita equivalent to $C_0(\Omega) \rtimes (G \times H)$, and a $C(X) \rtimes G - C_0(\Omega) \rtimes (G \times H)$-imprimitivity bimodule is given by $\fX = 1_X (C_0(\Omega) \rtimes (G \times H))$ (with respect to the identification $C(X) \rtimes G \cong 1_X (C_0(\Omega) \rtimes (G \times H)) 1_X$ provided by Lemma~\ref{XG=FullCorner}). We obtain (up to unitary equivalence) bijections between representations of $C(X) \rtimes G$ and representations of $C_0(\Omega) \rtimes (G \times H)$ and also between covariant representations of $G \curvearrowright X$ and $G \times H \curvearrowright \Omega$. We denote both of them by $\Ind_{\fX}$. Also, let $\fY$ be the $C(Y) \rtimes H - C_0(\Omega) \rtimes (G \times H)$-imprimitivity bimodule given by $1_Y (C_0(\Omega) \rtimes (G \times H))$ with respect to the identification $C(Y) \rtimes H \cong 1_Y (C_0(\Omega) \rtimes (G \times H)) 1_Y$ provided by Lemma~\ref{XG=FullCorner}. We define $\Ind_{\fY}$ similarly as $\Ind_{\fX}$. Now we have two ways to go from covariant representations of $G \curvearrowright X$ to covariant representations of $H \curvearrowright Y$: $\Ind_{\Phi^{-1}}$ introduced in Remark~\ref{Rem:Ind-Phi}, and $\Ind_{\fY}^{-1} \Ind_{\fX}$. It turns out that they coincide.
\bprop
\label{Ind=Ind}
In the situation described above, $\Ind_{\fY}^{-1} \Ind_{\fX} (\pi,\sigma)$ is unitarily equivalent to $\Ind_{\Phi^{-1}}(\pi,\sigma)$ for every covariant representation $(\pi,\sigma)$ of $G \curvearrowright X$.
\eprop
\bproof
Let $\Ind_{\Phi^{-1}}(\pi,\sigma) = (\rho,\tau)$, and let $\Ind_{\fY}^{-1} \Ind_{\fX} (\pi,\sigma) = (\rho',\tau')$. Let $i^X: \: C^*((X \rtimes G) \vert A) \into C(X) \rtimes G$ and $i^Y: \: C^*((Y \rtimes H) \vert A) \into C(Y) \rtimes H$ be the canonical embeddings. Also, let $i_X: \: C(X) \rtimes G \into C_0(\Omega) \rtimes (G \times H)$ and $i_Y: \: C(Y) \rtimes H \into C_0(\Omega) \rtimes (G \times H)$ be the embeddings obtained with the help of Lemma~\ref{XG=FullCorner}. Then $(\rho,\tau)$ is uniquely determined by $(\pi \rtimes \sigma) \circ i^X \circ \Phi^{-1} \sim_u (\rho \rtimes \tau) \circ i^Y$. We want to show that $\rho' \rtimes \tau'$ has the same property. $(\rho',\tau')$ is uniquely determined by the existence of a representation $\Pi$ of $C_0(\Omega) \rtimes (G \times H)$ with $\Pi \circ i_X \sim_u \pi \rtimes \sigma$ and $\Pi \circ i_Y \sim_u \rho' \rtimes \tau'$. Hence $(\rho' \rtimes \tau') \circ i^Y \sim_u \Pi \circ i_Y \circ i^Y$. On the groupoid level, $i_Y \circ i^Y$ is given by
$$
  Y \rtimes H \vert B \to Y \rtimes H \to \Omega \rtimes (G \times H), \, (y,h) \ma (y,\beta(y,h),h),
$$
where $\beta$ is defined in \S~\ref{tc->coc}. At the same time, $i_X \circ i^X \circ \Phi^{-1}$ on the groupoid level is given by
\begin{align*}
  & (Y \rtimes H) \vert B \to (X \rtimes G) \vert A \to X \rtimes G \to \Omega \rtimes (G \times H)\\
  & (y,h) \ma (y,b(h^{-1},y)^{-1}) \ma (y,b(h^{-1},y)^{-1},\alpha(b(h^{-1},y),y)^{-1}),
\end{align*}
where $b$ comes from the groupoid isomorphism $(X \rtimes G) \vert A \cong (Y \rtimes H) \vert A$ (see Remark~\ref{Kaku=coe} and \cite[Definition~2.6]{Li2}) and $\alpha$ is defined in \S~\ref{tc->coc}. We have $\alpha(b(h^{-1},y),y) = h^{-1}$ by \cite[Lemma~2.10]{Li1} (or rather its analogue for partial actions). Hence $i_Y \circ i^Y = i_X \circ i^X \circ \Phi^{-1}$, so that $(\rho' \rtimes \tau') \circ i^Y \sim_u \Pi \circ i_Y \circ i^Y = \Pi \circ i_X \circ i^X \circ \Phi^{-1} \sim_u (\pi \rtimes \sigma) \circ i^X \circ \Phi^{-1}$. Our claim follows.
\eproof
Let $G \curvearrowright {}_Y \Omega {}_X \curvearrowleft H$ and $G \curvearrowright X$, $H \curvearrowright Y$ be as above. Let $A \subseteq X \cap Y$ be a clopen subspace with $G.A = X$, $H.A = Y$ and $(X \rtimes G) \vert A \cong (Y \rtimes H) \vert A$. Let $\Phi: \: C^*((X \rtimes G) \vert A) \cong C^*((Y \rtimes H) \vert A)$ be the induced C*-isomorphism. Let $\Pi = (\Pi^X,\Pi^G)$ be a covariant representation of $G \curvearrowright X$ on the Hilbert space $\cH$. Let $\sigma$ be a unitary representation of $G$ on $\cH_{\sigma}$. It is clear that $(1 \otimes \Pi^X, \sigma \otimes \Pi^G)$ is a covariant representation of $G \curvearrowright X$ on $\cH_{\sigma} \otimes \cH$. Let $\Ind_{\Phi^{-1}}(\sigma,\Pi)$ be the unitary representation of $H$ which is part of the covariant representation $\Ind_{\Phi^{-1}}(1 \otimes \Pi^X, \sigma \otimes \Pi^G)$. Moreover, let $\tau$ be a unitary representation of $H$ on $\cH_{\tau}$. Let $\Theta = (\Theta^Y,\Theta^H) = \Ind_{\Phi^{-1}}(\Pi^X,\Pi^G)$. Denote by $\Ind_{\Phi}(\Theta,\tau)$ the unitary representation of $G$ which is part of the covariant representation $\Ind_{\Phi}(\Theta^Y \otimes 1,\Theta^H \otimes \tau)$.

\blemma
\label{identify-covrep}
$(1 \otimes \Pi^X \otimes 1, \sigma \otimes \Ind_{\Phi}(\Theta,\tau)) = \Ind_{\Phi}(1 \otimes \Theta^Y \otimes 1, \Ind_{\Phi^{-1}}(\sigma,\Pi) \otimes \tau)$.
\elemma
\bproof
We have to show that
$$
  (1 \otimes \Pi^X \otimes 1) \rtimes (\sigma \otimes \Ind_{\Phi}(\Theta,\tau)) \vert_{C^*((X \rtimes G) \vert A)} = (1 \otimes \Theta^Y \otimes 1) \rtimes (\Ind_{\Phi^{-1}}(\sigma,\Pi) \otimes \tau) \vert_{C^*((Y \rtimes H) \vert B)} \circ \Phi.
$$
Fix $g \in G$ and $h \in H$. Let $f$ be the characteristic function of a compact subset of $(X \times \gekl{g}) \cap (X \rtimes G) \vert A$ whose image under $\chi$ lies in $(Y \times \gekl{h}) \cap (Y \rtimes H) \vert B$. It suffices to consider such $f$ as they span a dense subset in $C^*((X \rtimes G) \vert A)$. We have
\begin{align*}
  & (1 \otimes \Theta^Y \otimes 1) \rtimes (\Ind_{\Phi^{-1}}(\sigma,\Pi) \otimes \tau)(\Phi(f))
  = ((1 \otimes \Theta^Y) \rtimes \Ind_{\Phi^{-1}}(\sigma,\Pi))(\Phi(f)) \otimes \tau(h)\\
  &= ((1 \otimes \Pi^X) \rtimes (\sigma \otimes \Pi^G)(f)) \otimes \tau(h)
  = \sigma(g) \otimes \Pi(f) \otimes \tau(h)
  = \sigma(g) \otimes (\Theta(\Phi(f)) \otimes \tau(h))\\
  &= \sigma(g) \otimes ((\Pi^X \otimes 1) \rtimes \Ind_{\Phi}(\Theta,\tau)(f))
  = (1 \otimes \Pi^X \otimes 1) \rtimes (\sigma \otimes \Ind_{\Phi}(\Theta,\tau))(f). \qedhere
\end{align*}
\eproof

Let $\Lambda$ be a representation of $C(X) \rtimes G$, and set $\ti{\Lambda} \defeq \Ind_{\fX} \Lambda$. Let
$$
  \cH_{\ti \Lambda,c} \defeq \menge{\eta \in \cH_{\ti \Lambda}}{\eta = \ti \Lambda(1_K) \eta \ \text{for some compact} \ K \subseteq \Omega},
$$
and let $\cL$ be the complex vector space of linear maps $\cH_{\ti \Lambda,c} \to \Cz$ which are bounded whenever restricted to a subspace of the form $\ti \Lambda(1_K) \cH_{\ti \Lambda}$, with $K \subseteq \Omega$ compact. Moreover, let $\Lambda^G$ be the unitary representation of $G$ on $\cH_{\Lambda}$ induced by $\Lambda$, and denote by ${\ti \Lambda}^G$ and ${\ti \Lambda}^H$ the unitary representations of $G$ and $H$ on $\cH_{\ti \Lambda}$ induced by $\ti \Lambda$. As $\cH_{\ti \Lambda,c}$ is obviously invariant under the $G$- and $H$-actions, we obtain by restriction $G$- and $H$-actions on $\cH_{\ti \Lambda,c}$. Finally, by dualizing, we obtain $G$- and $H$-actions on $\cL$.
\setlength{\parindent}{0cm} \setlength{\parskip}{0cm}
\blemma
\label{H=L^H}
There is a $G$-equivariant linear isomorphism $\cH_{\Lambda} \cong \cL^H$.
\elemma
\bproof
Up to unitary equivalence, we have $\cH_{\Lambda} = \ti \Lambda (1_{\bar X}) \cH_{\ti \Lambda}$, and $\Lambda^G$ is given by the composite
$$
  G \into C(X) \rtimes G \cong 1_{\bar X} C_0(\Omega) \rtimes (G \times H) 1_{\bar X} \overset{\ti \Lambda}{\lori} \cL \rukl{\ti{\Lambda}(1_{\bar X}) \cH_{\ti \Lambda}},
$$
where the first map is given by $G \into C(X) \rtimes G, \, g \ma u_g$.
\setlength{\parindent}{0.5cm} \setlength{\parskip}{0cm}

We define $L: \: \cH_{\Lambda} \to \cL$ by setting $L(\xi)(\eta) = \sum_{h \in H} \spkl{{\ti \Lambda}^H(h) \xi,\eta}$. Here $\spkl{\cdot,\cdot}$ is the inner product in $\cH_{\ti \Lambda}$, and our convention is that it is linear in the second component. Note that in the definition of $L(\xi)(\eta)$, the sum is always finite since $\eta$ lies in $\cH_{\ti \Lambda,c}$. It is clear that $L$ is linear. Moreover, we have
$$
  L(\xi)({\ti \Lambda}^H(h') \eta) = \sum_h \spkl{{\ti \Lambda}^H(h) \xi, {\ti \Lambda}^H(h') \eta} = \sum_h \spkl{{\ti \Lambda}^H((h')^{-1}h) \xi, \eta} = L(\xi)(\eta).
$$
Therefore, the image of $L$ lies in $\cL^H$, and we obtain a linear map $\cH_{\Lambda} \to \cL^H$. We claim that the inverse is given by $R: \: \cL^H \to \cH_{\Lambda}^* \cong \cH_{\Lambda}$, where the first map is given by restriction, $l \ma l \vert_{\ti \Lambda(1_{\bar X}) \cH_{\ti \Lambda}}$, and the second map is the canonical isomorphism, identifying $\zeta \in \cH_{\Lambda}$ with the element $\spkl{\zeta,\cdot} \in \cH_{\Lambda}^*$. Note that $l \vert_{\ti \Lambda(1_{\bar X}) \cH_{\ti \Lambda}}$ is bounded because of our definition of $\cH_{\ti \Lambda,c}$. Let us show that $R$ is the inverse of $L$. For $l \in \cL^H$, we have
\begin{align*}
  L(R(l))(\eta) &= \sum_h \spkl{{\ti \Lambda}^H(h) R(l),\eta} = \sum_h \spkl{R(l),{\ti \Lambda}^H(h^{-1}) \eta} 
  = \sum_h l(\ti \Lambda(1_{\bar X}) {\ti \Lambda}^H(h^{-1}) \eta)\\
  &= \sum_h l({\ti \Lambda}^H(h) \ti \Lambda(1_{\bar X}) {\ti \Lambda}^H(h^{-1}) \eta)
  = \sum_h l(\ti \Lambda(1_{\bar X \, h} \eta) = l(\eta).
\end{align*}
For $\xi \in \cH_{\Lambda} = \ti \Lambda (1_{\bar X}) \cH_{\ti \Lambda}$, we have $R(L(\xi)) = \xi$ since
\begin{align*}
  L(\xi)(\ti \Lambda(1_{\bar X}) \eta) = \sum_h \spkl{{\ti \Lambda}^H(h) \xi, \ti \Lambda(1_{\bar X}) \eta} = \sum_h \spkl{{\ti \Lambda}^H(h) \ti \Lambda(1_{\bar X \, h}) \xi, \eta} = \spkl{\xi,\eta}
\end{align*}
because 
$\ti \Lambda(1_{\bar X \, h}) \xi = \xi$ if $h = e$ and $\ti \Lambda(1_{\bar X \, h}) \xi = 0$ if $h \neq e$.
\setlength{\parindent}{0cm} \setlength{\parskip}{0cm}

Finally, let us show that $L$ is $G$-equivariant:
\begin{align*}
  L(\Lambda^G(g) \xi)(\eta) &= \sum_h \spkl{{\ti \Lambda}^H(h)(\Lambda^G(g) \xi),\eta}
  = \sum_h \sum_j \spkl{{\ti \Lambda}^H(h) {\ti \Lambda}(1_{g \, \bar X \, j^{-1} \cap \bar X}) {\ti \Lambda}^G(g) {\ti \Lambda}^H(j) \xi,\eta}\\
  &= \sum_{h,j} \spkl{{\ti \Lambda}^H(h) {\ti \Lambda}^H(j)^{-1} {\ti \Lambda}(1_{g \, \bar X \, j^{-1} \cap \bar X}) {\ti \Lambda}^H(j) {\ti \Lambda}^G(g) \xi,\eta}\\
  &= \sum_{h,j} \spkl{{\ti \Lambda}^H(h) {\ti \Lambda}(1_{g \, \bar X \cap \bar X \, j}) {\ti \Lambda}^G(g) \xi,\eta}
  = \sum_h \spkl{{\ti \Lambda}^H(h) {\ti \Lambda}(1_{g \, \bar X}) {\ti \Lambda}^G(g) \xi,\eta}\\
  &= \sum_h \spkl{{\ti \Lambda}^G(g) {\ti \Lambda}^H(h) {\ti \Lambda}(1_{\bar X}) \xi,\eta}
  = \sum_h \spkl{{\ti \Lambda}^H(h) \xi,{\ti \Lambda}^G(g)^{-1} \eta}
  = L(\xi)({\ti \Lambda}^G(g)^{-1} \eta). \qedhere
\end{align*}
\eproof
\setlength{\parindent}{0cm} \setlength{\parskip}{0cm}

\bcor
\label{inv=GxH-inv}
We have $\gekl{\Lambda^G \text{-invariant vectors}} = \cH_{\Lambda}^G \cong \cL^{G \times H}$.
\ecor
\btheo
There exists a one-to-one correspondence between $(\Ind_{\Phi^{-1}}(\sigma,\Pi) \otimes \tau)$-invariant vectors and $(\sigma \otimes \Ind_{\Phi}(\Theta,\tau))$-invariant vectors.
\etheo
\setlength{\parindent}{0cm} \setlength{\parskip}{0cm}
\bproof
Obviously, $(1 \otimes \Pi^X \otimes 1, \sigma \otimes \Ind_{\Phi}(\Theta,\tau))$ is a covariant representation of $G \curvearrowright X$. Let $\Lambda \defeq (1 \otimes \Pi \otimes 1) \rtimes (\sigma \otimes \Ind_{\Phi}(\Theta,\tau))$. Set $\ti{\Lambda} \defeq \Ind_{\fX} \Lambda$, and define $\cL$ as in Lemma~\ref{H=L^H}. Then Corollary~\ref{inv=GxH-inv} yields a one-to-one correspondence between $(\sigma \otimes \Ind_{\Phi}(\Theta,\tau))$-invariant vectors and $\cL^{G \times H}$.
\setlength{\parindent}{0.5cm} \setlength{\parskip}{0cm}

Let $\Ind_{\Phi^{-1}} \Lambda$ be the representation of $C(Y) \rtimes H$ corresponding to $\Ind_{\Phi^{-1}}(1 \otimes \Pi^X \otimes 1, \sigma \otimes \Ind_{\Phi}(\Theta,\tau))$. By Proposition~\ref{Ind=Ind}, $\Ind_{\fY} \Ind_{\Phi^{-1}} \Lambda \sim_u \ti{\Lambda}$. Hence, together with Lemma~\ref{identify-covrep}, Corollary~\ref{inv=GxH-inv} yields a one-to-one correspondence between $(\Ind_{\Phi^{-1}}(\sigma,\Pi) \otimes \tau)$-invariant vectors and $\cL^{G \times H}$.

Thus $\gekl{(\Ind_{\Phi^{-1}}(\sigma,\Pi) \otimes \tau) \text{-invariant vectors}}
  \overset{\text{1-1}}{\longleftrightarrow} \cL^{G \times H} \overset{\text{1-1}}{\longleftrightarrow}
  \gekl{(\sigma \otimes \Ind_{\Phi}(\Theta,\tau)) \text{-invariant vectors}}$.
\eproof
\setlength{\parindent}{0cm} \setlength{\parskip}{0.5cm}
\bcor
\label{ex-inv-vec}
$\Ind_{\Phi^{-1}}(\sigma,\Pi) \otimes \tau$ has an invariant vector if and only if $\sigma \otimes \Ind_{\Phi}(\Theta,\tau)$ has an invariant vector.
\ecor
\setlength{\parindent}{0cm} \setlength{\parskip}{0cm}

We now come to Shalom's property $H_{FD}$. Recall that a group $G$ has $H_{FD}$ if for every unitary representation $\sigma$ of $G$, $\bar{H}^1(G,\sigma) \ncong \gekl{0}$ implies that $\sigma$ contains a finite dimensional subrepresentation.

\bdefin
A $(G,H)$ continuous orbit couple is called $\bar{H}^1 G$-faithful if its $G$- and $H$-spaces are second countable compact, and its topological dynamical system $G \curvearrowright X$ has the property that for every unitary representation $\sigma$ of $G$ with $\bar{H}^1(G,\sigma) \ncong \gekl{0}$, there exists a covariant representation $(\Pi^X,\Pi^G)$ of $G \curvearrowright X$ such that $\bar{H}^1(G,\sigma \otimes \Pi^G) \ncong \gekl{0}$. 
\edefin
\btheo
\label{HFD-qi}
Let $G$, $H$ be countable discrete groups. Suppose there exists an $\bar{H}^1 G$-faithful topologically free $(G,H)$ continuous orbit couple. If $H$ has property $H_{FD}$, then $G$ has property $H_{FD}$.
\etheo
\setlength{\parindent}{0cm} \setlength{\parskip}{0cm}

For the proof, let us recall the following observation which is explained in \cite[\S~3.1]{Sha}:
\blemma
\label{fin-dim-subrep--inv-vec}
A unitary representation $\sigma$ of a countable discrete group $G$ contains a finite dimensional subrepresentation if and only if there is a unitary representation $\sigma'$ of $G$ such that $\sigma \otimes \sigma'$ has an invariant vector.
\elemma

\bproof[Proof of Theorem~\ref{HFD-qi}]
By Lemma~\ref{A=B}, we may assume that our $\bar{H}^1 G$-faithful topologically free $(G,H)$ continuous orbit couple corresponds to a topologically free $(G,H)$ coupling $G \curvearrowright {}_{Y}\Omega {}_{X} \curvearrowleft H$ with second countable compact spaces $X$ and $Y$, which leads to topological dynamical systems $G \curvearrowright X$ and $H \curvearrowright Y$ together with a clopen subspace $A \subseteq X \cap Y$ with $G.A = X$, $H.A = Y$ and $(X \rtimes G) \vert A \cong (Y \rtimes H) \vert A$. Now let $\sigma$ be a unitary representation of $G$ with $\bar{H}^1(G,\sigma) \ncong \gekl{0}$. By $\bar{H}^1 G$-faithfulness, there exists a covariant representation $(\Pi^X,\Pi^G)$ of $G \curvearrowright X$ with $\bar{H}^1(G,\sigma \otimes \Pi^G) \ncong \gekl{0}$. By Corollary~\ref{2nd,redcohom}, $\bar{H}^1(H, \Ind_{\Phi^{-1}}(\sigma,\Pi)) \cong \bar{H}^1(G,\sigma \otimes \Pi^G)$, so that $\bar{H}^1(H, \Ind_{\Phi^{-1}}(\sigma,\Pi)) \ncong \gekl{0}$. As $H$ has property $H_{FD}$, $\Ind_{\Phi^{-1}}(\sigma,\Pi)$ must have a finite dimensional subrepresentation. Thus Lemma~\ref{fin-dim-subrep--inv-vec} implies that there is a unitary representation $\tau$ of $H$ such that $\Ind_{\Phi^{-1}}(\sigma,\Pi) \otimes \tau$ has an invariant vector. By Corollary~\ref{ex-inv-vec}, $\sigma \otimes \Ind_{\Phi}(\Theta,\tau)$ must have an invariant vector. Again by Lemma~\ref{fin-dim-subrep--inv-vec}, this implies that $\sigma$ has a finite dimensional subrepresentation. Hence $G$ has property $H_{FD}$.
\eproof
\setlength{\parindent}{0cm} \setlength{\parskip}{0.5cm}

\bremark
\label{amenable-HFD}
A $(G,H)$ continuous orbit couple with second countable compact $G$- and $H$-spaces is $\bar{H}^1 G$-faithful if its topological dynamical system $G \curvearrowright X$ admits a $G$-invariant probability measure. To see this, let $\mu$ be such a measure. Let $(\Pi^X,\Pi^G)$ be the canonical covariant representation of $G \curvearrowright X$ on $L^2(\mu)$. Then $\Pi^G$ contains the trivial representation, so that $\sigma \otimes \Pi^G$ contains $\sigma$. This shows $\bar{H}^1 G$-faithfulness. In particular, this is the case when $G$ is amenable. Therefore, Theorem~\ref{u,q,l--DS} and Theorem~\ref{HFD-qi} imply \cite[Theorem~4.3.3]{Sha}. The case of amenable groups is not the only situation where invariant probability measures exist. It follows easily from \cite{Das} and Theorem~\ref{u,q,l--DS} that for residually finite groups $G$ and $H$ with coarsely equivalent box spaces, there exists a $(G,H)$ continuous orbit couple with second countable compact $G$- and $H$-spaces such that its topological dynamical system $G \curvearrowright X$ admits a $G$-invariant probability measure. A similar statement applies to sofic groups with coarsely equivalent spaces of graphs (see \cite{AF}).
\eremark

\section{Applications to (co)homology II}
\label{sec:HII}

We now turn to coarse invariants of (co)homological nature.
\setlength{\parindent}{0cm} \setlength{\parskip}{0cm}

\subsection{Coarse maps and res-invariant modules}
\label{ss:CM-resMod}

Let $G$ be a group, $R$ a commutative ring with unit and $W$ an $R$-module. Let $C(G,W)$ be the set of functions from $G$ to $W$. The $G$-action on itself by left multiplication induces a canonical left $RG$-module structure on $C(G,W)$. Explicitly, given $g \in G$ and $f \in C(G,W)$, $g.f$ is the element in $C(G,W)$ given by $(g.f)(x) = f(g^{-1}x)$ for all $x \in G$. We are interested in the following class of $RG$-submodules of $C(G,W)$. Given a subset $A$ of $G$, let $1_A$ be its indicator function, i.e., $1_A \in C(G,R)$ is given by $1_A(x) = 1$ if $x \in A$ and $1_A(x) = 0$ if $x \notin A$. Here $1$ is the unit of $R$. Given $f \in C(G,W)$ and $A \subseteq G$, we form the pointwise product $1_A \cdot f \in C(G,W)$. This is nothing else but the restriction of $f$ to $A$, extended by $0$ outside of $A$ to give a function $G \to W$.

\bdefin
An $RG$-submodule $L \subseteq C(G,W)$ is called $\res$-invariant if $1_A \cdot f$ lies in $L$ for all $f \in L$ and $A \subseteq G$.
\edefin

\bexs
\label{Examples}
For arbitrary $R$ and $W$, $C(G,W)$, $C_f(G,W) = \menge{f \in C(G,W)}{f \ {\rm takes} \ {\rm finitely} \ {\rm many} \ {\rm values}}$ and $RG \otimes_R W \cong \menge{f \in C(G,W)}{\supp(f) \ {\rm is} \ {\rm finite}}$ are $\res$-invariant.
\setlength{\parindent}{0.5cm} \setlength{\parskip}{0cm}

If $R = \Rz$ or $R = \Cz$, $W = R$, then $c_0 (G,W) = \menge{f \in C(G,W)}{\lim_{x \to \infty} \abs{f(x)} = 0}$ is $\res$-invariant, and for all $0 < p \leq \infty$, $\ell^p (G,W) = \menge{f \in C(G,W)}{\sum_{x \in G} \abs{f(x)}^p < \infty}$ is $\res$-invariant.

Let $G$ be a finitely generated discrete group and $\ell$ the right-invariant word length coming from a finite symmetric set of generators. Let $R = \Rz$ or $R = \Cz$ and $W = R$. As in \cite{Jol}, we define for $s \in \Rz$ and $1 \leq p \leq \infty$ the Sobolev space $H^{s,p}(G,W) \defeq \menge{f: \: G \to W}{f \cdot (1 + \ell)^s \in \ell^p (G,W)}$, and $H^{\infty,p}(G,W) \defeq \bigcap_{s \in \Rz} H^{s,p}(G,W)$. All these Sobolev spaces are $\res$-invariant.

In the last examples ($\ell^p$, $c_0$ and $H^{s,p}$), we can also replace $W$ by any normed space over $R$.
\eexs

We are also interested in the following topological setting: Let $R$ be a topological field and $W$ an $R$-module.
\bdefin
A topological $\res$-invariant $RG$-submodule $L$ of $C(G,W)$ is a $\res$-invariant $RG$-submodule of $C(G,W)$ together with the structure of a topological $R$-vector space on $L$ such that
\begin{align}
\label{res:cont}
  {\rm for} \ {\rm every} \ A \subseteq G, \ \ \ & L \to L, \, f \ma 1_A \cdot f \ {\rm is} \ {\rm continuous},\\
\label{g:cont}
  {\rm for} \ {\rm every} \ g \in G, \ \ \ & L \to L, \, f \ma g.f \ {\rm is} \ {\rm continuous}.
\end{align}
\edefin
\setlength{\parindent}{0cm} \setlength{\parskip}{0cm}
When we consider topological $\res$-invariant modules, $R$ will always be a topological field, though we might not mention this explicitly. For instance, in \ref{Examples}, $\ell^p(G,W)$ and $c_0(G,W)$ are topological $\res$-invariant modules. Also, $H^{s,p}(G,W)$ becomes a topological $\res$-invariant module with respect to the topology induced by the norm $\norm{f}_{s,p} = \norm{f \cdot (1 + \ell)^s}_{\ell^p (G,W)}$ for $s \in \Rz$, and with respect to the projective limit topology for $s = \infty$.

\setlength{\parindent}{0cm} \setlength{\parskip}{0.5cm}

In the following, we explain how coarse maps interact with $\res$-invariant modules. Recall that all our groups are countable and discrete, and that a map $\varphi: \: G \to H$ between groups $G$ and $H$ is a coarse map if for every $y \in H$, $\varphi^{-1}(\gekl{y})$ is finite, and for all $S \subseteq G \times G$, $\menge{\varphi(s) \varphi(t)^{-1}}{(s,t) \in S}$ must be finite if $\menge{st^{-1}}{(s,t) \in S}$ is finite (Definition~\ref{Def:CM}).
\bremark
\label{Rem:phi(g.)}
Let $\varphi: \: G \to H$ be a coarse map. Given $g \in G$, let $S = \menge{(g^{-1}x,x) \in G \times G}{x \in G}$. Then $\menge{st^{-1}}{(s,t) \in S} = \gekl{g^{-1}}$ is finite, so that $\menge{\varphi(g^{-1}x) \varphi(x)^{-1}}{x \in G}$ is finite. In other words, we can find a finite decomposition $G = \bigsqcup_{i \in I} X_i$, where $I$ is a finite index set, and a finite subset $\menge{h_i}{i \in I} \subseteq H$ such that $\varphi(g^{-1} x) = h_i^{-1} \varphi(x)$ for all $x \in X_i$ and $i \in I$.
\eremark
\setlength{\parindent}{0cm} \setlength{\parskip}{0cm}

Recall that two maps $\varphi, \, \phi: \: G \to H$ are close (written $\varphi \sim \phi$) if $\menge{\varphi(x) \phi(x)^{-1}}{x \in G}$ is finite (Definition~\ref{Def:CM}).
\setlength{\parindent}{0cm} \setlength{\parskip}{0.5cm}

\bremark
\label{Rem:sim}
If $\varphi, \, \phi: \: G \to H$ are close, then there is a finite decomposition $G = \bigsqcup_{i \in I} X_i$, where $I$ is a finite index set, and a finite subset $\menge{h_i}{i \in I} \subseteq H$ such that we have $\phi(x) = h_i \varphi(x)$ for all $x \in X_i$ and $i \in I$.
\eremark
\setlength{\parindent}{0cm} \setlength{\parskip}{0cm}

Let $R$ and $W$ be as above, and $\varphi: \: G \to H$ a coarse map. Given $f \in C(G,W)$, define $\varphi_*(f) \in C(H,W)$ by setting $\varphi_*(f)(y) = \sum_{\substack{x \in G \\ \varphi(x) = y}} f(x)$. Moreover, given $f \in C(H,W)$, define $\varphi^*(f) = f \circ \varphi \in C(G,W)$.
\bdefin
Given a $\res$-invariant $RG$-submodule $L$ of $C(G,W)$, let $\varphi_* L$ be the smallest $\res$-invariant $RH$-submodule of $C(H,W)$ containing $\menge{\varphi_*(f)}{f \in L}$. Given a $\res$-invariant $RH$-submodule $M$ of $C(G,W)$, let $\varphi^* M$ be the smallest $\res$-invariant $RG$-submodule of $C(H,W)$ containing $\menge{\varphi^*(f)}{f \in M}$.
\edefin
\blemma
We have
\begin{align}
\label{phi_L=}
  \varphi_* L &= \spkl{\menge{h.\varphi_*(f)}{h \in H, \, f \in L}}_R\\
\label{phi^M=}
  \varphi^* M &= \spkl{\menge{1_A \cdot \varphi^*(f)}{f \in M, \, A \subseteq G}}_R.
\end{align}
\elemma
\bproof
We obviously have \an{$\supseteq$} in \eqref{phi_L=}. To show \an{$\subseteq$}, it suffices to show that the right-hand side is $\res$-invariant as it is obviously an $RH$-submodule. Given $B \subseteq H$, we have for all $h \in H$ and $f \in L$ that
$$
   1_B \cdot (h.\varphi_*(f)) = h.\rukl{1_{h^{-1}B} \cdot \varphi_*(f)} = h. \rukl{\varphi_* \rukl{1_{\varphi^{-1}(h^{-1}B)} \cdot f}},
$$
which lies in the right-hand side as $L$ is $\res$-invariant.

For \eqref{phi^M=}, we again have \an{$\supseteq$} by construction. As the right-hand side is $\res$-invariant, it suffices to show that it is an $RG$-submodule in order to prove \an{$\subseteq$}. Given $g \in G$, by Remark~\ref{Rem:phi(g.)} we can find a finite decomposition $G = \bigsqcup_{i \in I} X_i$ and a finite subset $\menge{h_i}{i \in I} \subseteq H$ such that $\varphi(g^{-1}x) = h_i^{-1} \varphi(x)$ for all $x \in X_i$ and $i \in I$. Then, for all $A \subseteq G$, $g. \rukl{1_A \cdot \varphi^*(f)} = 1_{gA} \cdot \rukl{g. \varphi^*(f)} = \sum_{i \in I} 1_{X_i} \cdot 1_{gA} \cdot \rukl{g. \varphi^*(f)} = \sum_{i \in I} 1_{X_i} \cdot 1_{gA} \cdot \rukl{\varphi^* \rukl{h_i.f}}$ lies in the right-hand side of \eqref{phi^M=} as $M$ is an $RH$-submodule.
\eproof
\setlength{\parindent}{0cm} \setlength{\parskip}{0cm}

Note that in general, $\varphi_* L$ is not equal to $\menge{\varphi_*(f)}{f \in L}$, and $\varphi^* M$ is not equal to $\menge{\varphi^*(f)}{f \in M}$.

\blemma
\label{alg:sim,circ}
(i) If $\varphi, \, \phi: \: G \to H$ are coarse maps with $\varphi \sim \phi$, then $\varphi_* L = \phi_* L$ and $\varphi^* M = \phi^* M$ for all $L$, $M$.
\setlength{\parindent}{0.5cm} \setlength{\parskip}{0cm}

(ii) $\psi_* \varphi_* L = (\psi \circ \varphi)_* L$ and $\varphi^* \psi^* N = (\psi \circ \varphi)^* N$ for all $L$, $N$ and coarse maps $\varphi: \: G \to H$, $\psi: \: H \to K$.
\elemma
\bproof
(i) Let us show $\varphi_* L = \phi_* L$. By Remark~\ref{Rem:sim}, there is a finite decomposition $G = \bigsqcup_{i \in I} X_i$ and a finite subset $\menge{h_i}{i \in I} \subseteq H$ such that $\phi(x) = h_i \varphi(x)$ for all $x \in X_i$ and $i \in I$. Then
$$
  \phi_*(f) = \sum_{i \in I} \phi_*(1_{X_i} \cdot f) = \sum_{i \in I} h_i. \varphi_*(1_{X_i} \cdot f) \in \varphi_* L
$$
for all $f \in L$. Hence $\phi_* L \subseteq \varphi_* L$. By symmetry, we have $\phi_* L = \varphi_* L$.
\setlength{\parindent}{0.5cm} \setlength{\parskip}{0cm}

Let us show $\varphi^* M = \phi^* M$. Let $I$, $\menge{X_i}{i \in I}$ and $\menge{h_i}{i \in I}$ be as above. We have that
$$
  \varphi^*(f) = \sum_{i \in I} 1_{X_i} \cdot \varphi^*(f) = \sum_{i \in I} 1_{X_i} \cdot \phi^*(h_i.f) \in \phi^* M
$$
for all $f \in M$. Hence $\varphi^* M \subseteq \phi^* M$. By symmetry, we have $\varphi^* M = \phi^* M$.
\setlength{\parindent}{0cm} \setlength{\parskip}{0.25cm}

(ii) Let us show $\psi_* \varphi_* L = (\psi \circ \varphi)_* L$. Obviously, \an{$\supseteq$} holds as $\psi_*\varphi_* L \ni \psi_*(\varphi_*(f)) = (\psi \circ \varphi)_*(f)$ for all $f \in L$. Let us show \an{$\subseteq$}. By \eqref{phi_L=}, it suffices to show that $\psi_*(h.\varphi_*(f)) \in (\psi \circ \phi)_* L$ for all $h \in H$ and $f \in L$. By Remark~\ref{Rem:phi(g.)}, we can find a finite decomposition $H = \bigsqcup_{i \in I} Y_i$ and a finite subset $\menge{k_i}{i \in I} \subseteq K$ such that $\psi(h^{-1}y) = k_i^{-1} \psi(y)$ for all $y \in Y_i$ and $i \in I$. Then
$$
  \psi_*(h. \varphi_*(f)) = \sum_{i \in I} \psi_* \rukl{1_{Y_i} \cdot \rukl{h. \varphi_*(f)}} = \sum_{i \in I} k_i. \psi_* \rukl{1_{h^{-1} Y_i} \cdot \rukl{\varphi_*(f)}} = \sum_{i \in I} k_i. (\psi \circ \varphi)_* \rukl{1_{\varphi^{-1}(h^{-1} Y_i)} \cdot f}
$$
lies in $(\psi \circ \varphi)_* L$ for all $f \in L$ as $L$ is $\res$-invariant. This shows \an{$\subseteq$}.
\setlength{\parindent}{0.5cm} \setlength{\parskip}{0cm}

Let us show $\varphi^*\psi^* N = (\psi \circ \varphi)^* N$. \an{$\supseteq$} holds as $\varphi^* \psi^* N \ni \varphi^*(\psi^*(f))$ for all $f \in N$. Let us prove \an{$\subseteq$}. By \eqref{phi^M=}, it suffices to prove that $\varphi^*(1_B \cdot \psi^*(f)) \in (\psi \circ \varphi)^* N$ for all $B \subseteq H$ and $f \in N$. We have
$$
  \varphi^*(1_B \cdot \psi^*(f)) = 1_{\varphi^{-1}(B)} \cdot \varphi^*(\psi^*(f)) =  1_{\varphi^{-1}(B)} \cdot (\psi \circ \varphi)^*(f),
$$
which lies in $(\psi \circ \varphi)^* N$ as the latter is $\res$-invariant. This shows \an{$\subseteq$}.
\eproof

\subsection{Coarse embeddings and $\res$-invariant modules}
\label{ss:CE-resMod}

Recall that a map $\varphi: \: G \to H$ between groups $G$ and $H$ is a coarse embedding if for every $S \subseteq G \times G$, $\menge{st^{-1}}{(s,t) \in S}$ is finite if and only if $\menge{\varphi(s) \varphi(t)^{-1}}{(s,t) \in S}$ is finite  (Definition~\ref{Def:CM}).

\blemma
\label{Lem:dec-phi}
Let $\varphi: \: G \to H$ be a coarse embedding, and let $Y \defeq \varphi(G)$. Then we can find $X \subseteq G$ such that $X \to Y, \, x \ma \varphi(x)$, is a bijection. In addition, we can find a finite decomposition $G = \bigsqcup_{i=1}^I X_i$, $g(i) \in G$ for $1 \leq i \leq I$ and $h(i) \in H$ for $1 \leq i \leq I$, such that $X_i = g(i)^{-1} X(i)$ for some $X(i) \subseteq X$, with $g(1) = e$ (identity in $G$), $h(1) = e$ (identity in $H$), $X_1 = X(1) = X$, and $\varphi(x) = h(i) \varphi(g(i)x)$ for all $x \in X_i$ and $1 \leq i \leq I$.
\elemma
\setlength{\parindent}{0cm} \setlength{\parskip}{0cm}

\bproof
By Lemma~\ref{Lem:UE:X->Y}, we can find $X$ such that the restriction of $\varphi$ to $X$ is bijective onto its image and that there are finitely many $g(i) \in G$, $1 \leq i \leq I$, such that $G = \bigcup_{i=1}^I g(i)^{-1} X$, where we can certainly arrange $g(1) = e$. Now define recursively $X_1 \defeq X$ and $X(i) = X \setminus g(i) \rukl{g(1)^{-1} X_1 \cup \dotso \cup g(i-1)^{-1} X_{i-1}}$. Then $G = \bigsqcup_{i=1}^I g(i)^{-1} X(i)$. Using Remark~\ref{Rem:phi(g.)}, we can further decompose each $X(i)$ to guarantee that there exist $h(i) \in H$ for $1 \leq i \leq I$ such that $\varphi(x) = h(i) \varphi(g(i) x)$ for all $x \in g(i)^{-1} X(i)$ and $1 \leq i \leq I$. Setting $X_i \defeq g(i)^{-1} X(i)$, we are done.
\eproof

\blemma
\label{Lem:phi=phi}
Let $\varphi: \: G \to H$ be a coarse embedding, and fix $h \in H$. There exists a finite subset $F \subseteq G$ such that for all $x, \ti{x} \in G$ with $\varphi(\ti{x}) = h^{-1} \varphi(x)$, we must have $\ti{x} \in F x$.
\elemma
\bproof
Let $S = \menge{(s,t) \in G}{\varphi(s) = h^{-1} \varphi(t)}$. Then $\menge{\varphi(s) \varphi(t)^{-1}}{(s,t) \in S} = \gekl{h^{-1}}$ is finite, so that $F = \menge{st^{-1}}{(s,t) \in S}$ is finite since $\varphi$ is a coarse embedding.
\eproof

Let $\varphi: \: G \to H$ be a coarse embedding, and set $Y \defeq \varphi(G)$. Lemma~\ref{Lem:UE:X->Y} yields a subset $X \subseteq G$ such that the restriction of $\varphi$ to $X$ is a bijection $\ti{\varphi}: \: X \cong Y, \, x \ma \varphi(x)$. It is clear that $H = \bigcup_{h \in H} hY$. Enumerate $H$, say $H = \gekl{h_1, h_2, \dotsc}$, where $h_1 = e$ is the identity. Define recursively $Y_1 \defeq Y$ and $Y_j \defeq Y \setminus h_j^{-1} \rukl{h_1 Y_1 \cup \dotso \cup h_{j-1} Y_{j-1}}$. By construction, we have a decomposition as a disjoint union $H = \bigsqcup_{j=1}^{\infty} h_j Y_j$. Clearly, for all $h \in H$,
\bgl
\label{hYhjYj}
  hY \cap h_j Y_j = \empty \ {\rm for} \ {\rm all} \ {\rm but} \ {\rm finitely} \ {\rm many} \ j.
\egl
\bdefin
\label{Def:omega}
Define $\omega: \: H \to G$ by setting $\omega(y) = \ti{\varphi}^{-1}(h_j^{-1}y)$ for $y \in h_j Y_j$.
\edefin
By construction,
\bgl
\label{phiomega}
  (\varphi \circ \omega)(y) = h_j^{-1} y \ {\rm for} \ y \in h_j Y_j.
\egl
Take $F$ as in Lemma~\ref{Lem:phi=phi} for $h = e$. $(\omega \circ \varphi)(x) \in F x$ for all $x \in G$, so $\menge{(\omega \circ \varphi)(x) x^{-1}}{x \in G}$ is finite, i.e.,
\bgl
\label{omegaphi}
  \omega \circ \varphi \sim \id_G.
\egl
In general, pre-images under $\omega$ can be infinite, so that for an arbitrary $f \in C(H,W)$, $\omega_*(f)$ may not be defined. However, we can define $\omega_*(f)$ for $f \in \varphi_*L$, where $L \subseteq C(G,W)$ is a $\res$-invariant $RG$-submodule. We need some preparation. The following is an immediate consequence of \eqref{phi_L=} and \eqref{hYhjYj}: 
\blemma
\label{Lem:phi_L=bigoplus}
We have $\varphi_* L = \bigoplus_{j=1}^{\infty} 1_{h_j Y_j} \cdot  (\varphi_* L)$ as $R$-modules.
\elemma
Let $F$ be as in Lemma~\ref{Lem:phi=phi} for $h = e$. For every $x \in G$, define $F_x \subseteq F$ by $\menge{\ti{x} \in G}{\varphi(\ti{x}) = \varphi(x)} = F_x x$. For every subset $F_i \subseteq F$, define $X_i = \menge{x \in G}{F_x = F_i}$. Then $G = \bigsqcup_{F_i \subseteq F} X_i$, and by construction, we have the following
\blemma
\label{Lem:phi^phi_=}
$\varphi^*(\varphi_*(f)) = \sum_{F_i \subseteq F} 1_{X_i} \cdot \rukl{\sum_{g \in F_i} g^{-1}.f}$.
\elemma
Similarly, let $F$ be as in Lemma~\ref{Lem:phi=phi} for some fixed $h \in H$. Let $X \subseteq G$ be as above. For all $x \in X$, define $F_x \subseteq F$ by setting $\menge{\ti{x} \in G}{\varphi(\ti{x}) = h^{-1} \varphi(x)} = F_x x$. For a subset $F_i \subseteq F$, let $X_i = \menge{x \in X}{F_x = F_i}$. We have $X = \bigsqcup_{F_i \subseteq F} X_i$ and, by construction,
\blemma
\label{Lem:hphi_=}
$1_Y \cdot (h.\varphi_*(f)) = \varphi_* \rukl{\sum_{F_i \subseteq F} 1_{X_i} \cdot \rukl{\sum_{g \in F_i} g^{-1}.f}}$.
\elemma
Now we are ready for the following
\blemma
Let $L \subseteq C(G,W)$ be an $\res$-invariant $RG$-submodule. Then $\varphi_* L \to L, \, f \ma \omega_*(f)$ is well-defined, where $\omega_*(f)(x) = \sum_{\substack{y \in H \\ \omega(y) = x}} f(y)$.
\elemma
\bproof
By Lemma~\ref{Lem:phi_L=bigoplus}, it suffices to show that for every $j$ and $f \in 1_{h_j Y_j} \cdot \varphi_* L$, $\omega_*(f)$ lies in $L$. For such $f$, we know that $\omega_*(f) = 1_X \cdot \varphi^*(h_j^{-1}.f)$. As $f$ lies in $1_{h_j Y_j} \cdot \varphi_* L$, $h_j^{-1}.f$ lies in $1_{Y_j} \cdot \varphi_* L \subseteq 1_Y \cdot \varphi_* L$. Hence it suffices to show that $1_X \cdot \varphi^*(\varphi_* L) \subseteq L$. By \eqref{phi_L=}, it is enough to show that $\varphi^*(h.\varphi_*(f)) \in L$ for all $f \in L$. This follows immediately from Lemma~\ref{Lem:phi^phi_=} and Lemma~\ref{Lem:hphi_=}.
\eproof
\bdefin
For $\varphi$ and $L$ as above, set $\varphi^{*-1} L \defeq \menge{f \in C(H,W)}{\varphi^*(h.f) \in L \ {\rm for} \ {\rm all} \ h \in H}$.
\edefin
We collect a few properties of $\varphi^{*-1} L$:
\blemma
\label{Lem:phi^-1}
\begin{enumerate}
\item[a)] $\varphi^{*-1} L$ is an $\res$-invariant $RH$-submodule of $C(H,W)$.
\item[b)] For $f \in C(H,W)$, $f \in \varphi^{*-1} L$ if and only if for all $h \in H$, $1_{hY} \cdot f \in \varphi^{*-1} L$.
\item[c)] $\varphi^{*-1} L$ is the biggest $\res$-invariant $RH$-submodule $M$ of $C(G,W)$ such that $\varphi^*(f) \in L$ for all $f \in M$.
\item[d)] Let $\omega$ be as in Definition~\ref{Def:omega}. Then $\omega^*(f) \in \varphi^{*-1} L$ for all $f \in L$.
\item[e)] $\varphi^* \varphi^{*-1} L = L$.
\end{enumerate}
\elemma
\bproof
a) $\varphi^{*-1} L$ is $H$-invariant by definition. To see that $\varphi^{*-1} L$ is $\res$-invariant, take $B \subseteq H$ and $f \in \varphi^{*-1} L$. Then, for all $h \in H$, $\varphi^*(h.(1_B \cdot f)) = \varphi^*(1_{hB} \cdot (h.f)) = 1_{\varphi^{-1}(hB)} \varphi^*(h.f) \in L$, so $1_B \cdot f \in \varphi^{*-1} L$.
\setlength{\parindent}{0.5cm} \setlength{\parskip}{0cm}

b) follows from $\varphi^*(h.f) = \varphi^*(1_Y \cdot (h.f)) = \varphi^*(h.(1_{h^{-1}Y} \cdot f))$ for all $f \in C(H,W)$.

c) If $M$ is an $\res$-invariant $RH$-submodule of $C(H,W)$, then $f \in M$ implies $h.f \in M$ for all $h \in H$, and hence, by b), we conclude that $f \in \varphi^{*-1} L$.

d) By b), it suffices to prove $1_{hY} \cdot \omega^*(f) \in \varphi^{*-1} L$ for all $h \in H$. By \eqref{hYhjYj}, it suffices to prove $1_{h_j Y_j} \cdot \omega^*(f) \in \varphi^{*-1} L$ for all $j$. For all $y \in h_j Y_j$, $1_{h_j Y_j} \cdot \omega^*(f)(y) = f(\omega(y)) = f(\ti{\varphi}^{-1}(h_j^{-1}y)) = \varphi_*(1_X \cdot f)(h_j^{-1}y)$, hence $1_{h_j Y_j} \cdot \omega^*(f) = h_j.\varphi_*(1_X \cdot f)$. Let $h \in H$ be arbitrary. Lemma~\ref{Lem:hphi_=} and Lemma~\ref{Lem:phi^phi_=} imply that $\varphi^*(h h_j.\varphi_*(1_X \cdot f))$ lies in $L$. Hence $\omega^*(f)$ lies in $\varphi^{*-1} L$.

e) We have $\varphi^*(f) \subseteq L$ for all $f \in \varphi^{*-1} L$ by construction (see also c)). Hence $\varphi^* \varphi^{*-1} L \subseteq L$ by minimality of $\varphi^* \varphi^{*-1} L$. To show $L \subseteq \varphi^* \varphi^{*-1} L$, it suffices to show that $1_X \cdot L \subseteq \varphi^* \varphi^{*-1} L$ as $L = \sum_j g(i)^{-1}.(1_X \cdot L)$ by Lemma~\ref{Lem:dec-phi}. Let $f \in 1_X \cdot L$. Then $\omega^*(f) \in \varphi^{*-1} L$ by d), and $1_X \cdot \varphi^*(\omega^*(f)) \in \varphi^* \varphi^{*-1} L$. But we have $1_X \cdot \varphi^*(\omega^*(f)) = 1_X \cdot (\omega \circ \varphi)^*(f) = 1_X \cdot f = f$ as $\omega \circ \varphi = \id$ on $X$.
\eproof
\blemma
\label{Lem:phi^-1:sim,circ}
If $\varphi, \, \phi: \: G \to H$ are coarse embeddings with $\varphi \sim \phi$, then $\varphi^{*-1} L = \phi^{*-1} L$.

If $\varphi: \: G \to H$, $\psi: \: H \to K$ are coarse embeddings, then $\psi^{*-1} \varphi^{*-1} L = (\psi \circ \varphi)^{*-1} L$.
\elemma
\bproof
By Remark~\ref{Rem:sim}, we have for $f \in C(H,W)$: $\varphi^*(f) = \sum_i 1_{X_i} \cdot \phi^*(h_i.f)$. Hence $\phi^{*-1} L \subseteq \varphi^{*-1} L$. By symmetry, $\phi^{*-1} L = \varphi^{*-1} L$.
\setlength{\parindent}{0.5cm} \setlength{\parskip}{0cm}

If $f \in \psi^{*-1} \varphi^{*-1} L$, then $\psi^*(f) \in \varphi^{*-1} L$, and thus $(\psi \circ \varphi)^*(f) = \varphi^*(\psi^*(f)) \in L$. Lemma~\ref{Lem:phi^-1}~c) implies $f \in (\psi \circ \varphi)^{*-1} L$. To show $(\psi \circ \varphi)^{*-1} L \subseteq \psi^{*-1} \varphi^{*-1} L$, take $f \in (\psi \circ \varphi)^{*-1} L$. To show $f \in \psi^{*-1} \varphi^{*-1} L$, it suffices to show for all $k \in K$ and $h \in H$ that $\varphi^*(h.\psi^*(k.f)) \in L$. By Remark~\ref{Rem:phi(g.)}, we have $\psi(h^{-1}y) = k_j^{-1} \psi(y)$ for all $y \in Y_j$ and $j \in J$, for suitable $J$, $Y_j$ and $k_j$, so that $\varphi^*(h.\psi^*(k.f)) = \varphi^* \rukl{\sum_j 1_{Y_j} \cdot \psi^*(k_jk.f)} = \sum_j 1_{\varphi^{-1}(Y_j)} (\psi \circ \varphi)^*(k_jk.f)$, which lies in $L$ as $f$ lies in $(\psi \circ \varphi)^{*-1} L$.
\eproof

Our next goal is to define a suitable topology on $\varphi_* L$ in case $L$ is a topological $\res$-invariant $RG$-submodule of $C(G,W)$ and $\varphi$ is a coarse embedding. We start with some preparations.
\blemma
Let $\ti{Y} \subset Y$ and $\ti{X} = X \cap \varphi^{-1}(\ti{Y})$. Then $1_{\ti{X}} \cdot L \to 1_{\ti{Y}} \cdot (\varphi_* L), \, f \ma \varphi_*(f)$ is bijective.
\elemma
\bproof
Injectivity holds as we can recover $f$ from $\varphi_*(f)$ using
$$
  \varphi^*(\varphi_*(f))(\ti{x}) = \varphi_*(f)(\varphi(\ti{x})) = \sum_{\substack{x \in G \\ \varphi(x) = \varphi(\ti{x})}} f(x) = f(\ti{x})
$$
for $f \in 1_{\ti{X}} \cdot L$ and $\ti{x} \in \ti{X}$. For surjectivity, \eqref{phi_L=} implies that it suffices to show that for all $h \in H$ and $f \in L$, $1_{\ti{Y}} \cdot \rukl{h. \varphi_*(f)}$ lies in the image of our map. This follows immediately from Lemma~\ref{Lem:hphi_=}.
\eproof
For $j \in \Zz$, $j \geq 1$, set $X_j \defeq X \cap \varphi^{-1}(Y_j)$. Obviously, for all $j \geq 1$, we have $1_{h_j Y_j} \cdot (\varphi_* L) = h_j. \rukl{1_{Y_j} \cdot (\varphi_* L)}$. Thus $1_{X_j} \cdot L \to 1_{h_j Y_j} \cdot (\varphi_* L), \, f \ma h_j. \varphi_*(f)$ is an isomorphism. For $J \in \Zz$, $J \geq 1$, define
$$
  \Phi^J: \: \bigoplus_{j=1}^J 1_{X_j} \cdot L \to \varphi_* L, \, (f_j)_j \ma \sum_{j=1}^J h_j.\varphi_*(f_j).
$$
\bdefin
Let $L$ be a topological $\res$-invariant $RG$-submodule of $C(G,W)$. Let $\tau$ be the finest topology on $\varphi_* L$ such that for all $J \in \Zz$, $J \geq 1$, $\Phi^J$ is continuous. Here $1_{X_j} \cdot L$ is given the subspace topology from $L$, and $\bigoplus_{j=1}^J 1_{X_j} \cdot L$ is given the product topology.
\edefin

The proof of the following lemma is straightforward.
\blemma
$\tau$ is the finest topology on $\varphi_* L$ satisfying the following properties:
\setlength{\parindent}{0cm} \setlength{\parskip}{0cm}

\begin{enumerate}
\item[(T$_1$)] $(\varphi_* L, \tau)$ is a topological $\res$-invariant $RH$-submodule of $C(H,W)$.
\item[(T$_2$)] $L \to (\varphi_* L, \tau), \, f \ma \varphi_*(f)$ is continuous.
\end{enumerate}
\elemma
\blemma
\label{Lem:omega_--cont}
Let $\omega$ be as in Definition~\ref{Def:omega}. Then $\omega_*: \: \varphi_* L \to L$ is continuous.
\elemma
\bproof
By definition of the topology of $\varphi_* L$, it suffices to show that for every $j$, $1_{X_j} \cdot L \to 1_{h_j Y_j} \cdot (\varphi_* L), \, f \ma \omega_*(h_j.\varphi_*(f))$ is continuous. But it is easy to see that for $f \in 1_{X_j} \cdot L$, $\omega_*(h_j.\varphi_*(f)) = \varphi^*(\varphi_*(f))$. Continuity now follows from Lemma~\ref{Lem:phi^phi_=}.
\eproof

Now let us define a suitable topology on $\varphi^* M$ in case $M$ is a topological $\res$-invariant $RH$-submodule of $C(H,W)$ and $\varphi$ is a coarse embedding. Again, some preparations are necessary. Let $\varphi: \: G \to H$ be a coarse embedding and $M$ a $\res$-invariant $RH$-submodule of $C(H,W)$.
\blemma
\label{M1_Y=}
Let $\ti{X} \subseteq G$ be such that the restriction of $\varphi$ to $\ti{X}$ is injective. Let $\ti{Y} \defeq \varphi(\ti{X})$. Then $1_{\ti{Y}} \cdot M \to 1_{\ti{X}} \cdot (\varphi^* M), \, f \ma 1_{\ti{X}} \cdot \varphi^*(f)$ is a bijection.
\elemma
\bproof
For every $f \in 1_{\ti{Y}} \cdot M$ and $y \in H$, we have
$$
  \varphi_*(1_{\ti{X}} \cdot \varphi^*(f)) (y)
  = \sum_{\substack{x \in \ti{X} \\ \varphi(x) = y}} \varphi^*(f)(x) 
  = \sum_{\substack{x \in \ti{X} \\ \varphi(x) = y}} (f)(\varphi(x)) = f(y).
$$
Hence $\varphi_*(1_{\ti{X}} \cdot \varphi^*(f)) = f$, and our map is injective. To show surjectivity, it suffices by \eqref{phi^M=} to show that for every $f \in M$ and $A \subseteq G$, $1_{\ti{X}} \cdot (1_A \cdot \varphi^*(f))$ lies in the image of our map. This follows from
$1_{\ti{X}} \cdot \rukl{1_A \cdot \varphi^*(f)} = 1_{A \cap \ti{X}} \cdot \varphi^*(f) =  1_{\ti{X}} \cdot \varphi^*(1_{\varphi(A \cap \ti{X})} \cdot f)$.
\eproof
Now let $Y = \varphi(G)$. Lemma~\ref{Lem:UE:X->Y} gives us $X \subseteq G$ such that $\varphi \vert_X$ is a bijection $X \cong Y, \, x \ma \varphi(x)$. By Lemma~\ref{Lem:dec-phi}, we can find a finite decomposition $G = \bigsqcup_{i=1}^I X_i$ and finite subsets $\menge{g(i)}{1 \leq i \leq I} \subseteq G$, $\menge{h(i)}{i \leq i \leq I} \subseteq H$ such that $X_i = g(i)^{-1} X(i)$ for some $X(i) \subseteq X$ and $\varphi(x) = h(i) \varphi(g(i)x)$ for all $x \in X_i$ and $1 \leq i \leq I$. Let $Y_i \defeq \varphi(X_i)$ and $\Phi: \: \bigoplus_{i=1}^I 1_{Y_i} \cdot M \to \varphi^* M, \, (f_i)_i \ma \sum_{i=1}^I 1_{X_i} \cdot \varphi^*(f_i)$. As we obviously have $\varphi^*M = \bigoplus_{i=1}^I 1_{X_i} \cdot (\varphi^* M)$, $\Phi$ is surjective. And by Lemma~\ref{M1_Y=}, $\Phi$ is injective. Thus $\Phi$ is an isomorphism of $R$-modules.

\bdefin
Let $M$ be a topological $\res$-invariant $RH$-submodule of $C(H,W)$. Define the topology $\tau$ on $\varphi^* M$ so that $\Phi$ becomes a homeomorphism. Here $1_{Y_i} \cdot M$ is given the subspace topology from $M$, and $\bigoplus_{i=1}^I 1_{Y_i} \cdot M$ is given the product topology. 
\edefin

The following lemma is straightforward to prove.
\blemma
\label{Lem:top--phi^}
$\tau$ is the finest topology on $\varphi^* M$ satisfying the following properties:
\setlength{\parindent}{0cm} \setlength{\parskip}{0cm}

\begin{enumerate}
\item[(T$^1$)] $(\varphi^* M, \tau)$ is a topological $\res$-invariant $RG$-submodule of $C(G,W)$.
\item[(T$^2$)] $M \to (\varphi^* M, \tau), \, f \ma \varphi^*(f)$ is continuous.
\end{enumerate}
\elemma

Now we define a suitable topology on $\varphi^{*-1} L$ for a topological $\res$-invariant $RG$-submodule $L$ of $C(G,W)$ and a coarse embedding $\varphi$. Lemma~\ref{Lem:phi^-1}~b) implies that $\varphi^{*-1} L = \prod_j 1_{h_j Y_j} \cdot (\varphi^{*-1} L)$. The following is easy to verify:
\blemma
For every $j$, $\Phi^{(j)}: \: 1_{X_j} \cdot L \to 1_{h_j Y_j} \cdot (\varphi^{*-1} L), \, f \ma h_j.\varphi_*(f)$ is a bijection whose inverse is given by $1_{h_j Y_j} \cdot (\varphi^{*-1} L) \to 1_{X_j} \cdot L, \, f \ma 1_{X_j} \cdot \varphi^*(h_j^{-1}.f)$. 
\elemma
\bdefin
Let $L$ be a topological $\res$-invariant $RG$-submodule of $C(G,W)$. Define the topology $\tau$ on $\varphi^{*-1} L$ so that $\prod_j \Phi^{(j)}: \: \prod_j 1_{X_j} \cdot L \to \prod_j 1_{h_j Y_j} \cdot (\varphi^{*-1} L) = \varphi^{*-1} L$ becomes a homeomorphism. Here $1_{X_j} \cdot L$ is given the subspace topology coming from $L$, and $\prod_j 1_{X_j} \cdot L$ is given the product topology. 
\edefin

The following is straightforward to prove:
\blemma
$\tau$ is the coarsest topology on $\varphi^{*-1} L$ satisfying the following properties:
\setlength{\parindent}{0cm} \setlength{\parskip}{0cm}

\begin{enumerate}
\item[(T$^{-1}$)] $(\varphi^{*-1} L, \tau)$ is a topological $\res$-invariant $RH$-submodule of $C(H,W)$.
\item[(T$^{-2}$)] $(\varphi^{*-1} L, \tau) \to L, \, f \ma \varphi^*(f)$ is continuous.
\end{enumerate}
\elemma

\blemma
Let $L$, $\varphi$, $\omega$ and $\varphi^{*-1} L$ be as above. Then $\omega^*: \: L \to \varphi^{*-1} L$ is continuous. 
\elemma
\bproof
It suffices to show continuity of $L \to 1_{X_j} \cdot L, \, f \ma 1_{X_j} \cdot \varphi^*(h_j^{-1}.\omega^*(f))$ for all $j$. $1_{X_j} \cdot \varphi^*(h_j^{-1}.\omega^*(f)) = 1_{X_j} \cdot \varphi^*\rukl{h_j^{-1}.\rukl{1_{h_j Y_j} \cdot \omega^*(f)}} = 1_{X_j} \cdot \rukl{\varphi^* \varphi_*(1_{X_j} \cdot f)} = 1_{X_j} \cdot f$, which clearly depends continuously on $f$.
\eproof

\blemma
Let $L$, $\varphi$ and $\varphi^{*-1} L$ be as above. We have $\varphi^* \varphi^{*-1} L = L$ as topological $\res$-invariant modules.
\elemma
\bproof
Let $\tau$ be the topology of $L$ and $\ti{\tau}$ the topology of $\varphi^* \varphi^{*-1} L$. As $\varphi^*: \: \varphi^{*-1} L \to (L,\tau)$ is continuous by (T$^{-2}$), we must have $\tau \subseteq \ti{\tau}$ by Lemma~\ref{Lem:top--phi^}. To prove $\ti{\tau} \subseteq \tau$, we show that $\id: \: (L,\tau) \to (\varphi^* \varphi^{*-1} L, \ti{\tau})$ is continuous. By construction of $\ti{\tau}$ it suffices to show that $L \to 1_{Y_i} \cdot \varphi^{*-1} L, \, f \ma \varphi_*(1_{X_i} \cdot f)$ is continuous for all $i$. By construction of the topology on $\varphi^{*-1} L$, it is enough to show that $L \to L, \, f \ma 1_{X_j} \cdot \varphi^*(h_j^{-1}.(\varphi_*(1_{X_i} \cdot f)))$ is continuous. This now follows from Lemma~\ref{Lem:hphi_=} and Lemma~\ref{Lem:phi^phi_=}.
\eproof

We have the following topological analogue of Lemma~\ref{alg:sim,circ}, which is straightforward to prove.
\blemma
(i) If $\varphi, \, \phi: \: G \to H$ are coarse embeddings with $\varphi \sim \phi$, then $\varphi_* L = \phi_* L$, $\varphi^*M = \phi^* M$ and $\varphi^{*-1}L = \phi^{*-1} L$ as topological $\res$-invariant modules, for all topological $\res$-invariant $RG$-submodules $L$ of $C(G,W)$ and all topological $\res$-invariant $RH$-submodules $M$ of $C(H,W)$.
\setlength{\parindent}{0.5cm} \setlength{\parskip}{0cm}

(ii) If $\varphi: \: G \to H$ and $\psi: \: H \to K$ are coarse embeddings, then $\psi_* \varphi_* L = (\psi \circ \varphi)_* L$, $\varphi^* \psi^* N = (\psi \circ \varphi)^* N$ and $\psi^{*-1} \varphi^{*-1} L = (\psi \circ \varphi)^{*-1} L$ as topological $\res$-invariant modules, for all topological $\res$-invariant $RG$-submodules $L$ of $C(G,W)$ and all topological $\res$-invariant $RK$-submodules $N$ of $C(K,W)$.
\elemma

\subsection{Coarse maps and (co)homology}

Let us explain how coarse maps induce maps in group (co)homology. We first need to write group (co)homology in terms of groupoids.
\setlength{\parindent}{0cm} \setlength{\parskip}{0.5cm}

Let $G$ be a group, $R$ a commutative ring with unit, $L$ an $RG$-module. We write $g.f$ for the action of $g \in G$ on $f \in L$. We recall the chain and cochain complexes coming from the bar resolution (see \cite[Chapter~III, \S~1]{Bro}): Let $(C_*(L), \partial_*)$ be the chain complex $\dotso \overset{\partial_3}{\lori} C_2(L) \overset{\partial_2}{\lori} C_1(L) \overset{\partial_1}{\lori} C_0(L)$ with $C_0(L) = L$ and $C_n(L) = C_f(G^n,L) \cong R[G^n] \otimes_R L$, where $C_f$ stands for maps with finite support, and $\partial_n = \sum_{i=0}^n (-1)^i \partial_n^{(i)}$, where
\begin{align*}
  & \partial_n^{(0)}(f)(g_1, \dotsc, g_{n-1}) = \sum_{g_0 \in G} g_0^{-1}.f(g_0, g_1, \dotsc, g_{n-1}),\\
  &\partial_n^{(i)}(f)(g_1, \dotsc, g_{n-1}) = \sum_{\substack{g, \bar{g} \in G \\ g \bar{g} = g_i}} f(g_1, \dotsc, g_{i-1}, g, \bar{g}, g_{i+1}, \dotsc, g_{n-1}) \ {\rm for} \ 1 \leq i \leq n-1,\\
  & \partial_n^{(n)}(f)(g_1, \dotsc, g_{n-1}) = \sum_{g_n \in G} f(g_1, \dotsc, g_{n-1}, g_n).
\end{align*}
Let $(C^*(L), \partial^*)$ be the cochain complex $C^0(L) \overset{\partial^0}{\lori} C^1(L) \overset{\partial^1}{\lori} C^2(L) \overset{\partial^2}{\lori} \dotso$ where $C^0(L) = L$, $C^n(L) = C(G^n,L)$ for $n \geq 1$, and $\partial^n = \sum_{i=0}^{n+1} (-1)^i \partial^n_{(i)}$, with:
\begin{align*}
  & \partial^n_{(0)}(f)(g_0, \dotsc, g_n) = g_0.f(g_1, \dotsc, g_n),\\
  & \delta^n_{(i)}(f)(g_0, \dotsc, g_n) = f(g_0, \dotsc, g_{i-1}g_i, \dotsc, g_n) \ {\rm for} \ 1 \leq i \leq n,\\
  & \delta^n_{(n+1)}(f)(g_0, \dotsc, g_n) = f(g_0, \dotsc, g_{n-1}).
\end{align*}

Now let $W$ be an $R$-module and $L \subseteq C(G,W)$ be an $RG$-submodule. Consider the transformation groupoid $\cG \defeq G \rtimes G$ attached to the left multiplication action of $G$ on $G$. By definition, $\cG = \menge{(x,g)}{x \in G, \, g \in G}$, and the range and source maps are given by $r(x,g) = x$, $s(x,g) = g^{-1}x$, whereas the multiplication is given by $(x,g_1) (g_1^{-1}x,g_2) = (x,g_1g_2)$. Define $\sigma: \: \cG \to G, \, (x,g) \ma g$. Let $\cG^{(0)} = G$, and for $n \geq 1$, set
$$
  \cG^{(n)} \defeq \menge{(\gamma_1, \dotsc, \gamma_n) \in \cG^n}{s(\gamma_i) = r(\gamma_{i+1}) \ {\rm for} \ {\rm all} \ 1 \leq i \leq n-1},
$$
and define, for $n \geq 1$, $\sigma: \: \cG^{(n)} \to G^n$ as the restriction of $\sigma^n: \: \cG^n \to G^n$ to $\cG^{(n)}$.

Note that $\cG^{(n)} = \menge{((x_1,g_1), \dotsc, (x_n,g_n)) \in \cG^n}{g_i^{-1}x_i = x_{i+1}) \ {\rm for} \ {\rm all} \ 1 \leq i \leq n-1}$, so that we have a bijection
\begin{align}
\label{cG^=G^}
  \cG^{(n)} \cong G \times G^n, \, ((x_1,g_1), \dotsc, (x_n,g_n)) \ma (x_1, g_1, \dotsc, g_n).
\end{align}
This is because for $2 \leq i \leq n$, $x_i$ is determined by the equation $x_i = g_{i-1}^{-1} \dotsm g_1^{-1} x_1$. We will often use this identification of $\cG^{(n)}$ with $G \times G^n$ without explicitly mentioning it.

Now, given $f \in C(\cG^{(n)},W)$ and $\vecg \in G^n$, we view $f \vert_{\sigma^{-1}(\vecg)}$ as the map in $C(G,W)$ given by $x \ma f(x,\vecg)$. Set $\supp(f) \defeq \menge{\vecg \in G^n}{f \vert_{\sigma^{-1}(\vecg)} \neq 0}$.

Let us define a chain complex $(D_*(L), d_*)$ as follows: For $n = 0, 1, 2, \dotsc$, set
$$
  D_n(L) \defeq \menge{f \in C(\cG^{(n)},W)}{\supp(f) \ {\rm is} \ {\rm finite}, \, f \vert_{\sigma^{-1}(\vecg)} \in L \ {\rm for} \ {\rm all} \ \vecg \in G^n}.
$$
Moreover, for all $n \geq 1$, define maps $d_n: \: D_n(L) \to D_{n-1}(L)$ by setting $d_n = \sum_{i=0}^{n} (-1)^i d_n^{(i)}$ with $d_n^{(i)} = (\delta_n^{(i)})_*$, where $\delta_1^{(0)} = s$, $\delta_1^{(1)} = r$, and for $n \geq 2$,
\begin{align*}
  & \delta_n^{(0)}(\gamma_1, \dotsc, \gamma_n) = (\gamma_2, \dotsc, \gamma_n),\\
  & \delta_n^{(i)}(\gamma_1, \dotsc, \gamma_n) = (\gamma_1, \dotsc, \gamma_i \gamma_{i+1}, \dotsc, \gamma_n) \ {\rm for} \ 1 \leq i \leq n-1,\\
  & \delta_n^{(n)}(\gamma_1, \dotsc, \gamma_n) = (\gamma_1, \dotsc, \gamma_{n-1}).
\end{align*}
Here, we use the same notation as in \S~\ref{ss:CM-resMod}, i.e.,
$(\delta_n^{(i)})_*(f)(\veceta) = \sum_{\substack{\vecgamma \in \cG^{(n)} \\ \delta_n^{(i)}(\vecgamma) = \veceta}} f(\vecgamma)$.

Let us define a cochain complex $(D^*(L), d^*)$ by setting, for all $n = 0, 1, 2, \dotsc$,
$$
  D^n(L) \defeq \menge{f \in C(\cG^{(n)},W)}{f \vert_{\sigma^{-1}(\vecg)} \in L \ {\rm for} \ {\rm all} \ \vecg \in G^n}.
$$
Moreover, for all $n$, define maps $d^n: \: D^n(L) \to D^{n+1}(L)$ by setting $d^n = \sum_{i=0}^{n+1} (-1)^i d^n_{(i)}$, with $d^n_{(i)} = (\delta^n_{(i)})^*$ (as in \S~\ref{ss:CM-resMod}, $ (\delta^n_{(i)})^*(f) = f \circ \delta^n_{(i)}$), where $\delta^0_{(0)} = s$, $\delta^0_{(1)} = r$, and for all $n \geq 1$,
\begin{align*}
  & \delta^n_{(0)}(\gamma_0, \dotsc, \gamma_n) = (\gamma_1, \dotsc, \gamma_n),\\
  & \delta^n_{(i)}(\gamma_0, \dotsc, \gamma_n) = (\gamma_0, \dotsc, \gamma_{i-1}\gamma_i, \dotsc, \gamma_n) \ {\rm for} \ 1 \leq i \leq n,\\
  & \delta^n_{(n+1)}(f)(\gamma_0, \dotsc, \gamma_n) = (\gamma_0, \dotsc, \gamma_{n-1}).
\end{align*}
We are also interested in the topological setting, where we assume that $R$ is a topological field, $L \subseteq C(G,W)$ a $RG$-submodule together with the structure of a topological $R$-vector space such that the $G$-action $G \curvearrowright L$ is by homeomorphisms. Equip the above chain and cochain complexes $C_*(L)$ and $C^*(L)$ with the topologies of pointwise convergence. We also equip $D_*(L)$ and $D^*(L)$ with the topologies of pointwise convergence, i.e., $f_i \in C(\cG^{(n)},W)$ converges to $f \in C(\cG^{(n)},W)$ if and only if $\lim_i f_i \vert_{\sigma^{-1}(\vecg)} = f \vert_{\sigma^{-1}(\vecg)}$ in $L$ for all $\vecg \in G^n$.

The following is now immediate:
\blemma
\label{G-cG}
(i) We have isomorphisms $\chi_*$ of chain complexes and $\chi^*$ of cochain complexes given by $\chi_n: \: C_n(L) \to D_n(L), \, \chi_n(f)(x,\vecg) = f(\vecg)(x)$ and $\chi^n: \: C^n(L) \to D^n(L), \, \chi^n(f)(x,\vecg) = f(\vecg)(x)$.
\setlength{\parindent}{0.5cm} \setlength{\parskip}{0cm}

(ii) In the topological setting, $\chi_*$ and $\chi^*$ from (i) are topological isomorphisms.
\elemma
\setlength{\parindent}{0cm} \setlength{\parskip}{0cm}

By definition of group (co)homology, we have $H_n(G,L) = H_n(C_*(L))$ and $H^n(G,L) = H^n(C^*(L))$. By definition of reduced group (co)homology, we have ${\bar H}_n(G,L) = {\bar H}_n(C_*(L))$ and ${\bar H}^n(G,L) = {\bar H}^n(C^*(L))$ in the topological setting (recall that ${\bar H}_n(C_*(L)) = \ker(\partial_n) / \overline{\img(\partial_{n+1})}$ and ${\bar H}^n(C_*(L)) = \ker(\partial^n) / \overline{\img(\partial^{n-1})}$). Hence we obtain
\bcor
\label{Cor:G-cG}
(i) $\chi_*$ and $\chi^*$ from Lemma~\ref{G-cG} induce isomorphisms $H_n(\chi_*): \: H_n(G,L) \cong H_n(D_*(L))$ and $H^n(\chi^*): \: H^n(G,L) \cong H^n(D^*(L))$ for all $n$.
\setlength{\parindent}{0.5cm} \setlength{\parskip}{0cm}

(ii) In the topological setting, $\chi_*$ and $\chi^*$ from Lemma~\ref{G-cG} induce isomorphisms ${\bar H}_n(\chi_*): \: {\bar H}_n(G,L) \cong {\bar H}_n(D_*(L))$ and ${\bar H}^n(\chi^*): \: {\bar H}^n(G,L) \cong {\bar H}^n(D^*(L))$ for all $n$.
\ecor

In this groupoid picture of group (co)homology, let us now explain how coarse maps induce chain and cochain maps. Let $\varphi: \: G \to H$ be a coarse map. Let $\cG = G \rtimes G$ and $\cH = H \rtimes H$. Define $\varphi^1: \: \cG \to \cH, \, (x,g) \ma (\varphi(x), \varphi(x) \varphi(g^{-1}x)^{-1})$. It is easy to see that $\varphi^1$ is a groupoid homomorphism. This means that if $\gamma_1$ and $\gamma_2$ are composable, then so are $\varphi^1(\gamma_1)$ and $\varphi^1(\gamma_2)$, and we have $\varphi^1(\gamma_1 \gamma_2) = \varphi^1(\gamma_1) \varphi^1(\gamma_2)$. For all $n \geq 1$, define $\varphi^n: \: \cG^{(n)} \to \cH^{(n)}, \, (\gamma_1, \dotsc, \gamma_n) \ma (\varphi^1(\gamma_1), \dotsc, \varphi^1(\gamma_n))$. Moreover, if $\varphi: \: G \to H$ is a coarse embedding, let $\omega: \: H \to G$ be as above, and define $\omega^1: \: \cH \to \cG, \, (y,h) \ma (\omega(y), \omega(y) \omega(h^{-1}y)^{-1})$, and for all $n \geq 1$, define $\omega^n: \: \cH^{(n)} \to \cG^{(n)}, \, (\eta_1, \dotsc, \eta_n) \ma (\omega^1(\eta_1), \dotsc, \omega^1(\eta_n))$. Now let $L$ be a $\res$-invariant $RG$-submodule of $C(G,W)$. For $f \in D_n(L)$, consider $(\varphi^n)_*(f)(\veceta) = \sum_{\substack{\vecgamma \in \cG^{(n)} \\ \varphi^n(\vecgamma) = \veceta}} f(\vecgamma)$. In case $\varphi$ is a coarse embedding and $\omega$ is as above, set for $f \in D_n(\varphi_* L)$ $(\omega^n)_*(f)(\vecgamma) = \sum_{\substack{\veceta \in \cH^{(n)} \\ \omega^n(\veceta) = \vecgamma}} f(\veceta)$.
\blemma
\label{Lem:D_}
(i) Let $\varphi: \: G \to H$ be a coarse map. For all $n$, $D_n(\varphi): \: D_n(L) \to D_n(\varphi_* L), \, f \ma (\varphi^n)_*(f)$ is well-defined and gives rise to a chain map $D_*(\varphi): \: D_*(L) \to D_*(\varphi_* L)$. If $\psi: \: H \to K$ is another coarse map, then we have
\begin{align}
\label{D_(psiphi)}
  D_*(\psi \circ \varphi) = D_*(\psi) \circ D_*(\varphi).
\end{align}
If $L$ is a topological $\res$-invariant $RG$-submodule of $C(G,W)$ and $\varphi$ is a coarse embedding, then for all $n$, $D_n(\varphi)$ is continuous.

\setlength{\parindent}{0.5cm}
(ii) If $\varphi$ is a coarse embedding, then $D_n(\omega): \: D_n(\varphi_* L) \to D_n(L), \, f \ma (\omega^n)_*(f)$ is well-defined and gives rise to a chain map $D_*(\omega): \: D_*(\varphi_* L) \to D_*(L)$. If $L$ is a topological $\res$-invariant module, then $D_n(\omega)$ is continuous for all $n$.
\elemma
\setlength{\parindent}{0cm} \setlength{\parskip}{0cm}

Note that for \eqref{D_(psiphi)} to make sense, we implicitly use Lemma~\ref{alg:sim,circ}~(ii).
\bproof
(i) To show that $D_n(\varphi)$ is well-defined, we have to show that $(\varphi^n)_*(f) \in D_n(\varphi_* L)$ for all $f \in D_n(L)$. It suffices to treat the case that $\supp(f) = \gekl{\vecg}$ for a single $\vecg = (g_1, \dotsc, g_n) \in G^n$, as a general element in $D_n(L)$ is a finite sum of such $f$. Let us first show that $(\varphi^n)_*(f)$ has finite support. As $\varphi$ is a coarse map,
\begin{align}
\label{D_:F:finite}
  F \defeq \menge{\varphi(x) \varphi(g_i^{-1}x)^{-1}}{x \in G, \, 1 \leq i \leq n} \ {\rm is} \ {\rm finite}.
\end{align}
Clearly, $\supp \rukl{(\varphi^n)_*(f)} \subseteq F^n$. To show that for every $\vech = (h_1, \dotsc, h_n) \in H^n$, $(\varphi^n)_*(f) \vert_{\sigma^{-1}(\vech)}$ lies in $\varphi_* L$, define
$$
  A \defeq \menge{x \in G}{\varphi(g_{i-1}^{-1} \dotsm g_1^{-1}x) \varphi(g_i^{-1} \dotsm g_1^{-1}x)^{-1} = h_i \ {\rm for} \ {\rm all} \ 1 \leq i \leq n}.
$$
Then $\varphi^n(x,\vecg) \in \sigma^{-1}(\vech)$ if and only if $x \in A$. Hence
$$
  (\varphi^n)_*(f)(y,\vech) = \sum_{\substack{x \in A \\ \varphi^n(x,\vecg) = (y,\vech)}} f(x,\vecg) = \sum_{\substack{x \in A \\ \varphi(x) = y}} f(x,\vecg) = \varphi_* \rukl{1_A \cdot \rukl{f \vert_{\sigma^{-1}(\vecg)}}} (y),
$$
so that
\begin{align}
\label{D_(phi):formula}
  (\varphi^n)_*(f) \vert_{\sigma^{-1}(\vech)} = \varphi_*\rukl{1_A \cdot \rukl{f \vert_{\sigma^{-1}(\vecg)}}}.
\end{align}
As $f \vert_{\sigma^{-1}(\vecg)}$ lies in $L$, $L$ is $\res$-invariant and $\varphi_*(\ti{f}) \in \varphi_* L$ for all $\ti{f} \in L$, this shows that $(\varphi^n)_*(f) \vert_{\sigma^{-1}(\vech)} \in \varphi_* L$. Hence $D_n(\varphi)$ is well-defined for all $n$. $(D_n(\varphi))_n$ is a chain map because $\varphi^n$ is a groupoid homomorphism for all $n$. \eqref{D_(psiphi)} holds because we have $(\psi^n)_* \circ (\varphi^n)_* = ((\psi \circ \varphi)^n)_*$ for all $n$. \eqref{D_(phi):formula} shows continuity of $D_n(\varphi)$ for all $n$ as the right-hand side depends continuously on $f$. This is because $L$ satisfies \eqref{res:cont} and the topology on $\varphi_* L$ satisfies (T$_2$).
\setlength{\parindent}{0.5cm} \setlength{\parskip}{0cm}

(ii) To show that $D_n(\omega)$ is well-defined, take $f \in D_n(\varphi_* L)$. We may assume $\supp(f) = \gekl{\vech}$ for $\vech = (\bar{h}_1, \dotsc, \bar{h}_n)$ and $f \vert_{\sigma^{-1}(\vech)} \in 1_{hY} \cdot (\varphi_* L)$. By \eqref{hYhjYj}, $hY \cup \bar{h}_1^{-1} h Y \cup \dotso \cup \bar{h}_n^{-1} \dotsm \bar{h}_1^{-1} h Y \subseteq \bigcup_{j=1}^J h_j Y_j$ for some $J$. Thus, for all $y \in hY$ and $1 \leq i \leq n$, $\omega(\bar{h}_i^{-1} \dotsm \bar{h}_1^{-1} y) = \ti{\varphi}^{-1}(h_j^{-1} y)$ for some $1 \leq j \leq J$. Now consider $S = \menge{(\bar{h}_{i-1}^{-1} \dotsm \bar{h}_1^{-1} y, \bar{h}_i^{-1} \dotsm \bar{h}_1^{-1} y)}{y \in hY, \, 1 \leq i \leq n}$. $\menge{\varphi(\omega(s)) \varphi(\omega(t))^{-1}}{(s,t) \in S} \subseteq \menge{h_j^{-1} h_k}{1 \leq j, k \leq J}$ is finite, so that $F \defeq \menge{\omega(s) \omega(t)^{-1}}{(s,t) \in S}$ is finite as $\varphi$ is a coarse embedding. Hence $\supp((\omega^n)_*(f)) \subseteq F^n$. A similar formula as \eqref{D_(phi):formula} shows that $(\omega^n)_*$ is well-defined, and continuous in the topological setting.
\eproof
\setlength{\parindent}{0cm} \setlength{\parskip}{0cm}

Now let $M$ be a $\res$-invariant $RH$-submodule of $C(H,W)$. For $f \in D^n(M)$, consider $(\varphi^n)^*(f) = f \circ \varphi^n$. If $\varphi$ is a coarse embedding, $L$ an $\res$-invariant $RG$-submodule of $C(G,W)$, set for $f \in D^n(L)$: $(\omega^n)^*(f) = f \circ \omega^n$.
\blemma
\label{Lem:D^}
(i) Let $\varphi$ be a coarse map. For all $n$, $D^n(\varphi): \: D^n(M) \to D^n(\varphi^* M), \, f \ma (\varphi^n)^*(f)$ is well-defined and gives rise to a cochain map $D^*(\varphi): \: D^*(M) \to D^*(\varphi^* M)$. If $\psi: \: H \to K$ is another coarse map, we have
\begin{align}
\label{D^(psiphi)}
  D^*(\psi \circ \varphi) = D^*(\varphi) \circ D^*(\psi).
\end{align}
If $M$ is a topological $\res$-invariant $RH$-submodule of $C(H,W)$ and $\varphi$ is a coarse embedding, then $D^n(\varphi)$ is continuous for all $n$.

\setlength{\parindent}{0.5cm} \setlength{\parskip}{0cm}
(ii) If $\varphi$ is a coarse embedding, then $D^n(\omega): \: D^n(L) \to D^n(\varphi^{*-1} L), \, f \ma (\omega^n)^*(f)$ is well-defined and gives rise to a cochain map $D^*(\omega): \: D^*(L) \to D^*(\varphi^{*-1} L)$. If $L$ is a topological $\res$-invariant module, then $D^n(\omega)$ is continuous for all $n$.
\elemma
\setlength{\parindent}{0cm} \setlength{\parskip}{0cm}

For \eqref{D^(psiphi)} to make sense, we implicitly use (ii) in Lemma~\ref{alg:sim,circ}.

\bproof
(i) To show that $D^n(\varphi)$ is well-defined, we have to show that for all $f \in D^n(M)$, $(\varphi^n)^*(f) \in D^n(\varphi^* M)$, i.e., $(\varphi^n)^*(f) \vert_{\sigma^{-1}(\vecg)} \in \varphi^* M$ for all $\vecg = (g_1, \dotsc, g_n) \in G^n$. $F = \menge{\varphi(x) \varphi(g_i^{-1}x)^{-1}}{x \in G, \, 1 \leq i \leq n}$ is finite by \eqref{D_:F:finite}. We also know that $\varphi^n(x,\vecg) \in \sigma^{-1}(F^n)$ for all $x \in G$. For $\vech = (h_1, \dotsc, h_n) \in F^n$, let
$$
A_{\vech} \defeq \menge{x \in G}{\varphi(g_{i-1}^{-1} \dotsm g_1^{-1}x) \varphi(g_i^{-1} \dotsm g_1^{-1}x)^{-1} = h_i \ {\rm for} \ {\rm all} \ 1 \leq i \leq n}.
$$
Then $G = \bigsqcup_{\vech \in F^n} A_{\vech}$, and for $x \in A_{\vech}$, we have $\varphi^n(x,\vecg) = (\varphi(x),\vech)$. Hence
$$
  (\varphi^n)^*(f) \vert_{\sigma^{-1}(\vecg)} (x) = f(\varphi^n(x,\vecg)) = \sum_{\vech \in F^n} 1_{A_{\vech}}(x) \cdot \rukl{f \vert_{\sigma^{-1}(\vech)}}(\varphi(x)) ,
$$
and thus
\begin{align}
\label{D^(phi):formula}
  (\varphi^n)^*(f) \vert_{\sigma^{-1}(\vecg)} = \sum_{\vech \in F^n} 1_{A_{\vech}} \cdot \varphi^*\rukl{f \vert_{\sigma^{-1}(\vech)}}.
\end{align}
As $f \vert_{\sigma^{-1}(\vech)} \in M$, $\varphi^*(\ti{f}) \in \varphi^* M$ for all $\ti{f} \in M$ and $\varphi^* M$ is $\res$-invariant, this shows that $(\varphi^n)^*(f) \vert_{\sigma^{-1}(\vecg)} \in \varphi^* M$. Hence $D^n(\varphi)$ is well-defined for all $n$. $(D^n(\varphi))_n$ is a cochain map because $\varphi^n$ is a groupoid homomorphism for all $n$. \eqref{D^(psiphi)} holds because we have $(\varphi^n)^* \circ (\psi^n)^* = ((\psi \circ \varphi)^n)^*$ for all $n$. \eqref{D^(phi):formula} shows that $D^n(\varphi)$ is continuous for all $n$ as the right-hand side depends continuously on $f$ because the topology on $\varphi^* M$ satisfies (T$^1$) and (T$^2$).
\setlength{\parindent}{0.5cm} \setlength{\parskip}{0cm}

(ii) Given $f \in D^n(L)$ and $\vech = (\bar{h}_1, \dotsc, \bar{h}_n) \in H^n$, we show $(\omega^n)^*(f) \vert_{\sigma^{-1}(\vech)} \in \varphi^{*-1} L$. By Lemma~\ref{Lem:phi^-1}~b), it suffices to show $1_{hY} \cdot \rukl{(\omega^n)^*(f) \vert_{\sigma^{-1}(\vech)}} \in \varphi^{*-1} L$ for all $h \in H$. As we saw in the proof of Lemma~\ref{Lem:D_}~(ii), $F = \menge{\omega(\bar{h}_{i-1}^{-1} \dotsm \bar{h}_1^{-1} y) \omega(\bar{h}_i^{-1} \dotsm \bar{h}_1^{-1} y)^{-1}}{y \in hY, 1 \leq i \leq n}$ is finite. Thus $\omega^n(y,\vech) \in \sigma^{-1}(F^n)$ for all $y \in hY$. For $\vecg \in F^n$, let $B_{\vecg} = \menge{y \in hY}{\omega(\bar{h}_{i-1}^{-1} \dotsm \bar{h}_1^{-1} y) \omega(\bar{h}_i^{-1} \dotsm \bar{h}_1^{-1} y)^{-1} = g_i \ {\rm for} \ {\rm all} \ 1 \leq i \leq n}$. We then have $hY = \bigsqcup_{\vecg \in F^n} B_{\vecg}$, and for $y \in B_{\vecg}$, $\omega^n(y,\vech) = (\omega(y),\vecg)$, so that $1_{hY} \cdot \rukl{(\omega^n)^*(f) \vert_{\sigma^{-1}(\vech)}} = \sum_{\vecg \in F^n} 1_{B_{\vecg}} \cdot \omega^*(f \vert_{\sigma^{-1}(\vecg)})$, which lies in $\varphi^{*-1} L$ by Lemma~\ref{Lem:phi^-1}~d). This formula also shows continuity in the topological setting.
\eproof

Our next goal is to show that coarse maps which are close induce the same chain and cochain maps up to homotopy. Let $\varphi, \, \phi: \: G \to H$ be two coarse embeddings with $\varphi \sim \phi$. Let $L$ be a $\res$-invariant $RG$-submodule of $C(G,W)$ and $M$ a $\res$-invariant $RH$-submodule of $C(H,W)$. Let $\cG = G \rtimes G$ and $\cH = H \rtimes H$. Define $\theta: \: G \to \cH, \, x \ma (\varphi(x),\varphi(x) \phi(x)^{-1})$. For $n \geq 0$ and $1 \leq h \leq n+1$, let $\kappa_n^{(h)}: \: \cG^{(n)} \to \cH^{(n+1)}$ be given by $\kappa_0^{(1)} = \theta$, and for $n \geq 1$,
\begin{align*}
  & \kappa_n^{(h)}(\gamma_1, \dotsc, \gamma_n) = (\varphi^1(\gamma_1), \dotsc, \varphi^1(\gamma_{h-1}), \theta(r(\gamma_h)), \phi^1(\gamma_h), \dotsc, \phi^1(\gamma_n)) \ {\rm for} \ 1 \leq h \leq n,\\
  & \kappa_n^{(n+1)}(\gamma_1, \dotsc, \gamma_n) = (\phi^1(\gamma_1), \dotsc, \phi^1(\gamma_n), \theta(s(\gamma_n))).
\end{align*}
Moreover, for $n \geq 1$ and $1 \leq h \leq n$, let $\kappa^n_{(h)}: \: \cG^{(n-1)} \to \cH^{(n)}$ be given by $\kappa^1_{(1)} = \theta$, and for $n \geq 2$,
\begin{align*}
  & \kappa^n_{(h)}(\gamma_1, \dotsc, \gamma_{n-1}) = (\varphi^1(\gamma_1), \dotsc, \varphi^1(\gamma_{h-1}), \theta(r(\gamma_h)), \phi^1(\gamma_h), \dotsc, \phi^1(\gamma_{n-1})) \ {\rm for} \ 1 \leq h \leq n-1,\\
  & \kappa^n_{(n)}(\gamma_1, \dotsc, \gamma_{n-1}) = (\phi^1(\gamma_1), \dotsc, \phi^1(\gamma_{n-1}), \theta(s(\gamma_{n-1}))).
\end{align*}
\blemma
\label{Lem:simh}
(i) $k_n^{(h)} = (\kappa_n^{(h)})_*: \: D_n(L) \to D_{n+1}(\varphi_* L) = D_{n+1}(\phi_* L)$ is well-defined for all $n$ and $h$. $k_n \defeq \sum_{h=1}^{n+1} (-1)^{h+1} k_n^{(h)}$ gives a chain homotopy $D_*(\varphi) \sim_h D_*(\phi)$.
\setlength{\parindent}{0.5cm} \setlength{\parskip}{0cm}

(ii) $k^n_{(h)} = (\kappa^n_{(h)})^*: \: D^n(M) \to D_{n-1}(\varphi^* M) = D_{n-1}(\phi^* M)$ is well-defined for all $n$, $h$. $k^n \defeq \sum_{h=1}^{n} (-1)^{h+1} k^n_{(h)}$ gives a cochain homotopy $D^*(\varphi) \sim_h D^*(\phi)$.
\elemma
\setlength{\parindent}{0cm} \setlength{\parskip}{0cm}

\bproof
(i) Let us show that $k_n^{(h)}$ is well-defined, i.e., $(\kappa_n^{(h)})_*(f) \in D_{n+1}(\varphi_* L)$ for all $f \in D_n(L)$. We may assume $\supp(f) = \gekl{\vecg}$ for a single $\vecg = (g_1, \dotsc, g_n) \in G^n$, as a general element in $D_n(L)$ is a finite sum of such $f$. We first show that $\supp((\kappa_n^{(h)})_*(f))$ is finite. By \eqref{D_:F:finite} and because $\varphi \sim \phi$, we know that
$$
  F \defeq \menge{\varphi(x) \varphi(g_i^{-1}x)^{-1}}{x \in G, \, 1 \leq i \leq n}
  \cup \menge{\varphi(x) \phi(x)^{-1}}{x \in G}
  \cup \menge{\phi(x) \phi(g_i^{-1}x)^{-1}}{x \in G, \, 1 \leq i \leq n}
$$
is finite. As $\kappa_n^{(h)}(x,\vecg)$ lies in $\sigma^{-1}(F^{n+1})$ for all $x \in G$, we conclude that $\supp((\kappa_n^{(h)})_*(f))$ is contained in $F^{n+1}$, which is finite. Let us show that for every $\vech = (h_1, \dotsc, h_{n+1}) \in H^{n+1}$, $(\kappa_n^{(h)})_*(f) \vert_{\sigma^{-1}(\vech)}$ lies in $\varphi_* M$. Define
\begin{align*}
  A \defeq &\left\{x \in G : \: \varphi(g_{i-1}^{-1} \dotsm g_1^{-1}x) \varphi(g_i^{-1} \dotsm g_1^{-1}x)^{-1} = h_i \ {\rm for} \ {\rm all} \ 1 \leq i \leq h-1, \, \right.\\
  &\varphi(g_{h-1}^{-1} \dotsm g_1^{-1}x) \phi(g_{h-1}^{-1} \dotsm g_1^{-1}x)^{-1} = h_h, \,\\
  &\left. \phi(g_{i-1}^{-1} \dotsm g_1^{-1}x) \phi(g_i^{-1} \dotsm g_1^{-1}x)^{-1} = h_{i+1} \ {\rm for} \ {\rm all} \ h \leq i \leq n \right\}.
\end{align*}
Then $\kappa_n^{(h)}(x,\vecg) \in \sigma^{-1}(\vech)$ if and only if $x \in A$. Hence $(\kappa_n^{(h)})_*(f) \vert_{\sigma^{-1}(\vech)} = \varphi_* \rukl{1_A \cdot \rukl{f \vert_{\sigma^{-1}(\vecg)}}}$. As $f \vert_{\sigma^{-1}(\vecg)}$ lies in $L$, $L$ is $\res$-invariant, and $\varphi_*(\ti{f}) \in \varphi_* L$ for all $\ti{f} \in L$, we see that $(\kappa_n^{(h)})_*(f) \vert_{\sigma^{-1}(\vech)} \in \varphi_* L$. Hence $k_n^{(h)}$ is well-defined for all $n$ and $h$. A straightforward computation shows that $k_n$ indeed gives us the desired chain homotopy.
\setlength{\parindent}{0.5cm} \setlength{\parskip}{0cm}

(ii) Let us show that $k^n_{(h)}$ is well-defined, i.e., $(\kappa^n_{(h)})^*(f) \vert_{\sigma^{-1}(\vecg)} \in \varphi^* M$ for all $\vecg = (g_1, \dotsc, g_{n-1}) \in G^{n-1}$ and $f \in D^n(M)$. As in the proof of (i), note that 
\begin{align*}
  F \defeq &\menge{\varphi(x) \varphi(g_i^{-1}x)^{-1}}{x \in G, \, 1 \leq i \leq n-1}\\
  \cup &\menge{\varphi(x) \phi(x)^{-1}}{x \in G}\\
  \cup &\menge{\phi(x) \phi(g_i^{-1}x)^{-1}}{x \in G, \, 1 \leq i \leq n-1}
\end{align*}
is finite, and that $\kappa^n_{(h)}(x,\vecg) \in \sigma^{-1}(F^n)$. For $\vech = (h_1, \dotsc, h_n) \in F^n$, set
\begin{align*}
  A_{\vech} \defeq &\left\{x \in G : \: \varphi(g_{i-1}^{-1} \dotsm g_1^{-1}x) \varphi(g_i^{-1} \dotsm g_1^{-1}x)^{-1} = h_i \ {\rm for} \ {\rm all} \ 1 \leq i \leq h-1, \, \right.\\
  &\varphi(g_{h-1}^{-1} \dotsm g_1^{-1}x) \phi(g_{h-1}^{-1} \dotsm g_1^{-1}x)^{-1} = h_h, \,\\
  &\left. \phi(g_{i-1}^{-1} \dotsm g_1^{-1}x) \phi(g_i^{-1} \dotsm g_1^{-1}x)^{-1} = h_{i+1} \ {\rm for} \ {\rm all} \ h \leq i \leq n-1 \right\}.
\end{align*}
Then $G = \bigsqcup_{\vech \in F^n} A_{\vech}$, and for $x \in A_{\vech}$, $\kappa^n_{(h)}(x,\vecg) = (\varphi(x),\vech)$. Hence
$$
  (\kappa^n_{(h)})^*(f) \vert_{\sigma^{-1}(\vecg)} (x) = f(\kappa^n_{(h)}(x,\vecg)) = \sum_{\vech \in F^n} 1_{A_{\vech}}(x) \cdot \rukl{f \vert_{\sigma^{-1}(\vech)}}(\varphi(x))
$$
and thus $(\kappa^n_{(h)})^*(f) \vert_{\sigma^{-1}(\vecg)} = \sum_{\vech \in F^n} 1_{A_{\vech}} \cdot \rukl{\varphi^* \rukl{f \vert_{\sigma^{-1}(\vech)}}}$. Since $f \vert_{\sigma^{-1}(\vech)} \in M$, $\varphi^*(\ti{f}) \in \varphi^* M$ for all $\ti{f} \in M$ and $\varphi^* M$ is $\res$-invariant, this shows that $(\kappa^n_{(h)})^*(f) \vert_{\sigma^{-1}(\vecg)} \in \varphi^* M$. Hence $k^n_{(h)}$ is well-defined. It is straightforward to check that $k^n$ indeed gives us the desired cochain homotopy.
\eproof

Now let $\varphi: \: G \to H$ be a coarse embedding, $\omega: \: H \to G$ as above and $L$ an $\res$-invariant $RG$-submodule of $C(G,W)$. Define $\vartheta: \: H \to \cH, \, y \ma (y, y (\varphi \circ \omega)(y)^{-1})$. For $n \geq 0$ and $1 \leq h \leq n+1$, let $\lambda_n^{(h)}: \: \cH^{(n)} \to \cH^{(n+1)}$ be given by $\lambda_0^{(1)} = \vartheta$, and for $n \geq 1$,
\begin{align*}
  & \lambda_n^{(h)}(\eta_1, \dotsc, \eta_n) = (\eta_1, \dotsc, \eta_{h-1}, \vartheta(r(\eta_h)), (\varphi \circ \omega)^1(\eta_h), \dotsc, (\varphi \circ \omega)^1(\eta_n)) \ {\rm for} \ 1 \leq h \leq n,\\
  & \lambda_n^{(n+1)}(\eta_1, \dotsc, \eta_n) = (\eta_1, \dotsc, \eta_n, \vartheta(s(\eta_n))).
\end{align*}
Moreover, for $n \geq 1$ and $1 \leq h \leq n$, let $\lambda^n_{(h)}: \: \cH^{(n-1)} \to \cH^{(n)}$ be given by $\lambda^1_{(1)} = \vartheta$, and for $n \geq 2$,
\begin{align*}
  & \lambda^n_{(h)}(\eta_1, \dotsc, \eta_{n-1}) = (\eta_1, \dotsc, \eta_{h-1}, \vartheta(r(\eta_h)), (\varphi \circ \omega)^1(\eta_h), \dotsc, (\varphi \circ \omega)^1(\eta_{n-1})) \ {\rm for} \ 1 \leq h \leq n-1,\\
  & \lambda^n_{(n)}(\eta_1, \dotsc, \eta_{n-1}) = (\eta_1, \dotsc, \eta_{n-1}, \vartheta(s(\eta_{n-1}))).
\end{align*}
\blemma
\label{Lem:phiomega=id}
(i) We have $D_*(\omega \circ \varphi) \sim_h \id$. $l_n^{(h)} = (\lambda_n^{(h)})_*: \: D_n(\varphi_* L) \to D_n(\varphi_* L)$ is well-defined for all $n$ and $h$. $l_n \defeq \sum_{h=1}^{n+1} (-1)^{h+1} l_n^{(h)}$ gives a chain homotopy $D_*(\varphi \circ \omega) \sim_h \id$.
\setlength{\parindent}{0.5cm} \setlength{\parskip}{0cm}

(ii) We have $D^*(\omega \circ \varphi) \sim_h \id$. $l^n_{(h)} = (\lambda^n_{(h)})_*: \: D^n(\varphi^{*-1} L) \to D_n(\varphi^{*-1}  L)$ is well-defined for all $n$ and $h$. $l^n \defeq \sum_{h=1}^n (-1)^{h+1} l^n_{(h)}$ gives a chain homotopy $D^*(\varphi \circ \omega) \sim_h \id$.
\elemma
\bproof
(i) $D_*(\omega \circ \varphi) \sim_h \id$ follows from Lemma~\ref{Lem:simh}~(i) and \eqref{omegaphi}. That $l_n^{(h)}$ is well-defined can be proven as Lemma~\ref{Lem:D_}~(ii). It is straightforward to check that $l_n$ gives the desired chain homotopy.
\setlength{\parindent}{0.5cm} \setlength{\parskip}{0cm}

(ii) $D^*(\omega \circ \varphi) \sim_h \id$ follows from Lemma~\ref{Lem:simh}~(ii) and \eqref{omegaphi}. That $l^n_{(h)}$ is well-defined can be proven as Lemma~\ref{Lem:D^}~(ii). It is straightforward to check that $l^n$ gives the desired cochain homotopy.
\eproof

Combining Corollary~\ref{Cor:G-cG} and Lemmas~\ref{Lem:D_}, \ref{Lem:D^}, \ref{Lem:simh} and \ref{Lem:phiomega=id}, we obtain
\btheo
\label{H_,H^}
Let $\varphi: \: G \to H$ be a coarse map, $L$ a $\res$-invariant $RG$-submodule of $C(G,W)$ and $M$ a $\res$-invariant $RH$-submodule of $C(G,W)$.
\setlength{\parindent}{0.5cm} \setlength{\parskip}{0cm}

(i) $D_*(\varphi)$ induces homomorphisms $H_*(\varphi): \: H_*(G,L) \to H_*(H,\varphi_* L)$. If $\varphi$ is a coarse embedding, $H_*(\varphi)$ is an isomorphism with inverse $H_*(\omega)$. If in addition $L$ is a topological $\res$-invariant module, $D_*(\varphi)$ also induces topological isomorphisms ${\bar H}_*(\varphi): \: {\bar H}_*(G,L) \cong {\bar H}_*(H,\varphi_* L)$.

If $\phi: \: G \to H$ is a coarse map with $\varphi \sim \phi$, then $H_*(\varphi) = H_*(\phi)$, and ${\bar H}_*(\varphi) = {\bar H}_*(\phi)$ in the topological case.

If $\psi: \: H \to K$ is another coarse map, then $H_*(\psi \circ \varphi) = H_*(\psi) \circ H_*(\varphi)$, and ${\bar H}_*(\psi \circ \varphi) = {\bar H}_*(\psi) \circ {\bar H}_*(\varphi)$ in the topological case.

(ii) $D^*(\varphi)$ induces homomorphisms $H^*(\varphi): \: H^*(H,M) \to H^*(G,\varphi^* M)$. If $\varphi$ a coarse embedding, then $H^*(\varphi): \: H^*(H,\varphi^{*-1} L) \to H^*(G,L)$ is an isomorphism with inverse $H^*(\omega)$. If in addition $L$ and $M$ are topological $\res$-invariant modules, then $D^*(\varphi)$ also induces continuous homomorphisms ${\bar H}^*(\varphi): \: {\bar H}^*(H,M) \to {\bar H}^*(G,\varphi^* M)$ and topological isomorphisms ${\bar H}^*(\varphi): \: {\bar H}^*(H,\varphi^{*-1} L) \to {\bar H}^*(G,L)$.

If $\phi: \: G \to H$ is a coarse map with $\varphi \sim \phi$, then $H^*(\varphi) = H^*(\phi)$, and ${\bar H}^*(\varphi) = {\bar H}^*(\phi)$ in the topological case.

If $\psi: \: H \to K$ is another coarse map, then $H^*(\psi \circ \varphi) = H^*(\varphi) \circ H^*(\psi)$, and ${\bar H}^*(\psi \circ \varphi) = {\bar H}^*(\varphi) \circ {\bar H}^*(\psi)$ in the topological case.
\etheo

In particular, for coarse equivalences, i.e., coarse embeddings which are invertible modulo $\sim$, we get
\bcor
\label{Cor:phi--iso}
If $\varphi: \: G \to H$ is a coarse equivalence, then we obtain isomorphisms
$$
  H_*(\varphi): \: H_*(G,L) \cong H_*(H,\varphi_* L), \ \ \ H^*(\varphi): \: H^*(H,M) \cong H^*(G,\varphi^* M),
$$
and, in the topological case,
${\bar H}_*(\varphi): \: {\bar H}_*(G,L) \cong {\bar H}_*(H,\varphi_* L), \ \ \ {\bar H}^*(\varphi): \: {\bar H}^*(H,M) \cong {\bar H}^*(G,\varphi^* M)$.
\ecor

\bremark
\label{Rem:functor1}
Our constructions are functorial in $W$: Let $L_1 \subseteq C(G,W_1)$ and $L_2 \subseteq C(G,W_2)$ be $\res$-invariant $RG$-submodules, and assume that an $R$-linear map $\omega: \: W_1 \to W_2$ induces an $RG$-linear map $\lambda: \: L_1 \to L_2$. Then we also obtain an induced map $\varphi_* \lambda: \: \varphi_* L_1 \to \varphi_* L_2$, and we get commutative diagrams
\begin{align*}
  \xymatrix{
  D_*(L_1) \ar[d]^{D_*(\lambda)} \ar[r]^{D_*(\varphi)} & D_*(\varphi_* L_1) \ar[d]^{D_*(\varphi_* \lambda)}
  \\
  D_*(L_2) \ar[r]^{D_*(\varphi)} & D_*(\varphi_* L_2)
  }
  \ \ \ \ \ \ 
  \xymatrix{
  H_*(G, L_1) \ar[d]^{H_*(\lambda)} \ar[r]^{H_*(\varphi)} & H_*(H, \varphi_* L_1) \ar[d]^{H_*(\varphi_* \lambda)}
  \\
  H_*(G, L_2) \ar[r]^{H_*(\varphi)} & H_*(H, \varphi_* L_2)
  }  
\end{align*}
A similar statement applies to reduced homology in the topological setting, and to (reduced) cohomology.
\eremark

\subsection{Consequences}
\label{ss:Consequences}
Let us apply our results to the Examples in \ref{Examples}. Corollary~\ref{RGlpc0} (i) (c) below generalizes the result in \cite{Ger} that $H^n(G,RG)$ is a coarse invariant for groups with property $F_n$. The reader may also consult \cite[Example~5.21]{Roe}. Corollary~\ref{RGlpc0} (ii) (1) was known in special cases. For instance, in \cite{El}, group cohomology with $\ell^p$ coefficients has been identified with nonreduced $L^p$-cohomology, and in \cite{Pul,BMV,MV}, reduced group cohomology in degree $1$ (i.e., ${\bar H}^1$) with $\ell^p$ coefficients has been identified with $L^p$-cohomology, as studied in \cite{Gro,Pan}. Since $L^p$-cohomology is known to be a coarse invariant, this gives the special case of (ii)~(1) where $p \in [1,\infty[$ and our groups are finitely generated. Also, the case $p = \infty$ in (ii)~(1) was known since $H_*(G,\ell^{\infty} G)$ can be identified with uniformly finite homology (see \cite{BD,BNW}).
\bcor
\label{RGlpc0}
Let $G$ and $H$ be countable discrete groups and $\varphi: \: G \to H$ a coarse equivalence. 
\setlength{\parindent}{0.5cm} \setlength{\parskip}{0cm}

(i) For every commutative ring $R$ with unit and every $R$-module $W$, $\varphi$ induces isomorphisms
\begin{enumerate}
\item[(a)] $H_*(G,C(G,W)) \cong H_*(H,C(H,W))$,
\item[(b)] $H_*(G,C_f(G,W)) \cong H_*(H,C_f(H,W))$ and $H^*(G,C_f(G,W)) \cong H^*(H,C_f(H,W))$,
\item[(c)] $H^*(H,RH \otimes_R W) \cong H^*(G,RG \otimes_R W)$.
\end{enumerate}

(ii) Let $R = \Rz$ or $R = \Cz$ and $W=R$.
\begin{enumerate}
\item For all $0<p \leq \infty$, $\varphi$ induces isomorphisms
\begin{align*}
  &H_*(G,\ell^p (G,W)) \cong H_*(H,\ell^p (H,W)), \ \ \ H^*(H,\ell^p (H,W)) \cong H^*(G,\ell^p (G,W)),\\
  &{\bar H}_*(G,\ell^p (G,W)) \cong {\bar H}_*(H,\ell^p (H,W)), \ \ \ {\bar H}^*(H,\ell^p (H,W)) \cong {\bar H}^*(G,\ell^p (G,W)),
\end{align*}
\item $\varphi$ induces isomorphisms
\begin{align*}
  &H_*(G,c_0 (G,W)) \cong H_*(H,c_0 (H,W)), \ \ \ H^*(H,c_0 (H,W)) \cong H^*(G,c_0 (G,W)),\\
  &{\bar H}_*(G,c_0 (G,W)) \cong {\bar H}_*(H,c_0 (H,W)), \ \ \ {\bar H}^*(H,c_0 (H,W)) \cong {\bar H}^*(G,c_0 (G,W)).
\end{align*}
\item Let $G$ and $H$ be a finitely generated discrete groups. Then, for all $s \in \Rz \cup \gekl{\infty}$ and $1 \leq p \leq \infty$, $\varphi$ induces isomorphisms
\begin{align*}
  &H_*(G,H^{s,p}(G,W)) \cong H_*(H,H^{s,p}(H,W)), \ \ \  H^*(H,H^{s,p}(H,W)) \cong H^*(G,H^{s,p}(G,W)),\\
  &{\bar H}_*(G,H^{s,p}(G,W)) \cong {\bar H}_*(H,H^{s,p}(H,W)), \ \ \ {\bar H}^*(H,H^{s,p}(H,W)) \cong {\bar H}^*(G,H^{s,p}(G,W)).
\end{align*}
\end{enumerate}
\ecor
\setlength{\parindent}{0cm} \setlength{\parskip}{0cm}

\bproof
The point is that $L(G) = C(G,W)$, $C_f(G,W)$, $RG \otimes_R W$, $\ell^p(G,W)$, $c_0(G,W)$ or $H^{s,p}(G,W)$ have the property that for every coarse equivalence $\varphi: \: G \to H$, we have $\varphi_* L(G) = L(H)$ (and also topologically in the topological setting). Our claim now follows from Corollary~\ref{Cor:phi--iso}.
\eproof

As an immediate consequence, we obtain a new proof of the result in \cite{Sau} that homological and cohomological dimensions over $R$ are preserved by coarse embeddings as long as these dimensions are finite.
\bcor
\label{Cor:dim-finite}
Let $R$ be a commutative ring with unit. Let $G$ and $H$ be countable discrete groups, and assume that there is a coarse embedding $\varphi: \: G \to H$.
\setlength{\parindent}{0.5cm} \setlength{\parskip}{0cm}

If $G$ has finite homological dimension over $R$, i.e., ${\rm hd}_R \, G < \infty$, then ${\rm hd}_R \, G \leq {\rm hd}_R \, H$.

If $G$ has finite cohomological dimension over $R$, i.e., ${\rm cd}_R \, G < \infty$, then ${\rm cd}_R \, G \leq {\rm cd}_R \, H$.
\ecor
\setlength{\parindent}{0cm} \setlength{\parskip}{0cm}

\bproof
Assume that ${\rm hd}_R \, G = n < \infty$. Let $W$ be an $RG$-module such that $H_n(G,W) \ncong \gekl{0}$. Define $W \into C(G,W), \, w \ma f_w$, where $f_w (x) = x^{-1}.w$. It is easy to see that this is an embedding of $RG$-modules when we view $W$ as an $R$-module to construct $C(G,W)$ (i.e., we define the $RG$-module structure by setting $(g.f) (x) = f(g^{-1}.x)$ for $f \in C(G,W)$). The long exact sequence in homology gives us $0 \to H_n(G,W) \to H_n(G,C(G,W)) \to \dotso$
because the $(n+1)$-th group homology of $G$ vanishes for all coefficients by assumption. Hence $H_n(G,C(G,W)) \ncong \gekl{0}$. By Theorem~\ref{H_,H^}~(i), we have $H_n(H,\varphi_* C(G,W)) \cong H_n(G,C(G,W)) \ncong \gekl{0}$. Thus ${\rm hd}_R \, H \geq n$.
\setlength{\parindent}{0.5cm} \setlength{\parskip}{0cm}

Now assume ${\rm cd}_R \, G = n < \infty$. By \cite[Proposition~(2.3)]{Bro}, we know that $H^n(G,RG \otimes_R W) \ncong \gekl{0}$ for some $R$-module $W$. By Theorem~\ref{H_,H^}~(ii), $H^n(H,\varphi^{*-1} (RG \otimes_R W)) \cong H^n(G,RG \otimes_R W) \ncong \gekl{0}$. Thus ${\rm cd}_R \, H \geq n$.
\eproof

We also obtain a new proof for the following result, first proven in \cite{Sau}:
\bcor
\label{Cor:cdQ}
Let $R$, $G$ and $H$ be as above. Assume that $\varphi: \: G \to H$ is a coarse embedding. If $G$ is amenable and $\Qz \subseteq R$, then ${\rm hd}_R \, G \leq {\rm hd}_R \, H$ and ${\rm cd}_R \, G \leq {\rm cd}_R \, H$. 
\ecor
\bproof
As explained in \cite[\S~4]{Sau}, it was observed in \cite{Sha} that our assumptions on $G$ and $R$ imply the existence of an $RG$-linear split $C_f(G,R) \to R$ for the canonical homomorphism $R \to C_f(G,R)$ embedding $R$ as constant functions. Hence, given an arbitrary $RG$-module $V$, we obtain by tensoring with $V$ over $R$ that the canonical homomorphism $V \to C_f(G,V)$ splits. Note that $G$ acts on $C_f(G,V)$ diagonally, so that $C_f(G,V)$ is not a $\res$-invariant module in our sense. But $C_f(G,V) \cong C_f(G,V_{\rm triv})$, where $V_{\rm triv}$ is the $R$-module $V$ viewed as a $RG$-module with trivial $G$-action. Hence ${\rm hd}_R \, G = \sup_n \menge{n}{H_n(G,C_f(G,W)) \ncong \gekl{0} \ {\rm for} \ {\rm some} \ R \text{-} {\rm module} \ W}$. As $H_n(H,\varphi_* C_f(G,W)) \cong H_n(G,W)$ by Theorem~\ref{H_,H^}~(i), we conclude that ${\rm hd}_R \, G \leq {\rm hd}_R \, H$. The proof for ${\rm cd}_R$ is analogous.
\eproof

At this point, the following interesting question arises naturally:
\bquestion
\label{Q:dim}
Let $R$ be a commutative ring with unit, $G$ and $H$ countable discrete groups with no $R$-torsion. If $G$ and $H$ are coarsely equivalent, do we always have ${\rm hd}_R \, G = {\rm hd}_R \, H$ and ${\rm cd}_R \, G = {\rm cd}_R \, H$?
\equestion
Having no $R$-torsion means that orders of finite subgroups must be invertible in $R$, and this is certainly a hypothesis we have to include. For instance, \cite[Theorem~1.4]{Pet} implies that the answer to Question~\ref{Q:dim} is affirmative if our groups lie in the class $H \cF$. This class $H \cF$ has been introduced by Kropholler in \cite{Kro} and is defined as the smallest class of groups containing all finite groups and every group $G$ which acts cellularly on a finite dimensional contractible CW-complex with all isotropy subgroups already in $H \cF$. All countable elementary amenable groups and all countable linear groups lie in $H \cF$, and it is closed under subgroups, extensions, and countable direct unions. 
\bcor[to Theorem~1.4 in \cite{Petr}]
\label{Cor:dim=}
If $G$ and $H$ are in $H \cF$, then the answer to Question~\ref{Q:dim} is affirmative.
\ecor
\bproof
\cite[Theorem~1.4]{Petr} implies that
\bgl
\label{cd=supcdCE}
{\rm cd}_R \, G = \sup \menge{{\rm cd}_R \, G'}{G' \ {\rm coarsely} \ {\rm embeds} \ {\rm into} \ G \ {\rm and} \ {\rm cd}_R \, G' < \infty}.
\egl
Similarly for $H$. Now Corollary~\ref{Cor:dim-finite} implies ${\rm cd}_R \, G = {\rm cd}_R \, H$. Equality for ${\rm hd}_R$ follows because for countable groups, ${\rm cd}_R$ is infinite if and only if ${\rm hd}_R$ is infinite by \cite[Theorem~4.6]{Bie}.
\eproof

\bremark
The proof of Corollary~\ref{Cor:dim=} shows that Question~\ref{Q:dim} has an affirmative answer among all groups satisfying \eqref{cd=supcdCE}. In particular, for groups satisfying \cite[Conjecture~1.6]{Petr}, Question~\ref{Q:dim} has an affirmative answer. While counterexamples to \cite[Conjecture~1.6]{Petr} are presented in \cite{Gan}, these examples still satisfy \eqref{cd=supcdCE}, as becomes clear in \cite{Gan}. Hence also for them, Question~\ref{Q:dim} has an affirmative answer.
\eremark

Let us now show that being of type $FP_n$ over a ring $R$ is a coarse invariant. An alternative approach, based on \cite{KK}, has been sketched in \cite[Theorem~9.61]{DK}. The case $R = \Zz$ is treated in \cite{Alo}. Recall that for a commutative ring $R$ with unit, a group $G$ is of type $FP_n$ over $R$ if the trivial $RG$-module $R$ has a projective resolution $\dotso \to P_1 \to P_0 \to R \to 0$ where $P_i$ is finitely generated for all $0 \leq i \leq n$.
\bcor
\label{Cor:FPn}
Let $G$ and $H$ be two countable discrete groups. Assume that $G$ and $H$ are coarsely equivalent. Then $G$ is a of type $FP_n$ over $R$ if and only if $H$ is of type $FP_n$ over $R$.
\ecor
\bproof
By \cite[Proposition~2.3]{Bie}, $G$ is of type $FP_n$ over $R$ if and only if $G$ is finitely generated and $H_k(G,\prod_I RG) \cong \gekl{0}$ for all $1 \leq k \leq n$, where $I$ is an index set with $\abs{I} = \max(\aleph_0,\abs{R})$. The map $\prod_I RG \to C(G,\prod_I R), \, (f_i)_i \ma f$, where $(f(x))_i = f_i(x)$, identifies $\prod_I RG$ with the $RG$-submodule $L(G)$ of $C(G,\prod_I R)$ consisting of those functions $f$ with the property that for every $i \in I$, $(f(x))_i = 0$ for all but finitely many $x \in G$. Clearly, $L(G)$ is $\res$-invariant. Denote the analogous $\res$-invariant $RH$-submodule of $C(H,\prod_I R)$ by $L(H)$. It is easy to see that given a coarse equivalence $\varphi: \: G \to H$, we have $\varphi_* L(G) = L(H)$. Hence, by Theorem~\ref{H_,H^}~(i), we have $H_k(G,\prod_I RG) \cong H_k(G,L(G)) \cong H_k(H,L(H)) \cong H_k(H,\prod_I RH)$.
\eproof

As another consequence, we generalize the result in \cite{Ger} that for groups of type $F_{\infty}$, being a (Poincar{\'e}) duality group over $\Zz$ is a coarse invariant. We obtain an improvement since we can work over arbitrary rings $R$ and do not need the $F_{\infty}$ assumption. We only need our groups to have finite cohomological dimension over $R$. Recall that a group $G$ is called a duality group over $R$ if there is a right $RG$-module $C$ and an integer $n \geq 0$ with natural isomorphisms $H^k(G,A) \cong H_{n-k}(G,C \otimes_R A)$ for all $k \in \Zz$ and all $RG$-modules $A$ (see \cite[\S~9.2]{Bie}, \cite{Bie75}, and \cite[Chapter~VIII, \S~10]{Bro}). $G$ is called a Poincar{\'e} duality group over $R$ if $C \cong R$ as $R$-modules. The class of duality groups is closed under extensions and under taking graphs of groups, with certain hypotheses (see \cite{Bie,DT}). Examples of groups which are not duality groups over $\Zz$ but over some other ring can be found in \cite{Dav}, and examples of (Poincar{\'e}) duality groups which are not of type $F_{\infty}$ appear in \cite{Dav,Lea}. The second part of the following corollary generalizes \cite[Theorem~3.3.2]{Sha}. 
\bcor
\label{Cor:duality}
Let $R$ be a commutative ring with unit. Let $G$ and $H$ be countable discrete groups with finite cohomological dimension over $R$. If $G$ and $H$ are coarsely equivalent, then $G$ is a (Poincar{\'e}) duality group over $R$ if and only if $H$ is a (Poincar{\'e}) duality group over $R$.
\setlength{\parindent}{0.5cm} \setlength{\parskip}{0cm}

If $G$ and $H$ are amenable and $\Qz \subseteq R$, then $G$ is a (Poincar{\'e}) duality group over $R$ if and only if $H$ is a (Poincar{\'e}) duality group over $R$.
\ecor
\bproof
By \cite[Theorem~5.5.1 and Remark~5.5.2]{Bie75}, we know that a group $G$ is a duality group if and only if it has finite cohomological dimension, there is $n$ such that $H^k(G,A) \cong \gekl{0}$ for all $k \neq n$ and all induced $RG$-modules $A$, and $G$ is of type $FP_n$ over $R$. The second property is a coarse invariant by Corollary~\ref{RGlpc0}~(i)~(c). The third property is a coarse invariant by Corollary~\ref{Cor:FPn}.
Hence being a duality group is a coarse invariant. Being a Poincar{\'e} duality group means being a duality group and having dualizing module isomorphic to $R$. By Corollary~\ref{RGlpc0}~(i)~(c), the dualizing module is a coarse invariant. Thus being a Poincar{\'e} duality group is also a coarse invariant. The second part follows from the first part of the corollary and Corollary~\ref{Cor:cdQ}.
\eproof
If Question~\ref{Q:dim} has an affirmative answer, then we can replace the assumption of finite cohomological dimension by having no $R$-torsion in the first part of Corollary~\ref{Cor:duality}.
\setlength{\parindent}{0cm} \setlength{\parskip}{0.5cm}

As another consequence, we obtain the following rigidity result for coarse embeddings into Poincar{\'e} duality groups. The proof follows the one of \cite[Proposition~9.22]{Bie}.
\setlength{\parindent}{0cm} \setlength{\parskip}{0cm}
\bcor
\label{Cor:dim<dim}
Let $G$ and $H$ be countable discrete groups. Let $H$ be a Poincar{\'e} duality group over a commutative ring $R$ with unit. Assume that there is a coarse embedding $\varphi: \: G \to H$ which is not a coarse equivalence. If ${\rm hd}_R \, G < \infty$, then ${\rm hd}_R \, G < {\rm cd}_R \, H$. If, in addition, $G$ is of type $FP_{\infty}$ (i.e, $FP_n$ for all $n$), then ${\rm cd}_R \, G < {\rm cd}_R \, H$.
\setlength{\parindent}{0.5cm} \setlength{\parskip}{0cm}

In particular, every self coarse embedding of a Poincar{\'e} duality group over $R$ must be a coarse equivalence.
\ecor
\bproof
Let $n = {\rm cd}_R \, H$. Let $D = H^n(R,RH)$. As $H$ is a  Poincar{\'e} duality group over $R$, $D \cong R$ as $R$-modules, and the $RH$-module structure of $D$ is given by a group homomorphism $H \to R^*, \, h \ma u_h$. We know that ${\rm hd}_R \, G \leq {\rm cd}_R \, G \leq n$ by \cite[Theorem~4.6]{Bie} and Corollary~\ref{Cor:dim-finite}. Now let $L$ be a $\res$-invariant $RG$-submodule of $C(G,W)$. Then, by Theorem~\ref{H_,H^}~(i), $H_n(G,L) \cong H_n(H,\varphi_* L) \cong H^0(H,\Hom_R(D,\varphi_* L)) \cong \rukl{\Hom_R(D,\varphi_* L)}^H$, where we used that $H$ is a Poincar{\'e} duality group over $R$. Clearly, $\Hom_R(D,\varphi_* L) \cong \varphi_* L$ as $R$-modules, and the $H$-action of $\Hom_R(D,\varphi_* L)$ becomes $h _{\bullet} f = u_h \cdot (h.f)$ for $f \in \varphi_* L$. Now take $f \in (\varphi_* L)^H$. If $f \neq 0$, then $f(y) \neq 0$ for some $y \in H$, and it follows from $h _{\bullet} f = f$ for all $h \in H$ that $f(y) \neq 0$ for all $y \in H$. This, however, contradicts Lemma~\ref{Lem:phi_L=bigoplus} as $H$ cannot be contained in a finite union of $h_j Y_j$s if $\varphi$ is not a coarse equivalence. Hence $H_n(G,L) \cong (\varphi_* L)^H \cong \gekl{0}$. This implies ${\rm hd}_R \, G < n$ (compare also the proof of Corollary~\ref{Cor:dim-finite}). The rest follows from \cite[Theorem~4.6~(c)]{Bie} and that Poincar{\'e} duality groups are of type $FP_{\infty}$.
\eproof
\bquestion
In Corollary~\ref{Cor:dim<dim}, do we always get ${\rm cd}_R \, G < {\rm cd}_R \, H$, even without the $FP_{\infty}$ assumption? In other words, is the analogue of the main theorem in \cite{Str} true for coarse embeddings?
\equestion

We present one more application: Vanishing of $\ell^2$-Betti numbers is a coarse invariant. This was shown in \cite{Pan} for groups of type $F_{\infty}$, for more general groups in \cite{Ogu} (as explained in \cite{SauSch}), and for all countable discrete groups in \cite[Corollary~6.3]{MOSS}. Recently, Sauer and Schr{\"o}dl were even able to cover all unimodular locally compact second countable groups \cite{SauSch}. As vanishing of the $n$-th $\ell^2$-Betti number is equivalent to ${\bar H}^n(G,\ell^2 G) \cong \gekl{0}$ by \cite[Proposition~3.8]{Pet}, Corollary~\ref{RGlpc0}~(ii)~(1) gives another approach to the aforementioned result.
\bcor
Let $G$ and $H$ be countable discrete groups which are coarsely equivalent. Then, for all $n$, the $n$-th $\ell^2$-Betti number of $G$ vanishes if and only if the $n$-th $\ell^2$-Betti number of $H$ vanishes.
\ecor

\end{document}